\documentclass[leqno, myheadings, twoside]{amsart}
\usepackage{amsfonts}
\usepackage{amsmath, amsthm, amssymb, amscd, amsxtra,graphicx}
\usepackage{latexsym, amsfonts}
\usepackage{url}
\usepackage{texdraw}
\usepackage{epsfig}
\makeatletter \@addtoreset{equation}{section}

\makeatletter \renewcommand{\@biblabel}[1]{#1.}

\theoremstyle{remark}

\begin{document}
\title [One can hear the area and curvature] {One can hear the area and curvature of
boundary of a domain by hearing the Steklov eigenvalues}
\author{Genqian Liu}

\subjclass{35P20, 53C44, 58J35, 58J50\\   {\it Key words and phrases}. Sobolev trace inequality; Dirichlet-to-Neumann operator; Steklov eigenvalue;
    heat kernel associated to Dirichlet-to-Neumann operator; Poisson upper bound; Asymptotic expansion;
   Curvature}

\maketitle Department of Mathematics, Beijing Institute of
Technology,
 Beijing 100081, the People's Republic of China.
 \ \
E-mail address:  liugqz@bit.edu.cn

\vskip 0.46 true cm

\vskip 15 true cm

\begin{abstract}    For a given bounded domain $\Omega$ with smooth boundary in a smooth Riemannian manifold $(\mathcal{M},g)$, we show that the Poisson type upper-estimate of the heat kernel associated to the Dirichlet-to-Neumann operator, the Sobolev trace inequality, the Log-Sobolev trace inequality, the Nash trace inequality, and the Rozenblum-Lieb-Cwikel type inequality are all equivalent. Upon decomposing the Dirichlet-to-Neumann operator into a sum of the square root of the Laplacian and a pseudodifferntial operator and  by applying Grubb's method of symbolic calculus for the corresponding pseudodifferential heat kernel operators, we  establish a procedure to calculate all the coefficients of the asymptotic expansion of the trace of the heat kernel associated to Dirichlet-to-Neumann operator as $t\to 0^+$. In particular, we explicitly give the first four coefficients of this asymptotic expansion. These coefficients give precise information regarding the area and curvatures of the boundary of the domain in terms of the spectrum of the Steklov problem.
   \end{abstract}

\vskip 1.39 true cm

\section{ Introduction}

\vskip 0.45 true cm

An interesting question in Analysis is the following: what is the relationship between the geometrical quantitative characteristics of a bounded domain and the spectrum of the Steklov problem on the boundary of the domain? This problem originates from inverse spectral problems, where the known data is the spectrum of a differential or a pseudodifferential operator and one wishes to recover the geometry of a manifold. A rather efficient method to deal with this question is the asymptotic expansion of the trace of the heat kernel of the Dirichlet-to-Neumann (also called the Steklov-Poincar\'{e}) operator (the so-called ``heat kernel method''). The coefficients of the asymptotic expansion of the trace of such a heat kernel not only provides spectral invariants but also gives some geometric and topological information. The heat kernel method is closely related to the Calder\'{o}n problem and also has many applications in shape recognition, detection of distant physical objects (such as stars or atoms, or moving objects) from the light or sound they emit, identification of inhomogeneities (for instance in medical imaging) through the determination of different conductivities  (see, \cite{IMS} and \cite{JJo}, \cite{SU2}, \cite{JGP}).

\vskip 0.07 true cm

Let $\mathcal{M}$ be an $(n+1)$-dimensional, smooth, complete Riemannian
manifold
 with  metric tensor $g=(g_{jk})$, and let $\Omega$ be
 a bounded domain in $\mathcal{M}$ with smooth boundary $\partial \Omega$.
 The Laplace-Beltrami operator associated with the metric $g$ is given in
local coordinates by
 \begin{eqnarray} \label{1-1} \triangle_g u=\frac{1}{\sqrt{|g|}}\sum_{j,k=1}^{n+1}
   \frac{\partial}{\partial x_j} \left( \sqrt{|g|}\,g^{jk}
   \frac{\partial u}{\partial x_k}\right),\end{eqnarray}
 where $(g^{jk})$ is the inverse of the metric tensor
 $(g_{jk})$ and $|g| =\mbox{det}\,g$. For the Dirichlet problem associated with
(\ref{1-1}), \begin{eqnarray}  \label {1-2}\left\{ \begin{array}{ll}
\Delta_g u =0 \quad \; &
\mbox{in}\;\; \Omega, \\
 u=\phi \quad \;  &\mbox{on}\;\; \partial \Omega,\end{array} \right.\end{eqnarray}
  we denote the Dirichlet-to-Neumann
 operator in this case by \begin{eqnarray} \label {1-3} {\mathcal{N}}_g\phi=
\frac{\partial u}{\partial \nu}, \end{eqnarray}
 where $\nu=(\nu_1, \cdots, \nu_{n+1})$ denotes the unit inward normal to $\partial \Omega$.

The Dirichlet-to-Neumann operator is a self-adjoint, first order elliptic
pseudodifferential operator (see p.$\,$37-38 of \cite{Ta2}).
 Prototypical in inverse problems, the Dirichlet-to-Neumann operator is  related to the
Calder\'{o}n problem \cite{Cal} of determining the anisotropic conductivity of a body from current and
voltage measurements at its boundary. The Calder\'{o}n problem has been solved affirmatively for the isotropic conductivity of a body by J. Sylvester and G. Uhlmann in higher dimensional case ($\mbox{dim}\,M\ge 3$) \cite{SU1}, and by A. Nachman for two dimensional case \cite{Nac1}. Generalizations to less regular conductivities has been obtained by a number of authors (see, \cite{Ale}, \cite{Bro},  \cite{Cha}, \cite{GLU}, \cite{Nac1}, \cite{NSU}, \cite{PPU}, \cite{BT}, \cite{Nov}, \cite{Nac2}, \cite{BU}, \cite{KT}, \cite{SU2},  \cite{AP}).
  Unfortunately, ${\mathcal{N}}_g$ doesn't determine $g$ uniquely for general
 $n$-dimensional Riemannian manifold (This observation is due to L. Tartar, see
\cite{KV} for an account).
Because $\partial \Omega$ is
compact, the spectrum of ${\mathcal{N}}_g$ is nonnegative, discrete and
unbounded (see p.$\,$95 of \cite{Ban}, \cite{FSc}). The spectrum $\{\lambda_k\}_{k=1}^\infty$ of this operator
is just the Steklov spectrum of the domain $\Omega$.
More precisely, \begin{eqnarray*} \left\{ \begin{array}{ll}  \Delta_g u_k =0 \quad \; \quad \mbox{in}\;\; \Omega, \\
\frac{\partial u_k}{\partial \nu}=-\lambda_k u_k \quad \; \mbox{on}\;\; \partial \Omega,\end{array} \right.\end{eqnarray*}
where $u_k$ is the eigenfunction corresponding to the $k$-th Steklov eigenvalue $\lambda_k$.
 The study of
the spectrum of ${\mathcal{N}}_g$ was initiated by Steklov in 1902 (see \cite{St}). Eigenvalues and eigenfunctions
of this operator are used in fluid mechanics, heat transmission and vibration problems (see \cite{FK}
and \cite{KK}). Denote by $\omega_n$ the volume of the unit ball in ${\Bbb R}^{n}$.
A famous asymptotic formula of L. Sandgren \cite{Sa} (This formula was first established by Sandgren, and a sharp form was given by the author in \cite{Liu2}) states that
   \begin{eqnarray} \label{1-4}  N(\tau)=
 \#\{k\big|\lambda_k \le \tau\} =
   \frac{\omega_{n}\big(\mbox{vol}(\partial \Omega)\big)\tau^{n}}{(2\pi)^{n}}
  +o(\tau^{n}) \quad \;\mbox{as}\;\;
 \tau\to +\infty,\end{eqnarray}
 or, what is the same,
 \begin{eqnarray}
\label{1-.5} \quad \quad \quad \; Z= \mbox{Tr}\; e^{t{\mathcal{N}}_g} = \sum_{k=1}^\infty e^{-\lambda_k
t} \sim \frac{\Gamma(n+1)\, \omega_n (\mbox{vol} (\partial \Omega))}{(2\pi)^{n}\,t^n}\,
  =\frac{\Gamma(\frac{n+1}{2})\, \mbox{vol}(\partial \Omega)}{\pi^{\frac{n+1}{2}}t^n}\;  \quad \;
\mbox{as}\;\; t\to 0^+.\end{eqnarray}
Asymptotic formula (\ref{1-.5}) shows that you can
hear {\it the area of} the boundary $\partial \Omega$ by the first term of asymptotic expansion for Tr$\,e^{t{\mathcal{N}}_g}$. Therefore our problem just asks: how can one obtain more terms in the asymptotic expansion for the Tr$\,e^{t{\mathcal{N}}_g}$ ?  This problem is quite similar to the well-known Kac problem for the Laplacian on a domain (The Kac question asks:
is it possible to hear the shape of a domain just by ※hearing§
all of the eigenvalues of the Dirichlet Laplacian? see \cite{BB}, \cite{vdB},
\cite{R.C}, \cite{Iv}, \cite{Kac}, \cite{LFP}, \cite{AGMT}, \cite{GW}, \cite{GWW}, \cite{Lo}, \cite{ANPS}, \cite{Sar}, \cite{CH}, \cite{Ch1}, \cite{CLN} and the references therein).

In order to explain the main method of this paper, we briefly review the historical
background for the case of Dirichlet Laplacian on domains. In 1910, H. A. Lorentz  conjectured that
for a two-dimensional domain $\Omega\subset {\Bbb R}^2$, the
 asymptotics of the counting function of the Dirichlet eigenvalues $\{\mu_k\}$ are given by:
  \begin{eqnarray} \label{116}  N_D(\tau)= \# \{k\big|\mu_k \le \tau\}
  = \frac{\mbox{vol} (\Omega)}{2\pi}\tau +o(\tau)\quad \,\,
  \mbox{as}\;\; \tau\to \infty.\end{eqnarray}
This asymptotic in particular implies that $\mbox{vol}(\Omega)$
 is a spectral invariant. Lorentz's conjecture was proved in 1913 by Hermann
Weyl  (see \cite{We1} and \cite{We2}).
With  Weyl's  formula as a starting point, Pleijel \cite{Ple1} in 1954  showed that
\begin{eqnarray} \label{0.1} \sum_{k=1}^\infty e^{-\mu_k t} \sim \frac{\mbox{vol}(\Omega)}{2\pi t} -\frac{\mbox{vol}(\partial \Omega)}{4}
 \,\frac{1}{\sqrt{2\pi t}} \quad \; \mbox{as}\;\; t\to 0^+,\end{eqnarray}
where $\mbox{vol}\, (\partial \Omega)$ is the length of boundary $\partial \Omega$. By a Tauberian theorem the asymptotic formula for the first term on the right side of (\ref{0.1})
is equivalent to  Weyl's formula (\ref{116}). For simply connected domains Pleijel established the formula
\begin{eqnarray} \label{0.2}  \sum_{k=1}^\infty e^{-\mu_k t} \sim \frac{\mbox{vol}(\Omega)}{2\pi t} +\frac{\mbox{vol}(\partial \Omega)}{4} \,
 \frac{1}{\sqrt{2\pi t}} +\frac{1}{6} \quad \; \mbox{as} \;\; t\to 0^+. \end{eqnarray}
 Kac \cite{Kac} used a combination of probability techniques and heat equation methods to establish (\ref{0.1}) for convex domains, and he obtained (\ref{0.2})
as a limiting case of convex polygonal domains. Kac also conjectures that for  multiply connected domains in ${\Bbb R}^2$  with $r$ holes,
the number $\frac{1}{6}$ in (\ref{0.2}) should be replaced by $\frac{1}{6}(1-r)$. McKean and Singer  in a celebrated paper \cite{MS} gave an
affirmative answer to the conjecture of Kac with respect to the third term for multiply connected domains in $n$-dimensional Riemannian manifold (with or without boundary). McKean and Singer \cite{MS} also obtained information about the curvature of the boundary of
 $\Omega$, which showed that the Euler characteristic
 $\chi(\Omega)$  is also a spectral invariant. Gilkey \cite{Gil} explicitly calculated the first four coefficients of the expansion of the trace of the heat kernel.  In 1991, Gordon, Webb and Wolpert \cite{GWW}, found examples of pairs
of distinct plane domains with the same spectrum (there are higher dimensional examples for a similar problem
by Gordon-Webb \cite{GW}). In \cite{MS}, the heat kernel estimates of the Laplacian play an important role in the asymptotic expansion.
 The two-sided Gaussian estimates for the heat kernel
associated with a uniformly elliptic operator in ${\Bbb R}^n$ was proved by  Aronson \cite{Ar} (see also \cite{PE}, \cite{Dav3}, \cite{CLY},  \cite{Grig}, \cite{Ouh}, \cite{FS}):
\begin{eqnarray} \label{00-00} \frac{c}{t^{n/2}}e^{-|x-y|^2/ct} \le G(t,x,y)
\le \frac{C}{t^{n/2}} e^{-|x-y|^2/Ct},\end{eqnarray}
where $c$ and $C$ are two positive constants.
Varadhan \cite{Va1}, \cite{Va2} first realized that the Riemannian distance should be used instead.
 His
result implies that, on any manifold,
\begin{eqnarray*} \lim_{t\to 0^+} t\, \ln G(t,x,y)= -\frac{d^2(x,y)}{4},\end{eqnarray*}
where $d(x,y)$ is the Riemannian distance between $x$ and $y$ (see also, \cite{Nor}).
For a complete Riemannian manifold with Ricci curvature bounded below,
 Li-Yau \cite{LY} and Sturm \cite{Stur} obtained upper and lower Gaussian estimates on the heat kernel,
from  which Varadhan's asymptotic  result follows immediately.
 Varopoulos \cite{Var} proved that the Sobolev inequality is not only sufficient but also
necessary for the on-diagonal upper bound of heat kernel of Laplacian.
 Carlen, Kusuoka and Stroock \cite{CKSt} proved that the upper bound of heat kernel  is  equivalent to
the Nash inequality (see also \cite{Gri}).
 Davies \cite{Dav1}, \cite{Dav2}, \cite{Dav3} proved that the on-diagonal upper bound of the heat kernel is also equivalent
to the log-Sobolev inequality. In addition, the Sobolev inequality is equivalent to
the Rozenblum-Lieb-Cwikel inequality (see \cite{Roz}, \cite{Lie}, \cite{Cwi}, \cite{LY2}, \cite{Con}, \cite{LW} and \cite{LSo}).

\vskip 0.10 true cm

 Let us come back to the Steklov spectrum.
 Because ${\mathcal{N}}_g$ is a pseudodifferntial operator, the corresponding problems become much more difficult than those of the Laplacian.
 To get more geometric information from the Steklov eigenvalues on $\partial \Omega$, we consider
 the heat kernel associated to the Dirichlet-to-Neumann operator on $\partial \Omega$:
  \begin{eqnarray}  \label {1-5}\left\{ \begin{array}{ll} \frac{\partial u(t, x)}{\partial t}= {\mathcal{N}}_gu(t,x)  \;\;
\quad \; &\mbox{in}\;\;  [0, +\infty)\times \partial\Omega, \\
 u(0, x)=\phi(x)   \;\;\quad \; & \mbox{on}\;\; \partial \Omega,\end{array}\right.\end{eqnarray}
    where ${\mathcal{N}}_g$ is the Dirichlet-to-Neumann operator on $\partial \Omega$.
     By perturbation of a fractional Laplacian, Gimperlein and Grubb in \cite{GG} have recently proved the following important
  Poisson upper and lower bounds for the heat kernel ${\mathcal{K}}(t,x,y)$ of the Dirichlet-to-Neumann operator:
   \begin{eqnarray*}  &&  {\mathcal{K}}(t,x,y) \le \frac{Ct}{(t^2 +d^2(x,y))^{1/2}} \left[ \frac{1}{(t^2+d^2(x,y))^{n/2}}+1\right]
  \quad \mbox{for all}\;\; t>0\;\; \mbox{and}\;\; x,y\in \partial \Omega, \\ &&
   {\mathcal{K}}(t,x,y)\ge \frac{ct}{(t^2 +d^2(x,y))^{(n+1)/2}} \quad\;\, \mbox{for} \;\; t+d(x,y)<r, \end{eqnarray*}
  where $r>0$ is some constant (The Poisson upper bound estimate
  $${\mathcal{K}}(t,x,y) \le C(t\wedge 1)^{-n} \big(1+\frac{d(x,y)}{t}\big)^{-n-1}  \;\; \mbox{for all}\;\; x,y\in \partial \Omega \;\;\mbox{and} \,\, t>0,$$
 of the kernel was also obtained by Elst and Ouhabaz in \cite{EO}).

In the first part of this paper, combining classical methods and some new techniques we show
 that the Poisson type upper-estimate of the heat kernel associated to the Dirichlet-to-Neumann operator is equivalent to the Sobolev trace inequality,
   to the Log-Sobolev trace inequality, to the Nash trace inequality, and to the Rozenblum-Lieb-Cwikel type inequality.

For the coefficients of asymptotic expansion
of the trace of the heat kernel associated to the Dirichlet-to-Neumann operator, the first coefficient $a_0(n,x)$ had been known by (\ref{1-.5}); $a_1(2, x)$ had been obtain in \cite{EW}; the coefficients $a_1(n,x)$ and $a_2(n, x)$ were explicitly calculated in \cite{PS} and \cite{Liu3} in  different ways. In the second part of the paper, by decomposing the Ditrichlet-to-Neumann operator and by applying Grubb's method (see \cite{Gr}), we establish a procedure, from which all coefficients of the asymptotic expansion of the trace of the heat kernel associated to the Dirichlet-to-Neumann operator can be calculated as $t\to 0^+$. In particular, we explicitly give the first four coefficients of the asymptotic expansion.
  This provides the information for the area and curvature of the boundary.

  The main ideas are as follows:
  we first decompose the Dirichlet-to-Neumann operator ${\mathcal{N}}_g$ into a sum of the negative square root $-\sqrt{-\Delta_h}$ of the Laplacian  and a pseudodifferential operator
 $B$ on $\partial \Omega$, where $h$ is the induced metric on $\partial \Omega$ by $g$, then we calculate their principal, second, third and  $(M-1)$-th symbols for $-\sqrt{-\Delta_h}$ and $B$. Because the heat kernels of $\Delta_h$ and $-\sqrt{-\Delta_h}$ have an analytical formula, we can explicitly calculate the heat kernel ${\mathcal{K}}_V(t,x,y)$ of $-\sqrt{-\Delta_h}$. Therefore the heat kernel $\mathcal{K}(t,x,y)$ of ${\mathcal {N}}_g$ can be expanded into the form ${\mathcal{K}}_{V} (t,x,y) + {\mathcal{K}}_{V_{-2}} (t,x,y) +\cdots +{\mathcal{K}}_{V_{-M}} (t,x,y)
  + {\mathcal{K}}_{V'_M} (t,x,y)$ as $t\to 0^+$, where ${\mathcal{K}}_{V_{-2}} (t,x,y)$, ${\mathcal{K}}_{V_{-3}} (t,x,y)$ and ${\mathcal{K}}_{V_{M}} (t,x,y)$ corresponds to the principal, second and $(M-1)$-th symbols of $e^{tB}$. By considering the trace of the heat
 kernel of the Dirichlet-to-Neumann operator,
 \begin{eqnarray} \mbox{Tr} \, e^{t{\mathcal{N}}_g}=\int_{\partial\Omega} \mathcal{K}(t,x,x)dS(x) =\sum_{k=1}^\infty e^{-\lambda_k t},
  \end{eqnarray}
  and by a remainder estimate we eventually obtain the following asymptotic expansion:
    \begin{eqnarray*} &&\int_{\partial \Omega} {\mathcal{K}}(t,x,x)dS(x) = t^{-n} \int_{\partial \Omega}\frac{\Gamma(\frac{n+1}{2})}{\pi^{\frac{n+1}{2}}}\, dS(x)+ t^{1-n} \int_{\partial \Omega} a_1(n,x) \, dS(x) \\
   &&\quad\; + \cdots + t^{M-1-n} \int_{\partial \Omega} a_{M-1}(n,x)\, dS(x) +\left\{
           \begin{array}{ll} O(t^{M-n})\\ O(t\log t),\end{array}\right. \quad \;\mbox{as}\;\; t\to 0^+,\end{eqnarray*}
            where  $M$ can be taken as $2,3,4,\cdots$ with $n\ge M-1$, and the notation $O(t^m)$ (respectively, $O(t\log t)$) denotes a function which satisfies $|O(t^m)|\le c_0 t^m$ (respectively, $|O(t\log t)|\le c_0t\log t$) for some constant $c_0>0$ and all $t>0$.
           The first coefficient   $a_0(n,x)=\frac{\Gamma(\frac{n+1}{2})}{\pi^{\frac{n+1}{2}}}$ is independent of $n$ and $x$. The second coefficient $a_1(n,x)$ depends only on the ``area'' $\mbox{vol}(\partial \Omega)$  and the mean curvature of the boundary $\partial \Omega$. The coefficient $a_2(n,x)$ depends not only on the ``area＊' $\mbox{vol}(\partial \Omega)$ and the principal curvatures $\kappa_1,\cdots,\kappa_n$, but also on the scalar curvature ${\tilde{R}}_{\Omega}$ (respectively $R_{\partial \Omega}$) of $\Omega$ (respectively $\partial \Omega$). Finally the fourth coefficient $a_3(n,x)$ depends on $\mbox{vol}(\partial \Omega)$, the principal curvatures, the Ricci tensor ${\tilde{R}}_{jj}$ (respectively $R_{jj}$) and the scalar curvature ${\tilde R}_\Omega$ (respectively $R_{\partial \Omega}$ of $\Omega$ (respectively $\partial \Omega$) as well as the covariant derivative $\sum_{j=1}^n {\tilde R}_{j(n+1)j(n+1),(n+1)}$ of the Ricci curvature with respect to $(\bar \Omega, g)$ (see Theorem 6.1). Generally, $a_{M-1} (n,x)$ can be represented in terms of the geometical quantities (see, (ii) of Remark 6.2).
             This asymptotic expansion shows that one can hear the ``area'' of $\partial \Omega$ and $\int_{\partial \Omega} a_{M-1}(n,x)dS(x)$ ($M=1,2,3,\cdots$)  by ``hearing'' all of the Steklov eigenvalues.
             It also implies that  $\int_{\partial \Omega} a_{M-1}(n,x) dS(x)$ ($M\ge 1$) are all spectral invariants, from which we know that two domains with different spectral invariants just mentioned above, can never have the same Steklov spectrum. As a by-product, we explicitly give the first four coefficients of the asymptotic expansion of the trace of the corresponding heat kernel associated to the negative square root $-\sqrt{-\Delta_h}$ of the Laplacian, which also tells us that the first four terms of asymptotic expansion of the heat kernel can be obtained not only for any smooth elliptic partial differential operator of degree $2,3,4,5$ (see, p.$\,$45 of \cite{MS} and p.$\,$613 of \cite{Gil}) but also for the elliptic pseudodifferential operator of fractional-order $\frac{1}{2}$.

\vskip 1.35 true cm

\section{Preliminaries}

\vskip 0.45 true cm

Let $\Omega$ be a bounded domain with smooth boundary in a smooth Riemannian manifold $(\mathcal{M},g)$.
$\mathcal{K}$ is said to be a fundamental solution of $\frac{\partial u}{\partial t}={\mathcal{N}}_g u$,
  if for any fixed $y\in \partial \Omega$,
  \begin{eqnarray} \label{55-2} \left\{ \begin{array}{ll} \frac{\partial \mathcal{K}(t, x,y)}{\partial t} = {\mathcal{N}}_g \mathcal{K}(t, x, y),
  \quad  t>0, \,\, x\in \partial \Omega,\\
  \mathcal{K}(0, x, y)= \delta(x-y),\end{array}\right. \end{eqnarray}
  where ${\mathcal{N}}_g$ acts on the $x$ variable and $\delta (x-y)$ is the delta function concentrated at $y$.
  $\mathcal{K}$ is also called the heat kernel associated to the Dirichlet-to-Neumann operator (or the kernel of the semigroup $e^{t{\mathcal{N}}_g}$).
 The initial condition in (\ref{55-2}) means that for every continuous function $\phi(x)$ on $\partial \Omega$, if $x\in \partial \Omega$
  then \begin{eqnarray}  \label{4b2} \lim_{ t\to 0^+} \int_{\partial \Omega} \mathcal{K}(t,x,y) \phi(y) \,dS(y) =\phi(x). \end{eqnarray}

For the special manifold $({\Bbb R}^{n+1}_+, g)$, where ${\Bbb R}^{n+1}_+$ is the upper half-space $\{(x_1, \cdots, x_{n+1})\in {\Bbb R}^{n+1}
\big| x_{n+1}\ge 0\}$, and
      \begin{gather} \label{3--1}  \left(g_{jk}\right)_{(n+1)\times (n+1)}=\begin{pmatrix} h_{11} & h_{12} & \cdots  & h_{1n} &0\\
      h_{21} &  h_{22} & \cdots & h_{2n} & 0\\
      \vdots & \vdots & \ddots & \vdots & \vdots \\
      h_{n1} & h_{n2} & \cdots & h_{nn}  & 0\\
      0 & 0 & \cdots   & 0 &  1  \end{pmatrix} \end{gather}
    is a positive-definite, real symmetric
  $(n+1)\times (n+1)$  constant matrix.  It is easy to verify that the function
  \begin{eqnarray} \label{3.2}   u(x, x_{n+1})=\int_{{\Bbb R}^n} P(x-z, x_{n+1})
   \phi(z)\,\sqrt{|h|}\,dz  \end{eqnarray}
       is a solution of the problem:
       \begin{eqnarray} \label{3.1} \left\{ \begin{array}{ll} \sum_{j,k=1}^{n+1} h^{jk}
 \frac{\partial^2 u}{\partial x_j\partial x_k}=0, \quad \;\; &\mbox{in}\;\; R^{n+1}_+,\\
 u=\phi,  \quad \;\; & \mbox{on}\;\; \partial R^{n+1}_+, \end{array} \right.\end{eqnarray}
               where $(h^{jk})$ is the inverse of $\left( h_{jk}\right)$, and
               the Poisson kernel
       is \begin{eqnarray}P(x,x_{n+1})=  \frac{\Gamma(\frac{n+1}{2})}{\pi^{\frac{n+1}{2}}}
        \,  \frac{x_{n+1}}{( \sum_{i,j=1}^n  h_{jk} x_j x_k
 + x_{n+1}^2)^{\frac{n+1}{2}}}.\end{eqnarray}
       This implies \begin{eqnarray*} \frac{\partial u}{\partial \nu}= \frac{\partial u} {\partial x_{n+1}}
 \bigg|_{x_{n+1}=0}
 =  \int_{{\Bbb R}^n} \frac{\Gamma(\frac{n+1}{2})}{\pi^{\frac{n+1}{2}}}
        \,  \left(\sum_{j,k=1}^n
 h_{jk} (x_j-z_j) (x_k-z_k)\right)^{-\frac{n+1}{2}} \phi(z) \, \sqrt{|h|}\,dz.\end{eqnarray*}
 In other words, the Dirichlet-to-Neumann operator
${\mathcal{N}}_g: H^{\frac{1}{2}}(\partial {\Bbb R}^{n+1}_+)\rightarrow H^{-\frac{1}{2}}(\partial {\Bbb R}^{n+1}_+)$
  is given by
\begin{eqnarray} \label {2.4} \qquad  \quad \quad {\mathcal{N}}_g\phi(x) =\int_{{\Bbb R}^n} \frac{\Gamma(\frac{n+1}{2})}{\pi^{\frac{n+1}{2}}}
        \,  \left(\sum_{j,k=1}^n
 h_{jk} (x_j-z_j) (x_k-z_k)\right)^{-\frac{n+1}{2}} \phi(z)\,
  \sqrt{|h|}\, dz,\end{eqnarray}
 from which we can see
\begin{eqnarray}  \label {3,1} {\mathcal{N}}_g\phi = \frac{\Gamma(\frac{n+1}{2})}{\pi^{\frac{n+1}{2}}}
        \, \sqrt{-\sum_{j,k=1}^n
 h_{jk} \frac{\partial^2}{\partial x_j\partial x_k}} \;\; \phi. \end{eqnarray}
 Furthermore by applying Fourier transform, with the aid of the basic method of the linear transform and the following formula \begin{eqnarray*}   \left(\frac{1}{2\pi}\right)^{\frac{n}{2}}
       \int_{{\Bbb R}^n} e^{-t |\xi|} e^{\pm i x\cdot \xi} d\xi = \frac{\sqrt{2^n}\,\Gamma\big(\frac{n+1}{2}\big)}{\sqrt\pi}
  \,\frac{t }{\left(t^2 +|x|^2\right)^{\frac{n+1}{2}}},\end{eqnarray*}
  we can get that
 for  a positive definite and symmetric constant matrix $(h_{jk})_{n\times n}$,  \begin{eqnarray} \label{3..11} u(t,x)=\frac{\Gamma(\frac{n+1}{2})}{\pi^{\frac{n+1}{2}}}
        \,  \frac{t}{\left(t^2
 +\sum_{j,k=1}^n h_{jk} x_jx_k \right)^{\frac{n+1}{2}}}\end{eqnarray}
 is a fundamental solution of
   \begin{eqnarray} \label{3,,21} \left\{ \begin{array}{ll} \frac{\partial u(t,x)}{\partial t} ={\mathcal{N}}_gu(t,x),\quad x\in \partial {\Bbb R}^{n+1}_+, \;t>0, \\
   u(0, x)= \delta(x).\end{array}\right. \end{eqnarray}

\vskip 0.12 true cm
For the Laplacian $-\Delta_h$ on $\partial \Omega$,  we can define the square root operator $\sqrt{-\Delta_h}$ (see \cite{Ta2}).
Stinga and Torrea (see \cite{StTo}) showed that if $u$ satisfies
\begin{eqnarray} \label{y,1} \left\{\begin{array}{ll}  \Delta_h u + \frac{\partial^2u}{\partial x_{n+1}^2}=0 \;\; & \mbox{in}\;\; \partial \Omega \times (0,\infty)\\
 u(x,0)=\phi(x)  \;\; & \mbox{for}\;\; x\in \partial \Omega,\end{array}\right.\end{eqnarray}
 then \begin{eqnarray} \label{y,2} \frac{\partial u}{\partial x_{n+1}}= - \sqrt{-\Delta_h} \,\phi(x)
\end{eqnarray}
 and  \begin{eqnarray} -\sqrt{-\Delta_h} \,\phi(x)= \frac{1}{2\sqrt{\pi}} \int_0^\infty \left(e^{\mu\Delta_h} \phi(x) -\phi(x)\right)
\frac{d\mu}{\mu^{3/2}},\end{eqnarray}
where $e^{t\Delta_h}$  is the heat semigroup generated by $\Delta_h$ on $\partial \Omega$.

\vskip 0.10 true cm

 If $U$ is an open subset  of ${\Bbb R}^n$, we denote  by $S^m_{1,0}=S^m_{1,0} (U,
{\Bbb R}^n)$ the set of all $p\in C^\infty (U, {\Bbb R}^n)$ such
that for every compact set $K\subset U$
 we have
 \begin{eqnarray} \label {-2.3} |D^\beta_x D^\alpha_\xi p(x,\xi)|\le C_{K,\alpha,
 \beta}(1+|\xi|)^{m-|\alpha|}, \quad \; x\in K,\,\, \xi\in {\Bbb R}^n\end{eqnarray}
 for all $\alpha, \beta\in {\Bbb N}^n_+$.
 The  elements  of $S^m_{1,0}$  are  called  symbols (or total symbols) of order $m$.
 It is clear that $S^m_{1,0}$ is a
Fr\'{e}chet space with semi-norms given by the smallest constants
which can be used in (\ref{-2.3}) (i.e.,
\begin{eqnarray*} \|p\|_{K,\alpha, \beta}=
 \,\sup_{x\in K}\bigg|\left(D_x^\beta D_\xi^\alpha
 p(x, \xi)\right)(1+|\xi|)^{|\alpha|-m}\bigg|).\end{eqnarray*}
 Let $p(x, \xi)\in S^m_{1,0}$. A pseudo-differential operator
in an open set $U\subset {\Bbb R}^n$ is essentially defined by a
Fourier integral  operator (cf. \cite{Ho4}, \cite{Sh}):
\begin{eqnarray} \label{-2.1}  P(x,D) u(x) = \frac{1}{(2\pi)^{\frac{n}{2}}} \int_{{\Bbb R}^n} p(x,\xi)
 e^{i x\cdot\xi} \hat{u} (\xi)d\xi,\end{eqnarray}
and denoted by $OPS^m$.
Here $u\in C_0^\infty (U)$ and $\hat{u} (\xi)=\left(\frac{1}{2\pi}\right)^{\frac{n}{2}}
 \int_{{\Bbb R}^n} e^{-i\langle y, \xi\rangle} u(y)dy$ is the Fourier
transform of $u$.
 If there are smooth $p_{m-j} (x, \xi)$, homogeneous in $\xi$ of degree $m-j$ for $|\xi|\ge 1$, that is,
  $p_{m-j} (x, r\xi) =r^{m-j} p_{m-j} (x, \xi)$ for $r, |\xi|\ge 1$, and if
 \begin{eqnarray} \label{-2.2} p(x, \xi) \sim \sum_{j\ge 0} p_{m-j} (x, \xi)\end{eqnarray}
 in the sense that
 \begin{eqnarray}  p(x, \xi)-\sum_{j=0}^N p_{m-j} (x, \xi) \in S^{m-N-1}_{1, 0}, \end{eqnarray}
 for all $N$, then we say $p(x,\xi) \in S_{cl}^m$, or just $p(x, \xi)\in S^m$. We call
 $p_m(x, \xi)$, $p_{m-1}(x,\xi)$ and $p_{m-2} (x,\xi)$ the principal, second and third symbols of $P(x, D)$, respectively.
 Let $\mathcal{M}$ be a smooth $n$-dimensional Riemannian
 manifold (of class $C^\infty$). We will denote by $C^\infty(\mathcal{M})$ and $C_0^\infty(\mathcal{M})$
 the space of all smooth complex-valued functions on $\mathcal{M}$  and the subspace of all
 functions with compact support, respectively. Assume that we are given a linear
 operator
 \begin{eqnarray*} P: C^\infty_0(\mathcal{M}) \to C^\infty(\mathcal{M}).\end{eqnarray*}
 If $G$ is some chart in $\mathcal{M}$ (not necessarily connected) and $\kappa: G\to U$ its
diffeomorphism onto an open set $U \subset {\Bbb R}^n$, then let ${\tilde P}$ be defined by the diagram
\begin{eqnarray*}
\begin{CD}
C_0^\infty (G) @> P>>  C^\infty(G) \\
@A\kappa^* AA @AA\kappa^* A \\
 C_0^\infty(U) @> \tilde P>> C^\infty(U)
\end{CD}
\end{eqnarray*}
where
  $\kappa^*$ is the induced transformation
from $C^\infty(U)$ into $C^\infty(G)$, taking a
function $u$ to the function $u\circ \kappa$.
  (note, in the upper row is the operator $r_G\circ P\circ i_G$, where $i_G$ is the natural
embedding $i_G: C_0^\infty (G)\to C^\infty_0 (M)$ and $r_G$ is the natural restriction $r_G: \, C^\infty(M)\to
C^\infty(G)$; for brevity we denote this operator by the same letter $P$ as the
original operator).
 An operator $P: C_0^\infty (\mathcal{M}) \to C^\infty(\mathcal{M})$ is called a pseudodifferential
operator on $\mathcal{M}$ if for any chart diffeomorphism $\kappa: G\to U$, the
operator $\tilde P$ defined above is a pseudodifferential
operator on $U$.
 An operator $P$ is said to be an elliptic pseudodifferential operator of order $m$ if
  for every compact $K\subset \Omega$ there exists a  positive constant $c=c(K)$ such that \begin{eqnarray*} |p(x, \xi)\ge c|\xi|^m, \,x\in K,\, |\xi|\ge 1\end{eqnarray*} for any compact set $K\subset \Omega$.
  If $P$ is a non-negative elliptic pseudodifferential operator of order $m$, then the spectrum  of $P$ lies in a right half-plane and has a finite lower bound $\gamma(P) =
\inf\{\mbox{Re}\, \lambda\big| \lambda\in \sigma(P)\}$. We can modify $p_m(x, \xi)$ for small $\xi$ such that
$p_m(x,\xi)$  has a positive
lower bound throughout and lies in $\{\lambda=re^{i\theta} \big| r>0, |\theta|\le \theta_0\}$, where $\theta_0\in (0,\frac{\pi}{2})$.
According to \cite{Gr}, the resolvent $(P-\lambda)^{-1}$ exists and is holomorphic in $\lambda$ on a neighborhood of a set
\begin{eqnarray*}  W_{r_0,\epsilon} =\{\lambda\in {\Bbb C} \big| |\lambda|\ge r_0, \mbox{arg}\, \lambda \in [\theta_0+\epsilon, 2\pi -\theta_0-\epsilon], \,\mbox{Re}\, \lambda \le \gamma (P)-\epsilon\}\end{eqnarray*}
(with $\epsilon>0$). There exists a parametrix $Q'_\lambda$ on a neighborhood of a possibly larger set
(with $\delta>0,\epsilon>0$)
\begin{eqnarray*} V_{\delta,\epsilon} =\{ \lambda \in {\Bbb C} \big| |\lambda| \ge \delta \;\;\mbox{or arg}\, \lambda \in [\theta_0+\epsilon, 2\pi-\theta_0-\epsilon]\}\end{eqnarray*}
such that this parametrix coincides with $(P-\lambda)^{-1}$ on the intersection. Its symbol $q(x,\xi, \lambda)$
in local coordinates is holomorphic in $\lambda$ there and has the form (cf. Section 3.3 of \cite{Gr})
\begin{eqnarray} \label {-2.6} q(x, \xi,\lambda) \sim \sum_{l\ge 0} q_{-m-l} (x, \xi, \lambda)\end{eqnarray}
where \begin{eqnarray} \label{-2.7} && q_{-m}= (p_m(x, \xi) -\lambda)^{-1}, \quad \;  q_{-m-1} = b_{1,1}(x, \xi) q^2_{-m},\\
&& \, \cdots, \,  q_{-m-l} = \sum_{k=1}^{2l} b_{l,k} (x, \xi) q^{k+1}_{-m}, \cdots, \nonumber\end{eqnarray}
with symbols $b_{l,k}$ independent of $\lambda$  and homogeneous of degree $mk-l$ in $\xi$ for $|\xi|\ge1$.
 The semigroup $e^{-tP}$ can be defined from $P$ by the Cauchy integral formula (see p.$\,$4 of \cite{GG}):
 \begin{eqnarray*} e^{-tP} =\frac{i}{2\pi} \int_{\mathcal{C}} e^{-t\lambda} (P-\lambda)^{-1} d\lambda,\end{eqnarray*}
where $\mathcal{C}$ is a suitable curvature in the positive direction around the spectrum of $P$.

\vskip 1.58 true cm

\section{Upper estimates of the heat kernel, Sobolev trace inequality, Log-Sobolev trace inequality,  Nash trace inequality and
Rozenblum-Lieb-Cwikel type inequality}

\vskip 0.5 true  cm

We now can establish the equivalence for the Poisson type upper-estimate of the heat kernel associated to the Dirichlet-to-Neumann operator with a number of other inequalities.
\vskip 0.2 true cm

\noindent {\bf Theorem 3.1.}  \ \  {\it  Let $(\mathcal{M},g)$ be an $n+1$ dimensional Riemannian manifold ($n \ge 2$),
and  let $\Omega$ be a bounded domain with smooth boundary. Then the following statements are equivalent:

(i) \ \   Sobolev trace inequality:    there exist positive constants $A$ and $B$  such that for all $v\in W^{1,2}(\Omega)$,
 \begin{eqnarray} \label{5=1} \left( \int_{\partial \Omega} |v|^{\frac{2n}{n-1}} dS \right)^{\frac{n-1}{n}} \le A \int_\Omega |\nabla_g v|^2 dV   +B
 \int_{\partial \Omega} |v|^2 dS, \end{eqnarray}
 where $dV$ and $dS$ are the Riemannian elements of volume and area on $\Omega$ and $\partial \Omega$, respectively;

 (ii) \ \   Log-Sobolev trace inequality:  for all $v\in W^{1,2}(\Omega)$ and all $\epsilon >0$,
 \begin{eqnarray} \label{5=2}\int_{\partial \Omega} v^2 \ln |v|  \,dS \le \epsilon \int_{\Omega} |\nabla_g v|^2 dV  +
 \beta (\epsilon)  \|v\|_{L^2(\partial \Omega)}^2  + \|v\|_{L^2(\partial \Omega)}^2  \ln \|v\|_{L^2(\partial \Omega)},\end{eqnarray}
 where $\beta(\epsilon)= \frac{n}{2} \ln \frac{nA}{2e} +BA^{-1}\epsilon - \frac{n}{2}\ln \epsilon$;

  (iii) \ \  Nash trace inequality:  for all $v\in W^{1,2}(\Omega)$,
   \begin{eqnarray}  \label{5=3} \|v\|_{L^2(\partial \Omega)}^{2+\frac{2}{n}} \le \left( A\|\nabla_g v\|_{L^2(\Omega)}^2 +B\| v\|_{L^2(\partial \Omega)}^2
   \right) \| v\|_{L^1 (\partial \Omega)}^{2/n};\end{eqnarray}

    (iv) \ \  On-diagonal Dirichlet-to-Neumann heat kernel upper bound:
     \begin{eqnarray} \label{5=4} {\mathcal{K}} (t, x, y) \le  \left(nAe/4\right)^n \, \frac{e^{BA^{-1} t}}{t^n} \,\, \mbox{for all} \;\; x,y\in \partial \Omega
     \;\; \mbox{and}\;\; t>0; \end{eqnarray}
  $\quad \quad$

   (v) \ \  Off-diagonal Dirichlet-to-Neumann heat kernel upper bound:
   \begin{eqnarray} \label{08.08.1}{\mathcal{K}}(t,x,y) \le \frac{Ct}{(t^2+d^2(x,y))^{\frac{n+1}{2}}} \quad \, \mbox{for all}\;\; x,y\in \partial \Omega\,\,
   \mbox{and} \,\, 0\le t\le 1.\end{eqnarray}

 (vi) \ \  The Rozenblum-Lieb-Cwikel type inequality: \ let $q\in L^1_{loc} (\partial \Omega)$ and $q_- \in L^{n} (\partial \Omega)$, where $q_-(x)=-\min\{q(x),0\}$. For $\alpha\ge 0$, let $I_q(\alpha)$ be the number of eigenvalues $\lambda$ satisfying
  \begin{eqnarray} \label{5=6a} (A{\mathcal{N}}_g -B-q) \phi (x) =-\lambda \phi(x) \quad \, x\in \partial \Omega \end{eqnarray}
  with $\lambda \le -\alpha$, where $A$ and $B$ are the positive constants in the Sobolev trace inequality (\ref{5=1}). Then there exists a positive constant $C(n)$ depending only on $n$ such that \begin{eqnarray} \label{5=5a}  I_q(0)\le C(n) \int_{\partial \Omega} q_{-}^{n} \,dS.\end{eqnarray}}

\vskip 0.2 true cm

\noindent {\bf Remark 3.2.}  \   (a) \  The sharp Sobolev trace inequality was proved by Y.-Y. Li and M. Zhu (see Theorem 0.1 of \cite{LZ}). Li and Zhu proved that the constant $A=2(n-1)^{-1} \big((n+1)\omega_{n+1}\big)^{-\frac{1}{n}}$ in the front of $\int_{\partial \Omega} |\nabla_g v|^2 dV$ is optimal. It cannot be replaced by any smaller number.

(b) \  In the Rozenblum-Lieb-Cwikel type inequality, we can take $C(n)=e^n$ (see the proof for (i)$\Rightarrow$(vi)). If we replace the Sobolev constants $A$ and $B$ in (\ref{5=6a}) by any two positive real numbers ${\tilde A}$ and ${\tilde B}$, then the corresponding inequality (\ref{5=5a}) still holds; however, the constant $C(n)$ should be replaced by a new constant $C(n,{\tilde A},{\tilde B})$  which depends only on $n$, ${\tilde A}$ and ${\tilde B}$.

 (c) \   The estimate (\ref{5=5a}) is sharp, in the following sense.  Replacing $q$ by $\beta q$ in  (\ref{5=6a}), where $\beta>0$ is a large parameter, we derive from (\ref{5=5a})
 \begin{eqnarray*} I_{\beta q} (0) \le C(n) \beta^n \int_{\partial \Omega} q^n_-dS.\end{eqnarray*}
Similar to Rozenblum's method \cite{Roz} (see also, \cite{LSo}), it is easy to prove that for any $q\in C^\infty(\partial \Omega)$, $ I_{\beta q} (0)$ has the asymptotic formula:
\begin{eqnarray} \label{08-20.5}\lim_{\beta \to +\infty} \beta^{-n} I_{\beta q} (0) =\frac{\omega_n}{(2\pi)^n}\int_{\partial \Omega} \frac{q^n_-}{A^n} dS.\end{eqnarray} Here the asymptotic constant $\frac{\omega_n }{(2\pi)^n}$ is just the corresponding coefficient in Sandgren's asymptotic formula (\ref{1-4}).

\vskip 0.35 true cm

 \noindent  {\bf Proof of Theorem 3.1.} \ \    (i) $\Longrightarrow$ (ii):
  Given $f\in W^{1,2}(\Omega)$ such that $\|f\|_{L^2(\partial \Omega)}=1$, we introduce the measure
  \begin{eqnarray*} dW(x)= f^2(x)dS(x). \end{eqnarray*}
  It follows from  $\int_{\partial \Omega} f^2 dS =1$ that $\int_{\partial \Omega} dW=1$.
  Since $\ln \phi$ is a concave function of $\phi$, by applying Jensen's inequality we
 have $\int_{\partial \Omega} \ln \phi \,dW \le \ln \int_{\partial \Omega} \phi\, dW$.
    Putting $\phi= |f|^{q-2}$ with $q=\frac{2n}{n-1}$, we get
   \begin{eqnarray*} \int_{\partial \Omega}  (\ln |f|^{q-2}) f^2 dS \le \ln \int_{\partial \Omega} |f|^{q-2} f^2 dS =
  \ln \|f\|_{L^q(\partial \Omega)}^q,\end{eqnarray*}
     i.e.,
     \begin{eqnarray*} \int_{\partial \Omega}  f^2 \ln |f|\, dS \le \frac{q}{q-2} \ln \|f\|_{L^q(\partial \Omega)}.\end{eqnarray*}
     Since $\frac{q}{q-2} = n$, by applying the Sobolev trace inequality (\ref{5=1}) we have
     \begin{eqnarray*} \int_{\partial \Omega}  f^2 \ln f^2\, dS &\le & n \ln \|f\|_{L^q(\partial \Omega)}^2\le
      n \ln \left(A \int_{\Omega} |\nabla_g f|^2 dV + B\int_{\partial \Omega}
       f^2 dS\right)  \\
       &=& n \ln \left(A \int_{\Omega} |\nabla_g f|^2 dV +B \right).\end{eqnarray*}
      Note that for any fixed $a>0$, the following elementary inequality holds:
      \begin{eqnarray*}  \ln x \le a x- 1 -\ln a \; \quad \mbox{for all}\;\;  x>0. \end{eqnarray*}
Therefore
\begin{eqnarray*}  \int_{\partial \Omega} f^2 \ln f^2 dS &\le & n  \ln \big(A \int_{\Omega} |\nabla_g f|^2 dV  +B)\\
  &\le &  na \big(A \int_\Omega |\nabla_g f|^2 dV +B\big)  -n (1+\ln a).\end{eqnarray*}
        Taking $\epsilon= (naA)/2$, we deduce   \begin{eqnarray*} \int_{\partial \Omega} f^2 \ln f^2 dS \le 2 \epsilon \int_\Omega
     |\nabla_g f|^2 dV - n \ln \epsilon + 2BA^{-1} \epsilon +n \ln \frac{nA}{2e}. \end{eqnarray*}
      Finally,  by putting $f=\frac{v}{\|v\|_{L^2(\partial \Omega)}}$  we get
      \begin{eqnarray*} \int_{\partial \Omega} \frac{v^2}{\|v\|_{L^2(\partial \Omega)}^2}\, \ln \frac{v^2}{\|v\|_{L^2(\partial \Omega)}^2} \,dS \le 2 \epsilon \int_\Omega
     \frac{|\nabla_g v|^2}{\|v\|_{L^2(\partial \Omega)}^2} dV - n \ln \epsilon + 2BA^{-1} \epsilon +n \ln \frac{nA}{2e}, \end{eqnarray*}
     so that \begin{eqnarray*} \int_{\partial \Omega}  v^2 \ln \frac{|v|}{\|v\|_{L^2(\partial \Omega)}} dS \le \epsilon \int_\Omega
     |\nabla_g v|^2dV + \left( \frac{n}{2}\, \ln \frac{nA}{2e} +BA^{-1} \epsilon -\frac{n}{2} \ln\, \epsilon\right)\|v\|^2_{L^2(\partial \Omega)}.\end{eqnarray*}
      This is just the desired  Log-Sobolev trace inequality (\ref{5=2}).

     \vskip 0.25 true cm

 (i) $\Longrightarrow$ (iii): \
Using H\"{o}lder interpolation inequality, we have \begin{eqnarray*}  \int_{\partial \Omega}
 v^2 dS &=&\int_{\partial \Omega} v^{2- \frac{2}{n+1}} v^{\frac{2}{n+1}} dS = \int_{\partial \Omega}
v^{\frac{2n}{n+1}} v^{\frac{2}{n+1}} dS \\
&\le & \left(\int_{\partial \Omega}  v^{\frac{2np}{n+1}} dS \right)^\frac{1}{p} \left( \int_{\partial \Omega} v^{\frac{2p'}{n+1}} dS\right)^{\frac{1}{p'}}.\end{eqnarray*}
 By choosing $p=\frac{n+1}{n-1}$,  $\, p'=\frac{n+1}{2}$, we obtain \begin{eqnarray*}
\int_{\partial \Omega} v^2 dS \le \left( \int_{\partial \Omega}  v^{\frac{2n}{n-1}} dS \right)^{\frac{n-1}{n+1}}
\left( \int_{\partial \Omega}  |v|dS\right)^{\frac{2}{n+1}}, \end{eqnarray*}
so  that
\begin{eqnarray*}  \| v\|_{L^2(\partial \Omega)}^{2+\frac{2}{n}} \le  \left(\int_{\partial \Omega} v^{\frac{2n}{n-1}} dS\right)^{\frac{n-1}{n}}
\left( \int_{\partial \Omega} |v|dS \right)^{\frac{2}{n}}.\end{eqnarray*}
   Applying the Sobolev trace inequality (\ref{5=1}), we have
 \begin{eqnarray*}  \| v\|_{L^2(\partial \Omega)}^{2+\frac{2}{n}} \le \left(A\|\nabla_g
v\|_{L^2 (\Omega)}^2 +B \| v\|_{L^2 (\partial \Omega)}^2 \right) \|v\|_{L^1(\partial \Omega)}^{\frac{2}{n}}, \end{eqnarray*}
which is the  Nash trace inequality (\ref{5=3}).

 \vskip 0.25 true cm

(ii)  $\Longrightarrow$ (iv): \ Let $u(t,x)$ be a smooth solution of the heat equation for the Dirichlet-to-Neumann operator.
Clearly, $u$ can be written as
\begin{eqnarray*} u(t, x)=\int_{\partial \Omega} {\mathcal{K}}(t, x, y) u(0, y) dS(y),\end{eqnarray*}
where ${\mathcal{K}}$ is the fundamental solution of (\ref{55-2}).
Note that \begin{eqnarray*}  \sup_{u\ne 0} \frac{\|u(t, \cdot)\|_{L^\infty(\partial \Omega)}}{\|u(0,\cdot)\|_{L^1(\partial \Omega)}} =
\sup_{x,y} {\mathcal{K}}(t,x, y).\end{eqnarray*}
For any fixed $t>0$, we can take  $p(s)= \frac{t}{t-s}$. It is obvious that $p(0)=1$ and $p(t)=\infty$.
Let $u=u(t,x)$ be a positive solution to the
heat equation (\ref{1-5}) of the Dirichlet-to-Neumann operator. By  a direct calculation, we have
\begin{eqnarray} &&\frac{\partial \|u\|_{L^{p(s)}(\partial \Omega)}}{\partial s} = - \|u\|_{L^{p(s)}(\partial \Omega)} \cdot \frac{p'(s)}{p^2(s)}
\cdot \ln \|u\|_{L^{p(s)}(\partial \Omega)}^{p(s)}
\\ &&\quad \quad \; + \frac{\|u\|_{L^{p(s)}(\partial \Omega)}^{1-p(s)}}{ p(s)} \left(p'(s) \int_{\partial \Omega} u^{p(s)} \ln u\, dS +p(s) \int_{\partial \Omega}
u^{p(s)-1} ({\mathcal{N}}_gu) dS\right).\nonumber\end{eqnarray}
 Multiplying by $p^2(s) \|u\|_{L^{p(s)}(\partial \Omega)}^{p(s)}$ on both sides, and applying the property of ${\mathcal{N}}_g$, we obtain
 \begin{eqnarray*} p(s)^2 \|u\|_{L^{p(s)}(\partial \Omega)}^{p(s)}\,  \frac{\partial \|u\|_{L^{p(s)}(\partial \Omega)}}{\partial s} &=& -p'(s)
\|u\|_{L^{p(s)}(\partial \Omega)}^{1+p(s)}
  \ln \|u\|_{L^{p(s)}(\partial \Omega)}^{p(s)} \\
  &&+ p(s) p'(s) \|u\|_{L^{p(s)}(\partial \Omega)} \int_{\partial \Omega} u^{p(s)} \ln u\, dS\\
  &&- p^2(s) (p(s) -1) \|u\|_{L^{p(s)}(\partial \Omega)} \int_\Omega u^{p(s) -2} |\nabla_g u|^2 dV. \end{eqnarray*}
 Dividing both sides by $\|u\|_{L^{p(s)}(\partial \Omega)}$, we find that
 \begin{eqnarray} \label {55-5}&& p^2(s) \|u\|_{L^{p(s)}(\partial \Omega)}^{p(s)} \frac{\partial \left(\ln \|u\|_{L^{p(s)}(\partial \Omega)}\right)}{\partial s}
 = -p'(s) \|u\|_{L^{p(s)}(\partial \Omega)}^{p(s)}
  \ln \|u\|_{L^{p(s)}(\partial \Omega)}^{p(s)} \\
  && \qquad\quad \quad \;\; +p(s) p'(s) \int_{\partial \Omega} u^{p(s)} \ln u\, dS- 4(p(s)-1) \int_{\Omega} |\nabla_g(u^{p(s)/2})|^2 dV.\nonumber\end{eqnarray}
  Putting \begin{eqnarray*} v=\frac{u^{p(s)/2} (s,x)}{\|u^{p(s)/2}\|_{L^2(\partial \Omega)}},\end{eqnarray*}
 we get  \begin{eqnarray*} \|v\|_{L^2(\partial \Omega)}=1  \;\; \mbox{and}\;\; \ln v^2 =\ln u^{p(s)} -\ln \|u\|_{L^{p(s)}(\partial \Omega)}^{p(s)},\end{eqnarray*}
 so that  \begin{eqnarray*} p'(s) \int_{\partial \Omega} v^2 \ln v^2\, dS
 &=&p'(s) \int_{\partial \Omega} \frac{u^{p(s)}}{\|u\|_{L^{p(s)}(\partial \Omega)}^{p(s)}} \left( \ln u^{p(s)}
- \ln \|u\|_{L^{p(s)}(\partial \Omega)}^{p(s)}\right) dS\\
 &=&\frac{p(s) p'(s) }{\|u\|_{L^{p(s)}(\partial \Omega)}^{p(s)}} \int_{\partial \Omega} u^{p(s)} \ln u\, dS
-p'(s) \ln \|u\|_{L^{p(s)}(\partial \Omega)}^{p(s)}.\end{eqnarray*}
   We substitute $v$ into the right-hand side of (\ref{55-5}) to obtain
   \begin{eqnarray*} p^2(s) \,\frac{\partial \big(\ln \|u\|_{L^{p(s)}(\partial \Omega)}\big)}{\partial s} = p'(s) \left( \int_{\partial \Omega}
    v^2 \ln v^2 dS -\frac{4(p(s)-1)}{p'(s)} \int_\Omega |\nabla_gv|^2 dV\right).\end{eqnarray*}
    Choose \begin{eqnarray*} 2\epsilon= \frac{4(p(s)-1)}{p'(s)} =\frac{4s(t-s)}{t}\le t. \end{eqnarray*}
  By the Log-Sobolev trace inequality, we get
  \begin{eqnarray*} p^2(s) \frac{\partial \big(\ln \|u\|_{L^{p(s)}(\partial \Omega)}\big)}{\partial s}
  \le p'(s) \left[-n\ln \big(\frac{2s(t-s)}{t}\big) +A^{-1} Bt +n\ln \frac{nA}{2e}\right].\end{eqnarray*}
 Since $\frac{p'(s)}{p^2(s)}=\frac{1}{t}$, we deduce
  \begin{eqnarray*} \frac{\partial \big(\ln \|u\|_{L^{p(s)}(\partial\Omega)}\big)} {\partial s}\le \frac{1}{t} \left[ -n \ln \frac{s(t-s)}{t}
  + A^{-1} B t+n \ln \frac{nA}{4e}\right].\end{eqnarray*}
Integrating the above inequality from $0$ to $t$, we have
\begin{eqnarray*}
\ln \frac{\|u(t,x)\|_{L^{p(t)}(\partial \Omega)}}{\|u(0, x)\|_{L^{p(0)}(\partial \Omega)}} \le -n \ln t +2n + A^{-1} B t +n\ln \frac{nA}{4e}.\end{eqnarray*}
Noting that $p(0)=1, p(t)=\infty$, we get \begin{eqnarray*} \|u(t,x)\|_{L^\infty(\partial \Omega)} \le
\left( \frac{\exp(A^{-1} B t+ 2n+n\ln \frac{nA}{4e})}{t^n}\right)\|u(0, x)\|_{L^1(\partial \Omega)},\end{eqnarray*} which implies
\begin{eqnarray*} {\mathcal{K}}(t,x,y)\le \left(\frac{nAe}{4}\right)^{n} \frac{e^{A^{-1} B t}}{t^n}.\end{eqnarray*}

\vskip 0.3 true cm

 (iii)   $\Longrightarrow$  (iv): \
Let ${\mathcal{K}}={\mathcal{K}}(t,x,y)$  be the heat kernel associated to the Dirichlet-to-Neumann operator. For any fixed $y\in \partial \Omega$, we define
$v=v(t,x)={\mathcal{K}}(t,x,y)$. It is clear that  $\int_{\partial \Omega} v(t,x)dS(x) =1$, i.e., $\|v\|_{L^2(\partial \Omega)}=1$. Since \begin{eqnarray*}
 \frac{\partial }{\partial t} \left(\int_{\partial \Omega} v^2 (t,x) dS\right) &=&\int_{\partial \Omega}
2v \, \frac{\partial v}{\partial t} dS =\int_{\partial \Omega} 2 v({\mathcal{N}}_gv) dS \\
&=& \int_{\partial \Omega} 2v \,\frac{\partial v}{\partial \nu} \, dS =-2\int_\Omega
|\nabla_g v|^2 dV. \end{eqnarray*}
 From the Nash trace inequality with $\|v\|_{L^2(\partial \Omega)} =1$, we have
\begin{eqnarray*} \|v\|_{L^2(\partial \Omega)}^{2+\frac{2}{n}} \le \big( A\| \nabla_g v\|_{L^2(\Omega)}^2 +B\|v\|_{L^2(\partial \Omega)}^2 \big) \end{eqnarray*}
 so that  \begin{eqnarray*}  -\|\nabla_g v\|_{L^2(\Omega)}^2 \le -\frac{1}{A} \|v\|_{L^2(\partial \Omega)}^{2+\frac{2}{n}
} +\frac{B}{A} \|v\|_{L^2(\partial \Omega)}^2.\end{eqnarray*}
  Therefore \begin{eqnarray*} \frac{\partial}{\partial t} \left(\int_{\partial \Omega} v^2 (t, x) dS\right)
\le -\alpha \|v\|_{L^2(\partial \Omega)}^{2+\frac{2}{n}} +\gamma \|v\|_{L^2(\partial \Omega)}^2,\end{eqnarray*} where $\alpha= \frac{2}{A}, \gamma= \frac{2B}{A}$.
Write \begin{eqnarray*} f(t)= \int_{\partial \Omega} v^2(t,x)dS, \quad \, g(t)= e^{-\gamma t} f(t),\end{eqnarray*} we have
\begin{eqnarray*} \frac{\partial f(t)}{\partial t} \le -\alpha [f(t)]^{1+\frac{1}{n}} +\gamma f(t)
\end{eqnarray*}
and  \begin{eqnarray*} \frac{\partial g(t)}{\partial t} &=&-\gamma e^{-\gamma t} f(t) + e^{-\gamma t} \,\frac{\partial f(t)}{\partial t} \\
  &\le &-\gamma e^{-\gamma t} f(t) + e^{-\gamma t} \gamma f(t)-\alpha e^{-\gamma t}\left[f(t)\right]^{1+\frac{1}{n}}.  \end{eqnarray*}
 Thus  \begin{eqnarray*} \frac{\partial g(s)}{\partial s} \le -\alpha e^{\frac{\gamma s}{n}} [g(s)]^{1+\frac{1}{n}},  \quad \, s\in (0, t].\end{eqnarray*}
Integrating the above inequality from $\frac{t}{2}$ to $t$, we get
\begin{eqnarray*} -n\left([g(s)]^{-\frac{1}{n}}\right)\bigg|_{s=\frac{t}{2}}^{s=t}\le -\frac{n\alpha}
{\gamma} \left[ e^{\frac{\gamma t}{n}}- e^{\frac{\gamma t}{2n}} \right],\end{eqnarray*}
so that \begin{eqnarray*} ne^{\frac{\gamma t}{n}} [f(t)]^{-\frac{1}{n}} \ge \frac{n\alpha}{\gamma}
\left(e^{\frac{\gamma t}{n}} - e^{\frac{\gamma t}{2n}} \right). \end{eqnarray*}
 Therefore \begin{eqnarray} \label{5--0.} e^{\frac{\gamma t}{2n} } [f(t)]^{-\frac{1}{n}} \ge \frac{\alpha}{\gamma}
\left(e^{\frac{\gamma t}{2n}}-1\right) \ge \frac{\alpha}{\gamma}\cdot \frac{\gamma t}{2n} =\frac{\alpha t}{2n},\end{eqnarray}
 i.e.,  \begin{eqnarray*} f(t)\le \frac{1}{t^n} \left(\frac{2n}{\alpha}\right)^n e^{\frac{\gamma t}{2}}.\end{eqnarray*}
This implies  \begin{eqnarray*} \int_{\partial \Omega} {\mathcal{K}}(t,x,y) {\mathcal{K}}(t,x,y) dS(x)\le \frac{\big(2n/\alpha)^n}{t^n} e^{\gamma t/2}.\end{eqnarray*}
Noticing that  ${\mathcal{K}}(t,x,y)$ is symmetric with respect to $x$ and $y$, we find by the reproducing property  that  for any $x\in \partial \Omega$,
\begin{eqnarray*} {\mathcal{K}}(2t, x,x)=\int_{\partial \Omega} {\mathcal{K}}(t,x,y) {\mathcal{K}}(t,y,x) dS(y) \le \frac{(2n/\alpha)^n}{t^n} e^{\gamma t/2}.\end{eqnarray*}
  It follows that \begin{eqnarray*} {\mathcal{K}}(t,x,y)& =&\int_{\partial \Omega} {\mathcal{K}}(t/2, x,z) {\mathcal{K}}(t/2, z,y) dS(z) \\
  &\le & \left(\int_{\partial \Omega} {\mathcal{K}}^2 (t/2, x, z)dS(z) \right)^{\frac{1}{2}}
  \left(\int_{\partial \Omega} {\mathcal{K}}^2(t/2, z, y) dS(z)\right)^{\frac{1}{2}}\\
&=&\big({\mathcal{K}}(t, x,x)\big)^{1/2} \big({\mathcal{K}}(t,y,y)\big)^{1/2} \le \frac{(4n/\alpha)^n}{t^n} e^{\gamma t/4},\end{eqnarray*}
i.e., \begin{eqnarray*} {\mathcal{K}}(t,x,y)
 \le \frac{(2nA)^n}{t^n} e^{A^{-1} Bt/2}.\end{eqnarray*}

When $t$ is large, the above bound can be improved. Indeed, we find by the first inequality of (\ref{5--0.}) that for $t\ge 1$,
\begin{eqnarray*} f(t) ={\mathcal{K}}(2t, x,x) \le \frac{e^{\gamma t/2}\left(\frac{\gamma}{\alpha}\right)^n}{(e^{\frac{\gamma t}{2n}}-1)^n}\le
\frac{B^n e^{A^{-1}Bt}}{(e^{B/(nA)}-1)^n}.\end{eqnarray*}
 Therefore, there exists a positive constant $C$ depending only on $A$ and $B$ such that
 \begin{eqnarray*} {\mathcal{K}}(t,x,y) \le \min\left\{ \frac{(2nA)^n}{t^n} e^{A^{-1}Bt/2}, C\right\}, \;\; t>0, \;\; x,y\in \partial \Omega .\end{eqnarray*}

\vskip 0.3 true cm

(iv) $\Longrightarrow$ (i): \ Let $H=H(t, x,y)$ be the shift heat kernel of the  Dirichlet-to-Neumann operator:
\begin{eqnarray*}  {\mathcal{N}}_gu  -c_2 u-\frac{\partial u}{\partial t}= 0 \quad \;
\mbox{on}\;\; \partial \Omega.\end{eqnarray*}  Since \begin{eqnarray*} {\mathcal{K}}(t, x,y) \le  \frac{c_1 }{t^{n}} e^{c_2t},\end{eqnarray*}  we have
\begin{eqnarray*} H=H(t,x, y)=e^{-c_2t} {\mathcal{K}}(t,x,y)\le \frac{c_1}{t^n}.\end{eqnarray*}
 Note that \begin{eqnarray*} \int_{\partial \Omega} H(t,x,y) dS(y) \le \int_{\partial \Omega}
 {\mathcal{K}}(t, x,y) dS(y) =1.\end{eqnarray*}
It follows that for all $f\in L^2(\partial \Omega)$,
\begin{eqnarray*} \|H * f\|_{L^\infty (\partial \Omega)}& =&\sup_{x\in \partial \Omega} \bigg|\int_{\partial \Omega} H(t,x,y) f(y) dS(y)\bigg| \le
\sup_{x\in \partial \Omega} \left( \int_{\partial \Omega} H^2(t,x,y)dS(y)\right)^{1/2} \|f\|_{L^2(\partial \Omega)}\\
&\le&\frac{\sqrt{c_1}}{t^{n/2}} \left(\int_{\partial \Omega} H(t,x,y) dS(y) \right)^{1/2}
\|f\|_{L^2(\partial \Omega)} \le \frac{\sqrt{c_1}}{t^{n/2}} \|f\|_{L^2(\partial \Omega)} .\end{eqnarray*}
Also, by H\"{o}lder's inequality we find that for all $q\in [1,n)$ and $q'=\frac{q}{q-1}$,
\begin{eqnarray} \label{55-7} \|H*f\|_{L^\infty(\partial \Omega)} &\le &  \sup_{x\in \partial \Omega}
\left( \int_{\partial \Omega} H^{q'} (t, x, y)dS(y)\right)^{1/{q'}} \|f\|_{L^q(\partial \Omega)} \\
&\le& \sup_{x\in \partial \Omega}\left( \int_{\partial \Omega} H^{q'/q} (t,x,y)
\cdot H(t,x,y) dS(y) \right)^{1/q'} \|f\|_{L^q(\partial \Omega)} \nonumber\\
&\le & \frac{c_1^{1/q}}{t^{n/q}}
\|f\|_{L^q(\partial \Omega)}.\nonumber\end{eqnarray}
For the operator  $L= \left(\sqrt{(-{\mathcal{N}}_g+c_2)}\right)^{-1}$, by eigenfunction expansion
 (Laplace transform) we get that for $f\in  C^\infty(\partial \Omega)$,
\begin{eqnarray*}  (Lf)(x) &=&\big(\Gamma (1/2)\big)^{-1} \int_0^\infty
t^{-1/2} \big[e^{({\mathcal{N}}_g-c_2)t} f\big] (t, x) dt \\
&=&\big(\Gamma (1/2)\big)^{-1}\int_0^\infty t^{-\frac{1}{2}} (H*f) (t,x) dt,\end{eqnarray*}
where  $e^{({\mathcal{N}}_g-c_2)t} f$ is the semigroup notation for $H*f$.  For fixed $t>0$, we write $Lf=L_1f  +
 L_2f,$
   where  \begin{eqnarray*} L_1 f(x)= \Gamma(1/2)^{-1} \int_0^t s^{-\frac{1}{2}} [H*f](s,x) ds,\\
  L_2 f(x)= \Gamma(1/2)^{-1} \int_t^\infty  s^{-\frac{1}{2}} [H*f](s,x) ds.\end{eqnarray*}
Note that for any $\lambda>0$,
\begin{eqnarray*} \label{55-9} |\{x\big||Lf(x)|\ge \lambda \} |\le |\{x\big|| L_1 f(x)|\ge \frac{\lambda}{2}\} |+
|\{ x\big| |L_2 f(x)|>\frac{\lambda}{2}\}|.\end{eqnarray*} Form (\ref{55-7}) and the definition of $L_2f$,
 \begin{eqnarray*} \|L_2 f\|_{\infty} \le \frac{1}{\sqrt{\pi}}\,  c_1^{1/q} \int_t^\infty s^{-\frac{1}{2}- \frac{n}{q}}\|f\|_{L^q(\partial \Omega)} ds =\frac{2q}{(2n-q)\sqrt{\pi}} c_1^{\frac{1}{q}} t^{\frac{1}{2}-\frac{n}{q}}\|f\|_{L^q(\partial \Omega)}.\end{eqnarray*}
We choose $t$ such that \begin{eqnarray} \label{55-10} \frac{\lambda}{2} =\frac{2q}{(2n-q)\sqrt{\pi}}\,c_1^{\frac{1}{q}}
 t^{\frac{1}{2}-\frac{n}{q}} \|f\|_{L^q(\partial \Omega)}. \end{eqnarray}
It follows that \begin{eqnarray*} |\{ x\big||Lf(x)|\ge \lambda\} |\le |\{x\big| |L_1 f(x)|\ge \frac{\lambda}{2}\}|,
\end{eqnarray*} so that \begin{eqnarray*}
| \{ x\big| |Lf(x)|\ge \lambda \} | \le |\{x\big| |L_1f(x)|\ge \frac{\lambda}{2}\}|\le \left(\frac{\lambda}{2}\right)^{-q}
\int_{\partial \Omega} |L_1 f(x)|^q dS(x).\end{eqnarray*}
 From Minkowski's inequality and Young's inequality, we get
\begin{eqnarray*} \|L_1 f\|_{L^q(\partial \Omega)} &\le &\big(\Gamma (1/2)\big)^{-1}
\int_0^t s^{-\frac{1}{2}}\| H* f(s, \cdot) \|_{L^q(\partial \Omega)} ds \\
&\le & \big(\Gamma (1/2)\big)^{-1}  \int_0^t s^{-\frac{1}{2}} \left(\sup_x  \| H(s,x, \cdot)
\|_{L^1(\partial \Omega)} \right)\|f\|_{L^q(\partial \Omega)} ds \\
&\le& \frac{2}{\sqrt{\pi}}\, t^{\frac{1}{2}} \|f\|_{L^q(\partial \Omega)},\end{eqnarray*}
which shows \begin{eqnarray*} |\{x\big| Lf(x) |\ge \lambda \} | \le \left(\frac{2}{\sqrt{\pi}}\right)^{q} \left(\frac{\lambda}{2}\right)^{-q}
t^{\frac{q}{2}} \|f\|_{L^q(\partial \Omega)}^q.\end{eqnarray*}
According to the choice of $t$ in (\ref{55-10}), this is equivalent to
\begin{eqnarray*} \{x\big| |Lf (x)|\ge \lambda \}| \le \left(\frac{16}{\pi}\right)^{nq/(2n-q)} \left(\frac{q}{2n-q}\right)^{q^2/(2n-q)} (c_1)^{\frac{q}{2n-q}} \lambda^{-r}
  \|f\|_{L^q(\partial \Omega)}^r, \end{eqnarray*}
  where  $r=(2qn)/(2n-q)$.  By the Marcinkiewicz interpolation lemma, we know that $L$ is a bounded
  operator from $L^2(\partial \Omega)$ to $L^p(\partial \Omega)$ with $p=\frac{2n}{n-1}$ (taking $q=2$), i.e.,
   \begin{eqnarray} \label{55-11} \|Lu\|_{L^p(\partial \Omega)} \le c(c_1)^{\frac{1}{n-1}} \|u\|_{L^2(\partial \Omega)}
   \quad \,\mbox{for all}\;\; u\in C^\infty (\partial \Omega).\end{eqnarray}
Write $v=Lu$. Then $u=L^{-1}v$ and \begin{eqnarray*}  \|u\|_{L^2(\partial \Omega)}^2
&=&\langle L^{-1} v, L^{-1} v\rangle =\langle L^{-2} v, v\rangle \\
&=& \langle -{\mathcal{N}}_gv +c_2 v, v\rangle =\int_{\Omega} |\nabla_g v|^2 dV + c_2 \int_{\partial \Omega} v^2 ds.\end{eqnarray*}
Substitute this into (\ref{55-11}) we arrive at the Sobolev trace inequality:
\begin{eqnarray*} \|v\|_{L^2(\partial \Omega)}^p \le const. (c_1)^{\frac{1}{n-1}}
\left( \|\nabla_g v\|_{L^2(\Omega)}^2 +c_2\|v\|_{L^2(\partial \Omega)}^2\right).\end{eqnarray*}
 If we take $c_1 =\left(\frac{nAe}{4}\right)^{n-1}$ and $c_2=A^{-1}B$ as in (\ref{5=4}), then we get a  Sobolev trace inequality with the claimed constants:
 \begin{eqnarray*} \|v\|_{L^p(\partial \Omega)}^2 \le const. A \|\nabla_g v\|_{L^2(\Omega)}^2  +const. B \|v\|_{L^2(\partial \Omega)}^2.\end{eqnarray*}

\vskip 0.3 true cm

(v) $\Longrightarrow$ (iv): \  It is obvious because $t\in[0, 1]$ and $d(x,y)$ is bounded for all $x,y\in \partial \Omega$.

\vskip 0.3 true cm

(iv) $\Longrightarrow$ (v): \ By (iv) there exists a constant $c_1>0$ such that
\begin{eqnarray} \label{3-30}  {\mathcal{K}}(t,x,y) \le c_1 t^{-n} \quad \,\mbox{for all}\;\; t\in[0,1]\;\; \mbox{and}\;\;
x,y\in \partial \Omega.\end{eqnarray}
Clearly, $\partial \Omega$  has at most a finite number of connected components, say $\Lambda_1,\cdots, \Lambda_N$, with $\Lambda_i\ne \Lambda_j$ if $i\ne j$.
Let \begin{eqnarray*} \mathfrak{T}=\{ \psi\in C^\infty (\partial \Omega) \big| \max_{\underset{x_i\in \Lambda_i,x_j\in \Lambda_j}{i,j\in \{1,\cdots,N\}}} |\psi(x_i)-\psi(x_j)|+D\|\nabla \psi\|_{\infty} \le D\}, \end{eqnarray*}
where $D=1+\sum_{i=1}^N \mbox{diam}\, \Lambda_i$. We define the derivative $\delta_\psi$ on the space $\mathcal{L}(L^p(\partial \Omega))$ of all the linear operators on $L^p(\partial \Omega)$ by
$\delta_\psi (E) =[M_\psi, E]$, where $M_\psi$ denotes the multiplication operator with the function $\psi$.  Furthermore, we can define $\delta_\psi^{m}(E)=\delta_\psi (\delta_\psi^{m-1}(E))$ for all $m>1$ by induction.
It follows from (\ref{3-30}) and Proposition 4.1 of \cite{EO} that there exists a constants $c_2>0$ such that
\begin{eqnarray} \label{3.17a} \|\delta_\psi^{n+1} (T_t)\|_{1\to\infty} \le c_2 t \quad \; \mbox{for all}\;\; \psi\in\mathfrak{T}\,\, \mbox{and}\;\;
 t\in [0,1].\end{eqnarray}
 Note that \begin{eqnarray*}
  \delta_\psi^{n+1} (T_t) = [M_\psi, [\cdots, [M_\psi, T_t]\cdots]]= \sum_{k=0}^{n+1}(-1)^k \begin{pmatrix}k \\ n+1 \end{pmatrix} \psi^{n+1-k} T_t \psi^k .\end{eqnarray*} It follows from (\ref{3.17a}) that for any fixed $\psi\in\mathfrak{T}\,\, \mbox{and}\;\;
 t\in [0,1]$, \begin{eqnarray*}
  &&- c_2t \int_{\partial \Omega} v(y) \, dS(y)\le  \int_{\partial \Omega}  \sum_{k=0}^{n+1} (-1)^k \begin{pmatrix}k \\ n+1 \end{pmatrix} \psi^{n+1-k}(x)\, \mathcal{K} (t,x,y) \,\psi^k(y) \,v (y) \, dS(y)\\
  &&\qquad \qquad  \le  c_2t \int_{\partial \Omega} v(y) \, dS(y)\quad \, \mbox{for all}\;\; 0\le v\in C^\infty(\partial \Omega). \end{eqnarray*}
  This implies \begin{eqnarray*} |(\psi(x)-\psi(y))^{n+1} {\mathcal{K}}(t,x,y) |\le c_2 t\quad \, \mbox{for all}\;\; t\in[0,1],  x,y\in \partial \Omega\,\, \mbox{and}\;\; \psi\in \mathfrak{T}.\end{eqnarray*}   Since $d(x,y)=\sup \{|\psi(x)-\psi(y)|\big| \psi\in \mathfrak{T} \}$, we get \begin{eqnarray*} (d(x,y))^{n+1} {\mathcal{K}}(t,x,y) \le c_2t, \end{eqnarray*}
 i.e.,    \begin{eqnarray*} \left( \frac{d(x,y)}{t}\right)^{n+1} {\mathcal{K}}(t,x,y) \le c_2t^{-n} \quad\, \mbox{for all}\;\; t\in[0,1],  x,y\in \partial \Omega\,\, \mbox{and}\;\; \psi\in \mathfrak{T}.\end{eqnarray*}
  Combining this and (\ref{3-30}) we have
    \begin{eqnarray*} \left( 1+ \frac{d(x,y)}{t}\right)^{n+1} {\mathcal{K}}(t,x,y) \le 2^{n+1} (c_1+c_2) t^{-n} \quad\, \mbox{for all}\;\; t\in[0,1],  x,y\in \partial \Omega\,\, \mbox{and}\;\; \psi\in \mathfrak{T}.\end{eqnarray*}
This is equivalent to
\begin{eqnarray*}  {\mathcal{K}}(t,x,y) \le c t^{-n} \left( 1+\frac{d^2(x,y)}{t^2}\right)^{-(n+1)/2} \quad\, \mbox{for all}\;\; t\in[0,1],  x,y\in \partial \Omega\,\, \mbox{and}\;\; \psi\in \mathfrak{T}, \quad \;\,\end{eqnarray*}
which just is the inequality (\ref{08.08.1}).

\vskip 0.2 true cm

(vi)$\Longrightarrow$(i): \  If the inequality (\ref{5=5a}) holds for all $q\in C^\infty (\partial \Omega)$, then
$I_q(0)\le (1-\epsilon)$ for all $q\in C^\infty (\partial \Omega)$ with $\|q\|_{L^n (\partial \Omega)}^n \le \frac{1-\epsilon}{C(n)}$, where
$0<\epsilon<1$. This implies $I_q(0)=0$ for all such $q$, and consequently  $(-A{\mathcal{N}}_g +B+ q)$ is a nonnegative operator for all $q\in C^\infty (\partial \Omega)$ with $\|q\|_{L^n (\partial \Omega)}^n \le \frac{1-\epsilon}{C(n)}$. It follows that \begin{eqnarray*} 0\le \langle (-A{\mathcal{N}}_g +B+q)f, f \rangle = A\int_{\Omega} |\nabla_g u|^2 dV + B\int_{\partial \Omega} f^2 dS+\int_{\partial \Omega}
 q|f|^2 dS\end{eqnarray*}
 for all such $q$ and all $0\le f\in C^\infty (\partial \Omega)$,
 where $u\in C^\infty (\bar \Omega)$ satisfies
     $\Delta_g u=0$ in $ \Omega$ and $u=f$ on $\partial \Omega$, i.e., \begin{eqnarray*}  \sup_{\underset {-q\in C^\infty (\partial \Omega)}{\|q\|_{L^n (\partial \Omega) \le ((1-\epsilon)/C(n))^{1/n}}}} \bigg\{ \int_{\partial \Omega}- q|f|^2 dS \bigg\} \le A \int_{\Omega} |\nabla_g u|^2 dV +B\int_{\partial \Omega} f^2dS.\end{eqnarray*}
     As the dual of $L^{n} (\partial \Omega)$ is $L^{n/(n-1)}(\partial \Omega)$, this yields
    \begin{eqnarray*} \bigg(\int_{\partial \Omega} |f|^{2n/(n-1)} dx\bigg)^{(n-1)/n} \le \left(\frac{C(n)}{1-\epsilon}\right)^{1/n} \left(A\int_{\Omega} |\nabla_g u  |^2 dV+ B\int_{\partial \Omega} f^2 dS\right) \quad \; \mbox{for all}\;\; f\in C^\infty (\partial \Omega).\end{eqnarray*} By virtue of
\begin{eqnarray*}  \inf_{\underset {v=f\;\;on\;\;\partial \Omega} {v\in W^{1,2}(\Omega)}} \int_{\partial \Omega} |\nabla_g v|^2 dV =\int_{\Omega}
  |\nabla_g u|^2 dV, \end{eqnarray*}
    we immediately get that for all $v\in W^{1,2}(\Omega)$, the Sobolev  trace inequality (\ref{5=1}) holds with new constant constants $A'=A\left(\frac{C(n)}{1-\epsilon}\right)^{1/n},\; B'=B\left(\frac{C(n)}{1-\epsilon}\right)^{1/n}$.

\vskip 0.2 true cm

(i)$\Longrightarrow$(vi): \   By monotonicity of $I_q(0)$ with respect to the $q(x)$, we may assume $q(x)< 0$ for all $x\in \partial \Omega$.
 Then the number $I_q(0)$ of non-positive eigenvalues for (\ref{5=6a}) is equal to the number of eigenvalues less than or equal to $1$ for the problem
 \begin{eqnarray} \label{08-1.1} (A{\mathcal{N}}_g -B)\psi (x) =\mu q(x) \psi(x) \quad \, \mbox{for all}\;\; x\in \partial \Omega,\end{eqnarray}
 where $A$ and $B$ are  positive  constants in the Sobolev trace inequality (\ref{5=1}). In fact, since \begin{eqnarray*}
 \frac{ \int_{\Omega} A|\nabla_g {\tilde\psi}|^2 dV +(B+q)\int_{\partial \Omega} \psi^2dS}{\int_{\partial \Omega} \psi^2 dS} = \frac{\int_{\partial \Omega} |q|\psi^2 dS}{ \int_{\partial \Omega}\psi^2 dS} \bigg[ \frac{A\int_{\Omega} |\nabla_g {\tilde \psi}|^2 dV+B\int_{\partial \Omega} \psi^2 dS}{\int_{\partial \Omega} |q|\psi^2dS} -1\bigg],\end{eqnarray*}
we see  that the dimension of the subspace on which the left-hand side is non-positive is equal to the dimension on which the quadratic form
$\big(A\int_{\Omega} |\nabla_g {\tilde \psi}|^2dV+B\int_{\partial \Omega} \psi^2 dS\big)/\big(\int_{\partial \Omega} |q| \psi^2 dS\big)$ is less than or equal to $1$, where $\tilde \psi$ is the harmonic extension of $\psi$  to $\bar \Omega$.
Let $\{ \psi_i (x) \}_{i=1}^\infty$ be a set of orthonormal eigenfunctions satisfying \begin{eqnarray*} (A{\mathcal{N}}_g -B)\psi_i =\mu_i q \psi_i\quad \mbox{on}\;\; \partial \Omega\end{eqnarray*} with the  eigenvalues $\{\mu_i\}$. It follows from \cite{Sa} (or \cite{Liu2}) that
$\mu_i> 0$ for all $i\ge 1$. Then the kernel of the ``heat'' equation $$\left(\frac{-A{\mathcal{N}}_g +B}{q} -\frac{\partial }{\partial t}\right)u=0 \quad\, \mbox{in}\;\; [0,\infty)\times (\partial \Omega)$$ must take the form
\begin{eqnarray} 0<H(t,x,y)= \sum_{i=1}^\infty e^{-\mu_i t} \psi_i(x) \psi_i(y)\quad \; \mbox{on}\;\; (0,\infty) \times (\partial \Omega) \times (\partial \Omega),\end{eqnarray}
where  the $L^2$-norm is given by the volume form $-q(x) dS(x)$ instead of $dS(x)$.
Consider the function \begin{eqnarray} \label{08-01.7} h(t) =\sum_{i=1}^\infty e^{-2\mu_i t} = \int_{\partial \Omega} \int_{\partial \Omega} H^2(t,x,y) q(x)q(y)dS(x)\,dS(y).\end{eqnarray}
In view of $ \bigg(\frac{A{\mathcal{N}}_g -B}{(-q(y))} -\frac{\partial }{\partial t}\bigg)H(t,x,y) \equiv 0$, we have
\begin{eqnarray} \label{08-01.2}\quad \;\quad\quad\;  \frac{\partial h}{\partial t}&=& 2 \int_{\partial \Omega} \int_{\partial \Omega}  H(t, x,y) q(x)q(y)
 \frac{\partial H(t, x,y)}{\partial t} dS(x)\, dS(y) \\  &=&  2 \int_{\partial \Omega} \int_{\partial \Omega}  H(t, x,y)\big[\big(-A{\mathcal{N}}_g +B\big)  H(t, x,y)\big] q(x)dS(x)\, dS(y).\nonumber\\
 &=& -2\int_{\partial \Omega} (-q(x)) \bigg(A \int_{\Omega} |\nabla_y {\tilde H}(t,x,y)|^2 dV(y)  +B \int_{\partial \Omega}  H^2 (t,x,y) dS(y) \bigg) dS(x). \nonumber\end{eqnarray}
 Here, for any fixed $t>0$ and $x\in \partial \Omega$, ${\tilde H}(t,x,y)$ satisfies \begin{eqnarray*} \left\{ \begin{array}{ll} \Delta_y {\tilde{H}}(t,x,y) =0 \quad \; \mbox{for}\;\; y\in \Omega,\\
  {\tilde H} (t,x,y) =H(t,x,y) \quad \; \mbox{for} \;\; y\in\partial \Omega.\end{array}\right. \end{eqnarray*}
 Similar to the method of \cite{LY2}, we get
 \begin{eqnarray} \label{08-01.4}&& \quad \;h(t)= \int_{\partial \Omega} (-q(x)) \int_{\partial \Omega} H^2(t,x,y) (-q(y)) dS(y)\, dS(x) \\
 &\le& \int_{\partial \Omega} (-q(x)) \left[\left( \int_{\partial \Omega} H^{\frac{2n}{n-1}}(t,x,y) dS(y) \right)^{\frac{n-1}{n+1}}
 \left( \int_{\partial \Omega} H(t,x,y) (-q(y))^{\frac{n+1}{2}}  dS(y)\right)^{\frac{2}{n+1}}\right] dS(x) \nonumber\\
  &\le & \bigg[ \int_{\partial \Omega}  (-q(x)) \bigg( \int_{\partial \Omega} H^{\frac{2n}{n-1} }(t,x,y)dS(y) \bigg)^{\frac{n-1}{n}} dS(x)\bigg]^{\frac{n}{n+1}} \nonumber\\
  && \times \bigg[ \int_{\partial \Omega}  (-q(x)) \bigg( \int_{\partial \Omega} H (t,x,y) (-q(y))^{\frac{n+1}{2}}  dS(y)\bigg)^2 dS(x) \bigg]^{\frac{1}{n+1}}.\nonumber\end{eqnarray}
 Set \begin{eqnarray*} Q(t,x) =\int_{\partial \Omega} H(t, x,y) (-q(y))^{\frac{n+1}{2}} dS(y)\end{eqnarray*}
 It is obvious that  \begin{eqnarray*} \bigg( \frac{-A{\mathcal{N}}_g +B}{q(x)} -\frac{\partial }{\partial t}\bigg) Q(t,x)\equiv 0\end{eqnarray*}
  and $$Q(0,x) =\int_{\partial \Omega} \delta(x-y) (-q(y))^{\frac{n-1}{2}} (-q(y))dS(y) =(-q(x))^{\frac{n-1}{2}}.$$
Similar to the argument of (\ref{08-01.2}), we have
 $\frac{\partial}{\partial t} \int_{\partial \Omega} Q^2 (t,x) (-q(x)) dS(x) \le 0$, so that
 \begin{eqnarray*} \int_{\partial \Omega} Q^2(t,x) (-q(x)) dS(x) \le \int_{\partial \Omega} Q^2(0,x) (-q(x)) dS(x) = \int_{\partial \Omega}  (-q(x))^{n} dS(x).\end{eqnarray*}
 By this and (\ref{08-01.4}),  we obtain
  \begin{eqnarray}\label{08-01.5} && h^\frac{n+1}{n} (t) \bigg(\int_{\partial \Omega} (-q(x))^{n}  dS(x) \bigg)^{-1/n} \\
 && \quad \quad \;\;  \le \int_{\partial \Omega}  (-q(x)) \bigg(\int_{\partial \Omega}  H^{\frac{2n}{n-1}} (t,x,y) dS(y) \bigg)^{\frac{n-1}{n}}dS(x) .\nonumber\end{eqnarray}
  Combining  (\ref{08-01.2}),  (\ref{08-01.5}) and the Sobolev trace inequality (\ref{5=1}), we get \begin{eqnarray*} \frac{\partial h}{\partial t} \le  -2  \bigg(\int_{\partial \Omega}  (-q(x))^{n}  dS(x) \bigg)^{-1/n} h^{\frac{n+1}{n}}(t).\end{eqnarray*}
 Thus
 \begin{eqnarray}\label{08-01.6}
  h(t) \le \left(\frac{n}{2}\right)^n\bigg( \int_{\partial \Omega} (-q(x))^{n}dS(x) \bigg) t^{-n}. \end{eqnarray}
  From (\ref{08-01.7}) we obtain
  \begin{eqnarray*} \left(\frac{n}{2}\right)^{n} \bigg(\int_{\partial \Omega} (-q(x))^{n}  dS(x) \bigg) t^{-n} \ge \sum_{i=1}^\infty e^{-2\mu_i t}.\end{eqnarray*}
  Finally, we set $\mu_k$ to be the greatest eigenvalue less than or equal to $1$.
     By setting $t=\frac{n}{2\mu_k}$,
  we get  \begin{eqnarray} \label{08-01.8}  \left( \int_{\partial \Omega} |q(x)|^{n} dS(x) \right) \mu_k^{n}&=&\left(\frac{n}{2}\right)^{n}
    \bigg(\int_{\partial \Omega} (-q(x))^{n}dS(x)\bigg) \left(\frac{n}{2}\right)^{-n} \mu_k^{n} \\
    &\ge& \sum_{i=1}^\infty e^{-n\mu_i/\mu_k}
 \ge  ke^{-n},\nonumber\end{eqnarray}
 which leads to  \begin{eqnarray} \label{08-01.9} \int_{\partial \Omega} |q(x)|^{n} dS(x) &\ge& \mu_k^{n} \int_{\partial \Omega} |q(x)|^{n} dS(x)\\
 &\ge & k e^{-n} = I_q(0) e^{-n}.\nonumber\end{eqnarray}
   Therefore, the Rozenblum-Lieb-Cwikel  type inequality holds with $C(n)= e^n$.
$\quad \quad \square $

\vskip 1.49 true cm

\section{Symbol expression of the Dirichlet-to-Neumann operator }

\vskip 0.45 true cm

Let $\Omega$ be a bounded domain $\Omega\subset  (\mathcal{M},g)$ with smooth boundary $\partial \Omega$. Then $\partial \Omega$ has an induced Riemannian metric $h$, and $\partial \Omega\hookrightarrow \bar \Omega$ has a second fundamental form, with associated Weingarten map
  \begin{eqnarray*} A_\nu :T_x (\partial \Omega) \to  T_x (\partial \Omega).\end{eqnarray*} Let $A^*_{\nu}: T^*_x (\partial \Omega)\to T_x^* (\partial \Omega)$ be the adjoint of $A_\nu$, and let $\langle \cdot, \cdot\rangle$
  be the  inner product on $T^*_x (\partial \Omega)$ arising from the given Riemannian metric.
As shown in  \cite{LU} and \cite{CNS}, ${\mathcal{N}}_g$ is a negative-semidefinite, self-adjoint, elliptic pseudodifferential operator in $OPS^0(\partial \Omega)$ with the principal symbol $-\left(\sum_{j,k=1}^n h^{jk} \xi_j\xi_k\right)^{1/2}$, where $(h^{jk})$ is the inverse of $h$.
    It was proved (see, Proposition C.1 of \cite{Ta2} ) that   \begin{eqnarray*} \label {b.33}{\mathcal{N}}_g= -\sqrt{-\Delta_h} +B \quad \;\mbox{mod}\;\, OPS^{-1}(\partial \Omega),\end{eqnarray*}
  where $B\in OPS^{0} (\partial \Omega)$  has principal symbol
    \begin{eqnarray*} \label{b.34} p_0(x, \xi)= \frac{1}{2} \left(Tr\, A_\nu  -\frac{\langle A^*_\nu \xi, \xi\rangle}{\langle \xi, \xi\rangle}\right).\end{eqnarray*}
 The following two Lemmas provide more information for the operators $-\sqrt{-\Delta_h}$ and $B$.

\vskip 0.29 true cm

\noindent  {\bf Lemma 4.1.} \   {\it Let $p(x, \xi)$ be the symbol of the operator $-\sqrt{-\Delta_h}$ on $\partial \Omega$, and let $p(x, \xi)\sim \sum_{j\ge 0} p_{1-j} (x, \xi)$.  Then
\begin{eqnarray}  \label{-5b.1} && p_1(x, \xi)=-\left(\sum_{j,k=1}^n h^{jk}(x) \xi_j\xi_k\right)^{1/2},\\
&& \label{-5b.1'} p_0(x, \xi)= \frac{i}{2} \left(\sum_{j,k=1}^n h^{jk}(x) \xi_j \xi_k\right)^{-1/2}
 \left[ \sum_{j,k=1}^n \frac{1}{\sqrt{|h(x)|}}\,\frac{\partial (\sqrt{|h(x)|} \, h^{jk}(x))}{\partial x_k} \, \xi_j\right.\\
   &&\left. \quad \qquad \quad \;\; - \frac{1}{2} \left(\sum_{j,k=1}^n h^{jk}(x) \xi_j \xi_k\right)^{-1} \left(\sum_{j,k,l,m=1}^n
 h^{lm}(x)\, \frac{\partial h^{jk}(x)}{\partial x_l} \xi_j\xi_k \xi_m\right) \right]\nonumber\end{eqnarray}
 and \begin{eqnarray}  \label{-5b.2'} &&  p_{-1} (x,\xi)= \frac{1}{2p_1} \bigg[ -p_0^2 +i \sum_{l=1}^n \bigg(\frac{\partial p_1}{\partial \xi_l}\, \frac{\partial p_0}{\partial x_l} + \frac{\partial p_0}{\partial \xi_l}\, \frac{\partial p_1}{\partial x_l} \bigg) +\frac{1}{2} \sum_{j,k=1}^n \frac{\partial^2 p_1}{\partial \xi_j \partial \xi_k} \, \frac{\partial^2 p_1}{\partial x_j \partial x_k} \bigg].\end{eqnarray}
}

\vskip 0.16 true cm

\noindent  {\it Proof.} \  Denote by $c(x, \xi)$ the symbol of $-\Delta_h$ on $\partial \Omega$. It follows from the composition formula (see, for example, (3.17) of p.$\,$13 of \cite{Ta2}) that
\begin{eqnarray*} c(x, \xi)&\sim& \sum_{\alpha \ge 0} \frac{i^{|\alpha|}}{\alpha!}D^\alpha_\xi p(x, \xi) D_x^\alpha p(x, \xi)\\
 &=& \left(p_1 +p_0 +p_{-1}+r\right)^2 - i \sum_{l=1}^n \frac{\partial (p_1+p_0 +p_{-1} +r)}{\partial \xi_l} \, \frac{\partial (p_1+p_0+p_{-1}+ r)}{\partial x_l} \\
 && -\frac{1}{2} \, \sum_{j,k=1}^n \frac{\partial^2(p_1 +p_{0} +p_{-1}+ r)}{\partial \xi_j \partial \xi_k} \, \frac{\partial^2(p_1 +p_{0} +p_{-1}+ r)}{\partial x_j \partial x_k}
  +\cdots, \end{eqnarray*}
where $r=p-p_1-p_0-p_{-1}\in S_{1,0}^{-2}$,  $D^\alpha =D^{\alpha_1}_1 \cdots D^{\alpha_n}_n$ and $D_l =\frac{1}{i}\, \frac{\partial}{\partial x_l}$.
On the other hand, it is well known that \begin{eqnarray*}c(x, \xi)= \sum_{j,k=1}^n h^{jk}(x) \xi_j\xi_k  -i \sum_{j,k=1}^n \frac{1}{\sqrt{|h(x)|}}\,\frac{\partial (\sqrt{|h(x)|} \, h^{jk}(x))}{\partial x_j} \, \xi_k.\end{eqnarray*}
 Then the symbol equations in $\xi$ of degree $2$, $1$ and $0$, respectively, give  \begin{eqnarray*}&&  p_1^2(x,\xi)  =\sum_{j,k=1}^n h^{jk}(x) \xi_j \xi_k,\\
  &&  2p_0(x,\xi)p_1(x,\xi) - i\sum_{l=1}^n
\frac{\partial p_1(x, \xi)}{\partial \xi_l}\, \frac{\partial p_1(x, \xi)}{\partial x_l} = -i \sum_{j,k=1}^n \frac{1}{\sqrt{|h(x)|}}\,  \frac{\partial (\sqrt{|h(x)|} \, h^{jk}(x))}{\partial x_j} \, \xi_k,\\
&& p_0^2 +2 p_1 p_{-1} -i \sum_{l=1}^n \bigg(\frac{\partial p_1}{\partial \xi_l}\, \frac{\partial p_0}{\partial x_l} +\frac{\partial p_0}{\partial \xi_l}\,\frac{\partial p_1}{\partial x_l} \bigg) -\frac{1}{2} \sum_{j,k=1}^n \frac{\partial^2 p_1}{\partial \xi_j\partial \xi_k} \, \frac{\partial^2p_1}{\partial x_j\partial x_k}=0,
\end{eqnarray*}
which implies \begin{eqnarray*} \label{+4.4} && p_1(x, \xi)=\pm\bigg(\sum_{j,k=1}^n h^{jk}(x) \xi_j\xi_k\bigg)^{1/2},\\
\label{+4.5}&& p_0(x, \xi)= \frac{i}{2p_1} \bigg[ \sum_{l=1}^n \frac{\partial p_1(x,\xi)}{\partial \xi_l} \, \frac{\partial p_1(x,\xi)}{\partial x_l}-
\sum_{j,k=1}^n \frac{1}{\sqrt{|h(x)|}}\,\frac{\partial (\sqrt{|h(x)|}\, h^{jk} (x))}{\partial x_j} \, \xi_k\bigg],\end{eqnarray*}
 \begin{eqnarray*} &&p_{-1} (x,\xi)= \frac{1}{2p_1} \bigg[ -p_0^2 +i \sum_{l=1}^n \bigg(\frac{\partial p_1}{\partial \xi_l}\, \frac{\partial p_0}{\partial x_l} + \frac{\partial p_0}{\partial \xi_l}\, \frac{\partial p_1}{\partial x_l} \bigg) +\frac{1}{2} \sum_{j,k=1}^n \frac{\partial^2 p_1}{\partial \xi_j \partial \xi_k} \, \frac{\partial^2 p_1}{\partial x_j \partial x_k} \bigg].\end{eqnarray*}
Noticing that the principal symbol of $-\sqrt{-\Delta_h}$ is non-positive, we immediately obtain (\ref{-5b.1}), (\ref{-5b.1'}) and (\ref{-5b.2'}).
$\quad \quad \square$

\vskip 0.32 true cm

Denote by $(A_\nu)_{jk}$ (respectively, $(A^*_\nu)_{jk}$) the entry in the $j$-th row and $k$-th column of $A_\nu$ (respectively, $A_\nu^*$) under a basis of $T_x(\partial \Omega)$ (respectively, $T^*_{x} (\partial \Omega)$), and let $C$ be an $n\times n$ matrix satisfying $C^TC=(h^{jk})$.
Let $2Q_1$ and $6Q_2$ be the principal minor determinants  of orders $2$ and $3$ for the matrix $\big((-2A_\nu)_{jk}\big)$, respectively.  We denote  $$Q_3= \sum_{j\ne k} \left[ (-2A_\nu)_{jj}  \big(-2{\tilde R}_{k(n+1)k(n+1)} +2(A_\nu^2)_{kk}\big) -(-2A_\nu)_{jk} \big(-2{\tilde R}_{k(n+1)j(n+1)} +2(A_\nu^2)_{kj}\big)\right],$$ where ${\tilde R}_{j(n+1)k(n+1)}$ is the curvature tensor with respect to $g$, and $e_{n+1}=\nu$.

\vskip 0.30 true cm

 \noindent  {\bf Lemma 4.2.} \ \  {\it   Let $\Omega$ be a bounded domain with smooth boundary $\partial \Omega$ in an $(n+1)$-dimensional Riemannian manifold $(\mathcal{M}, g)$.  Then the Dirichlet-to-Neumann operator ${\mathcal{N}}_g$ is given by
  \begin{eqnarray} \label {b.33} {\mathcal{N}}_g= -\sqrt{-\Delta_h} + B=-\sqrt{-\Delta_h} -B_0- B_{-1}-B_{-2} \quad \;\mbox{mod}\;\, OPS^{-3}(\partial \Omega),\end{eqnarray}
  where \begin{eqnarray} \label{z111} B_0 = -\frac{1}{2} \left[
  (\mbox{Tr}\, A_\nu) \,Id + \frac{\Gamma\big(\frac{n-2}{2}\big)}{4 \pi^{\frac{n}{2}}} \,\frac{1}{|x|^{n-2}} *\sum_{j,k=1}^n \big((C^{-1})^T A_\nu^* C^{-1}\big)_{jk}\, \frac{\partial^2 }{\partial x_j \partial x_k}  \right]\end{eqnarray}
    \begin{eqnarray} \label{z111-1} &&B_{-1} =- \frac{1}{8} \left[\bigg(-
2Q_1 + 2 \sum_{j=1}^n {\tilde {R}}_{j(n+1)j(n+1)} -2\,\mbox{Tr}\, A_\nu^2 +3 (\mbox{Tr}\, A_\nu)^2\bigg)
  2^{-1+\frac{n}{2}} \frac{ \Gamma \big(\frac{n-1}{2}\big)}{\Gamma(\frac{1}{2})} \, \frac{1}{|x|^{n-1}}\right. \\
  &&\quad \quad\quad  \left.+ \frac{5 \Gamma\big( \frac{n-5}{2}\big)}{2^5 \pi^{\frac{n}{2}} \Gamma\big(\frac{5}{2}\big)}  \, \frac{1}{|x|^{n-5}}
  * \bigg( \sum_{j,k=1}^n \big((C^{-1})^T A_\nu^* C^{-1}\big)_{jk} \frac{\partial^2 }{\partial x_j \partial x_k}\bigg)^2 \right.\nonumber\\
  &&\quad \quad \quad \left. + \frac{ \Gamma \big(\frac{n-3}{2}\big)}{4\pi^{\frac{n}{2}} \Gamma\big(\frac{3}{2}\big)} \, \frac{1}{|x|^{n-3}}
  * \sum_{j,k=1}^n \bigg( 2{\tilde{R}}_{j(n+1)k(n+1)} +6(A_\nu^2)_{jk}\bigg) \, \frac{\partial^2 }{\partial x_j \partial x_k}\right],
  \nonumber\\ &&B_{-2}\in {\mbox{OPS}}^{-2} (\partial \Omega),
\nonumber\end{eqnarray}
where  $f*g$ is the convolution of two functions $f$ and $g$.
    Furthermore, the principal, second and third symbols of $B$  are
            \begin{eqnarray} \label{b.40} p_0^B= \frac{1}{2} \left(\mbox{Tr}\, A_\nu  -  \frac{\langle A^*_\nu \xi,\xi\rangle}{\langle \xi,\xi\rangle}  \right),\end{eqnarray}
            \begin{eqnarray}\label{y6} \quad \;\qquad  p_{-1}^B(x, \xi)&=&\frac{1}{8} \left[-2Q_1
 + 2 \sum_{j=1}^n {\tilde {R}}_{j(n+1)j(n+1)} -2\,\mbox{Tr}\, A_\nu^2 +3 (\mbox{Tr}\, A_\nu)^2 + \frac{5\langle A_\nu^* \xi, \xi\rangle^2}{\langle \xi, \xi\rangle^2} \right.\\
  && \left. - \frac{1}{\langle \xi,\xi\rangle}\sum_{j,k=1}^n  \left(2{\tilde{R}}_{j(n+1)k(n+1)} +6(A_\nu^2)_{jk} \right)\xi_j \xi_k  \right] \frac{1}{\sqrt{\langle \xi, \xi\rangle}}\nonumber\end{eqnarray} and
\begin{eqnarray} \label{4-080}   p_{-2}^B (x, \xi)= \frac{\langle A_\nu^* \xi, \xi\rangle}{8\langle \xi, \xi\rangle^2}   \bigg( 2Q_1 -2 \sum_{j=1}^n {\tilde R}_{j(n+1)j(n+1)} +2 \mbox{Tr}\, A_\nu^2
-3 (\mbox{Tr}\, A_\nu)^2 \end{eqnarray}
\begin{displaymath}  \,\, -\frac{5\langle A_\nu^* \xi, \xi\rangle^2}{\langle \xi, \xi\rangle^2}  +\frac{1}{\langle \xi, \xi\rangle}\, \sum_{j,k=1}^n \big(2{\tilde R}_{j(n+1)k(n+1)} +6 (A_\nu^2)_{jk} \big) \xi_j\xi_k \bigg)\qquad \qquad \qquad \quad
    \nonumber \\ \end{displaymath}
     \begin{eqnarray*} &&  \quad \;\;- \frac{1}{4 \langle \xi, \xi\rangle}\, \left\{ \frac{3}{2} (\mbox{Tr}\, A_\nu) \bigg( 2Q_1 -2 \sum_{j=1}^n {\tilde R}_{j(n+1)j(n+1)} + 2 \mbox{Tr}\, A_\nu^2 \bigg)-4(\mbox{Tr}\, A_\nu)^3 \right. \\ && \left.\quad \quad
     +\frac{1}{4} \bigg[6Q_{2} +3Q_{3} +\sum_{j=1}^n\bigg(-2 {\tilde R}_{j(n+1)j(n+1),(n+1)} +4 \sum_{l=1}^n {\tilde R}_{l(n+1)j(n+1)} (A_\nu)_{lj}\bigg)\bigg]\right.\end{eqnarray*}
           \begin{eqnarray*} && \left. -\frac{1}{2} \bigg[\frac{\langle A_\nu^* \xi, \xi\rangle}{ \langle \xi, \xi\rangle} + \mbox{Tr}\, A_\nu\bigg] \bigg[ \frac{1}{2\langle \xi, \xi\rangle} \bigg(\sum_{j,k=1}^n \big(2 {\tilde R}_{j(n+1)k(n+1)}+6 (A_\nu^2)_{jk} \big)  \xi_j\xi_k \bigg)\right. \quad\quad \;\\
      && \ \left.-\frac{2\langle A_\nu^* \xi, \xi\rangle^2}{\langle \xi, \xi\rangle^2} +2 (\mbox{Tr}\, A_\nu)^2 -\frac{1}{2}
      \bigg(2Q_1- 2\sum_{j=1}^n {\tilde R}_{j(n+1)j(n+1)} +2 \mbox{Tr}\, A_\nu^2\bigg)\bigg]\right.\quad \quad\;\end{eqnarray*} \begin{eqnarray*}
      && \quad \qquad\left. +\frac{1}{4\langle \xi,\xi\rangle}
            \bigg[\sum_{j,k=1}^n \bigg(2{\tilde R}_{j(n+1)k(n+1),(n+1)} +20 \sum_{l=1}^n {\tilde R}_{j(n+1)l(n+1)}(A_\nu)_{lk} +24 (A_\nu^3)_{jk}\bigg) \xi_j\xi_k\bigg]\right.\quad \;\; \\ && \quad \qquad\left.
            -\frac{3\langle A_\nu^* \xi, \xi\rangle}{2\langle \xi, \xi\rangle^2} \bigg[ \sum_{j,k=1}^n\bigg(2 {\tilde R}_{j(n+1)k(n+1)}
        +6 (A_\nu^2)_{jk}\bigg)\xi_j\xi_k\bigg]  \right.\quad \;\;\\   && \quad \qquad\left. +  \frac{4 \langle A_\nu^* \xi, \xi\rangle^3 }{\langle \xi,\xi\rangle^3} + \frac{\mbox{Tr}\, A_\nu}{2} \bigg[\frac{1}{2\langle \xi,\xi\rangle} \bigg( \sum_{j,k=1}^n
      \big( 2{\tilde R}_{j(n+1)k(n+1)} +6 (A_\nu^2)_{jk} \big)\xi_j\xi_k\bigg) \right.\\
      && \quad \qquad\left.-\frac{2 \langle A_\nu^* \xi, \xi\rangle^2}{\langle \xi, \xi\rangle^2} +2 (\mbox{Tr}\, A_\nu)^2 - \frac{1}{2} \bigg (2Q_1-2 \sum_{j=1}^n {\tilde R}_{j(n+1)j(n+1)} + 2\mbox{Tr}\, A_\nu^2\bigg)\bigg] \right.\\
       &&  \quad \qquad \left.-\frac{1}{4} \bigg[ \frac{\langle A_\nu^* \xi, \xi\rangle}{\langle \xi, \xi\rangle} +\mbox{Tr}\, A_\nu\bigg] \bigg[ 2Q_1 -2\sum_{j=1}^n
       {\tilde R}_{j(n+1)j(n+1)} + 2 \mbox{Tr}\, A_\nu^2 -4(\mbox{Tr}\; A_\nu)^2 \bigg] \right\}, \end{eqnarray*}
respectively. }

 \vskip 0.4 true cm

\noindent  {\it Proof.} \  We choose coordinates $x=(x_1,\cdots, x_n)$ on an open set in $\partial \Omega$ and then
 coordinates $(x, x_{n+1})$ on a neighborhood in $\bar\Omega$ such that
$x_{n+1}=0$ on $\partial \Omega$ and $|\nabla x_{n+1}|=1$ near $\partial \Omega$ while $x_{n+1}>0$ on $\Omega$ and such that $x$ is constant
on each geodesic segment in $\bar\Omega$ normal to $\partial \Omega$. Then the metric tensor on $\bar \Omega$ has
the form (see, p.$\,$532 of \cite{Ta2})
\begin{gather} \label{a-1} \big(g_{jk} (x,x_{n+1}) \big)_{(n+1)\times (n+1)} =\begin{pmatrix} ( h_{jk} (x,x_{n+1}))_{n\times n}& 0\\
      0& 1 \end{pmatrix}. \end{gather}
 The Laplace-Beltrami operator $\Delta_g$ on $\Omega$ is given in local coordinates by
\begin{eqnarray*} \label{a-2} \Delta_g u&=&\sum_{j,k=1}^{n+1} \frac{1}{\sqrt{|g|}}\, \frac{\partial}{\partial x_j}\left( \sqrt{|g|}\, g^{jk} \frac{\partial u}{\partial x_k}\right) \\
&=&\frac{1}{\sqrt{|h|}} \,\frac{\partial}{\partial x_{n+1}}\left( \sqrt{|h|} \,\frac{\partial u}{\partial x_{n+1}}\right) +\sum_{j,k=1}^{n} \frac{1}{\sqrt{|h|}} \frac{\partial}{\partial x_j}\left( \sqrt{|h|} h^{jk} \,\frac{\partial u}{\partial x_k}\right)\nonumber\\
&=&\frac{\partial^2 u}{\partial x_{n+1}^2} + a(x_{n+1}) \frac{\partial u}{\partial x_{n+1}}  +
 L(x, x_{n+1}, D_{x}) u, \quad \;\;  a(x_{n+1})=\frac{1}{2|h|}\, \frac{\partial |h|}{\partial x_{n+1}},  \nonumber\end{eqnarray*}
where, as usual,
$|g|=\mbox{det}(g_{jk}), \, \, |h|=\mbox{det}(h_{jk})$, and $L(x_{n+1})= L(x, x_{n+1}, D_x) $ is a family of Laplace-Beltrami operators on $\partial \Omega$, associated to the family of metric $\left(h_{jk}\right)$ on $\partial \Omega$, so $L(0)= \Delta_{h}$.
Similar to the proof of Proposition C.1 of \cite{Ta2}, we will construct smooth families of operators $A_j(x_{n+1})\in OPS^1(\partial \Omega)$ such that
\begin{eqnarray} \label{z1} &&\frac{\partial^2}{\partial x_{n+1}^2}  +a(x_{n+1}) \frac{\partial}{\partial x_{n+1}}  + L(x_{n+1})\\
 && \quad \; =\left( \frac{\partial }{\partial x_{n+1}} -A_1(x_{n+1})\right) \left(\frac{\partial }{\partial x_{n+1}} +A_2(x_{n+1})\right),\nonumber\end{eqnarray}
modulo a smooth operator. It will follow that the principal parts of $A_1(x_{n+1})$ and $A_2(x_{n+1})$  is
 $\sqrt{-L(x_{n+1})}$, and
 \begin{eqnarray} \label{z2} {\mathcal{N}}_g=-A_2(0)\quad \, \mbox{mod}\;\; OPS^{-\infty} (\partial \Omega).\end{eqnarray}
In view of $[\frac{\partial }{\partial x_{n+1}}, A_2(x_{n+1})]=
\frac{\partial A_2(x_{n+1})}{\partial x_{n+1}}$ (see (4.9) of \cite{Mol}), it follows that
  the right-hand side of (\ref{z1}) is equal to
\begin{eqnarray*} \label{z3} \frac{\partial^2}{\partial x_{n+1}^2} -A_1(x_{n+1}) \frac{\partial }{\partial x_{n+1}}  +A_2 (x_{n+1}) \frac{\partial}{\partial x_{n+1}}
+\frac{\partial A_2(x_{n+1})}{\partial x_{n+1}} -A_1(x_{n+1})A_2(x_{n+1}).\end{eqnarray*}
Thus we get  \begin{eqnarray*} \label{z4}\left\{ \begin{array}{ll} A_2(x_{n+1}) -A_1(x_{n+1}) = a(x_{n+1}),\\
 -A_1(x_{n+1}) A_2(x_{n+1}) +\frac{\partial A_2(x_{n+1})}{\partial x_{n+1}} = L(x_{n+1}),\end{array}\right.\end{eqnarray*}
 from which we have an equation for $A_1(x_{n+1})$:
\begin{eqnarray*} A_1(x_{n+1})^2 +A_1(x_{n+1}) a(x_{n+1}) -\frac{\partial A_1(x_{n+1})}{\partial x_{n+1}} =-L(x_{n+1})+ \frac{\partial a(x_{n+1})}{\partial x_{n+1}}.\end{eqnarray*}
 Setting \begin{eqnarray}\label{z9}
 A_1(x_{n+1})= \Upsilon(x_{n+1}) +B(x_{n+1}), \quad \, \Upsilon(x_{n+1})=\sqrt{-L(x_{n+1})},\end{eqnarray}
     we obtain an equation for $B(x_{n+1})$:
\begin{eqnarray} \label{z6} && 2B(x_{n+1}) \Upsilon (x_{n+1}) +[\Upsilon(x_{n+1}), B(x_{n+1})]+B(x_{n+1})^2 -\frac{\partial B(x_{n+1})}{\partial x_{n+1}} \\ && \quad \quad \;\;  +B(x_{n+1}) a(x_{n+1})
= \frac{\partial \Upsilon(x_{n+1})}{\partial x_{n+1}} -\Upsilon(x_{n+1}) a(x_{n+1}) +\frac{\partial a(x_{n+1})}{\partial x_{n+1}}.\nonumber \end{eqnarray}
We will inductively obtain terms $B_j(x_{n+1}) \in OPS^{-j}(\partial \Omega)$ and establish that, with $B(x_{n+1})\sim
\sum_{j\ge 0} B_j(x_{n+1})$, the operators
   \begin{eqnarray*} A_1(x_{n+1}) =\sqrt{-L(x_{n+1})} +B(x_{n+1}), \quad \; A_2(x_{n+1}) =\sqrt{-L(x_{n+1})} +B(x_{n+1})+a(x_{n+1})\end{eqnarray*}
   do yield (\ref{z1}) modulo a smooth operator.
  Since $B(x_{n+1})$ is a smooth family in $OPS^{0}(\partial \Omega)$, the principal part $B_0(x_{n+1})$ must satisfy
   $2B_0(x_{n+1}) \Upsilon(x_{n+1})= \frac{\partial \Upsilon (x_{n+1})}{\partial x_{n+1}} -a(x_{n+1}) \Upsilon (x_{n+1})$, i.e.,
   \begin{eqnarray} \label{z7} B_0(x_{n+1}) =\frac{1}{2} \,\frac{\partial \Upsilon(x_{n+1})}{\partial x_{n+1}} \Upsilon(x_{n+1})^{-1}  -\frac{1}{2} a(x_{n+1})\quad \;\, OPS^{-1} (\partial \Omega).\end{eqnarray}
   By the second identity of (\ref{z9}) we derive
    \begin{eqnarray}
\label{z7'}    \left(\frac{\partial \Upsilon(x_{n+1})}{\partial x_{n+1}}(0)\right)\Upsilon(0)^{-1} =\frac{1}{2} \left(\frac{\partial L(x_{n+1})}{\partial x_{n+1}}(0)\right) L(0)^{-1} \quad \; \mbox{mod}\;\;  OPS^{-1}(\partial \Omega).\end{eqnarray}
           Therefore,
           \begin{eqnarray} \label{z11}  B_0(0) = \frac{1}{4} \left[\left(\frac{\partial L(x_{n+1})}{\partial x_{n+1}}(0)\right) L(0)^{-1} -|h|^{-1}\frac{\partial |h|}{\partial  x_{n+1}} (x,0)\right].\end{eqnarray}
Next, noting that $B(x_{n+1})$ is also a smooth family of $OPS^{-1}(\partial \Omega)$, by compare the second symbol of $B(x_{n+1})$ we have
\begin{eqnarray*}  2B_{-1} (x_{n+1}) \Upsilon(x_{n+1}) + B_0(x_{n+1})^2 -\frac{\partial B_0(x_{n+1})}{\partial x_{n+1}} +
B_0(x_{n+1}) a(x_{n+1}) =\frac{\partial a(x_{n+1})}{\partial x_{n+1}},\end{eqnarray*}
i.e.,  \begin{eqnarray} \label{y1}  B_{-1} (x_{n+1})& =&\frac{1}{2} \left( \frac{\partial a(x_{n+1})}{\partial x_{n+1}} -B_0(x_{n+1})^2 +\frac{\partial B_0(x_{n+1})}{\partial x_{n+1}}\right.\\
&& \left. -B_0(x_{n+1})\frac{{}_{}}{{}_{}} a(x_{n+1})\right) \Upsilon (x_{n+1})^{-1}.\nonumber\end{eqnarray}
From (\ref{z7}) we get \begin{eqnarray*} \frac{\partial B_0(x_{n+1})}{\partial x_{n+1}} =\frac{1}{2} \left[ \frac{\partial^2 \Upsilon (x_{n+1})}{\partial x_{n+1}^2} \Upsilon (x_{n+1})^{-1} - \left(\frac{\partial \Upsilon (x_{n+1})}{\partial x_{n+1}}\right)^2 \Upsilon (x_{n+1})^{-2} - \frac{\partial a(x_{n+1})}{\partial x_{n+1}}\right],\end{eqnarray*}
so that
\begin{eqnarray*} && B_{-1} (x_{n+1}) =\frac{1}{2} \bigg[ \frac{1}{2} \,\frac{\partial a(x_{n+1})}{\partial x_{n+1}} -B_0(x_{n+1})^2 +
\frac{1}{2} \frac{\partial^2 \Upsilon (x_{n+1})}{\partial x_{n+1}^2} \Upsilon (x_{n+1})^{-1}\\
&&\quad \quad \;- \frac{1}{2}\, \bigg(\frac{\partial \Upsilon (x_{n+1})}{\partial x_{n+1}}\bigg)^2 \Upsilon (x_{n+1})^{-2}
 -B_0(x_{n+1}) a(x_{n+1}) \bigg]\Upsilon (x_{n+1})^{-1}.
\end{eqnarray*}
Again, by using the second identity of (\ref{z9}), we have
\begin{eqnarray*}
 && 2\Upsilon(x_{n+1}) \, \frac{\partial^2 \Upsilon (x_{n+1})}{\partial x_{n+1}^2} +2 \left( \frac{\partial \Upsilon (x_{n+1})}{\partial x_{n+1}}\right)^2 =-\frac{\partial^2 L(x_{n+1})}{\partial x_{n+1}^2},\\
 && \quad \, \,\left( \frac{\partial \Upsilon (x_{n+1})}{\partial x_{n+1}}\right)^2= \frac{1}{4} \left( \frac{\partial L (x_{n+1})}{\partial x_{n+1}}\right)^2\Upsilon (x_{n+1})^{-2}.\end{eqnarray*}
Combining these identities we obtain
\begin{eqnarray*}  B_{-1} (x_{n+1})& =&\frac{1}{2} \left[ \frac{1}{2} \,\frac{\partial a(x_{n+1})}{\partial x_{n+1}} -B_0(x_{n+1})^2-
\frac{1}{4} \, \frac{\partial^2 L(x_{n+1})}{\partial x_{n+1}^2} (-L(x_{n+1}))^{-1}\right. \\
 &&\left.- \frac{1}{4} \left(\frac{\partial L(x_{n+1})}{\partial x_{n+1}}\right)^2 L(x_{n+1})^{-2} -B_0(x_{n+1}) a(x_{n+1})\right] {\sqrt{-L(x_{n+1})}\,}^{-1},
\end{eqnarray*}
so that
\begin{eqnarray*}  B_{-1} (0)& =&\frac{1}{2} \left[ \frac{1}{2} \,\frac{\partial a(x_{n+1})}{\partial x_{n+1}}(x,0) -B_0(0)^2-
\frac{1}{4}  \left(\frac{\partial^2 L(x_{n+1})}{\partial x_{n+1}^2}(0)\right) (-L(0))^{-1} \right.\\
&& \left.- \frac{1}{4} \left(\frac{\partial L(x_{n+1})}{\partial x_{n+1}}(0)\right)^2 L(0)^{-2} -B_0(0) a(0)\right] {\sqrt{-L(0)}\,}^{-1}.
\end{eqnarray*}
Inserting (\ref{z11})
into the above identity, we get
\begin{eqnarray} \label{z60} B_{-1} (0) &=&\frac{1}{2} \left[ \frac{1}{4|h|} \,\frac{\partial^2 |h|}{\partial x_{n+1}^2}(x,0)-\frac{3}{16 |h|^2}\left(\frac{\partial |h|}{\partial x_{n+1}}(x,0)\right)^2\right. \nonumber \\ && \left.- \frac{5}{16} \left(\frac{\partial L(x_{n+1})}{\partial x_{n+1}}(0)\right)^2 L(0)^{-2} -\frac{1}{4}  \left(\frac{\partial^2 L(x_{n+1})}{\partial x_{n+1}^2}(0)\right) \big(-L(0)\big)^{-1} \right] {\sqrt{-L(0)}\,}^{-1}.\nonumber
\end{eqnarray}
Furthermore, by compare the third symbol of $B(x_{n+1})$ we get
\begin{eqnarray*} 2B_{-2}(x_{n+1})\Psi(x_{n+1}) +2B_0(x_{n+1})B_{-1}(x_{n+1}) -\frac{\partial B_{-1}(x_{n+1})}{\partial x_{n+1}} + a(x_{n+1}) B_{-1}(x_{n+1}) =0,\end{eqnarray*}
i.e., \begin{eqnarray} \label{4/04-1} B_{-2}(x_{n+1}) &=&\bigg( -B_0(x_{n+1})B_{-1}(x_{n+1}) +\frac{1}{2} \,\frac{\partial B_{-1}(x_{n+1})}{\partial x_{n+1}} \\  &&-\frac{a(x_{n+1})}{2} \,B_{-1}(x_{n+1}) \bigg) \Psi^{-1}(x_{n+1}),\nonumber\end{eqnarray}
where \begin{eqnarray} \label{4/05-2}  &&\frac{\partial B_{-1}(x_{n+1})}{\partial x_{n+1}}= \frac{1}{4} \bigg[ \frac{1}{2} \, \frac{\partial a(x_{n+1})}{\partial x_{n+1}} -B_0(x_{n+1})^2 -\frac{1}{4} \, \frac{\partial^2L(x_{n+1})}{\partial x_{n+1}^2} (-L(x_{n+1}))^{-1}\\
&& -\frac{1}{4} \left(\frac{\partial L(x_{n+1})}{\partial x_{n+1}}\right)^2 L(x_{n+1})^{-2} - a(x_{n+1})B_0(x_{n+1})  \bigg] \frac{\partial L(x_{n+1})} {\partial x_{n+1}}(-L(x_{n+1}))^{-\frac{3}{2}} \nonumber\\
&&+ \frac{1}{2} \bigg[ \frac{1}{2} \,\frac{\partial^2 a(x_{n+1})}{\partial x_{n+1}^2}-2B_0(x_{n+1}) \, \frac{\partial B_0(x_{n+1})}{\partial x_{n+1}}  +\frac{1}{4} \, \frac{\partial^3 L(x_{n+1})}{\partial x_{n+1}^3}  L(x_{n+1})^{-1} \nonumber \\&&-\frac{3}{4} \, \frac{\partial^2 L(x_{n+1})}{\partial x_{n+1}^2} \frac{\partial L(x_{n+1})}{\partial x_{n+1}}  L(x_{n+1})^{-2} +\frac{1}{2} \big(\frac{\partial L(x_{n+1})}{\partial x_{n+1}}\big)^2
\, \frac{\partial L(x_{n+1})}{\partial x_{n+1}} L(x_{n+1})^{-3} \nonumber\\
&& -a(x_{n+1}) \,\frac{\partial B(x_{n+1})}{\partial x_{n+1}}
 -B_0(x_{n+1}) \frac{\partial a(x_{n+1})}{\partial x_{n+1}}\bigg] \sqrt{-L(x_{n+1})}^{\;-1}.\nonumber\end{eqnarray}
 Hence \begin{eqnarray*}\label{y7} {\mathcal{N}}_g&=&-A_2(0) =- \sqrt{-L(0)} -\big(B_0(0) +a(0)\big)-B_{-1}(0)-B_{-2}(0)\\
   &=& -\sqrt{-\Delta_g}- \frac{1}{4} \left[\left(\frac{L(x_{n+1})}{\partial x_{n+1}}(0)\right)L(0)^{-1} +|h|^{-1}\frac{\partial |h|}{\partial x_{n+1}}(x,0)\right]\nonumber\\
   && -  \frac{1}{2} \left[ \frac{1}{4|h|} \,\frac{\partial^2 |h|}{\partial x_{n+1}^2}(x,0)-\frac{3}{16|h|^2}\left(\frac{\partial |h|}{\partial x_{n+1}}(x,0)\right)^2 \right. \nonumber\\ && \left. - \frac{5}{16} \left(\frac{\partial L(x_{n+1})}{\partial x_{n+1}}(0)\right)^2 L(0)^{-2}
  -\frac{1}{4}  \left(\frac{\partial^2 L(x_{n+1})}{\partial x_{n+1}^2}(0)\right) \big(-L(0)\big)^{-1} \right] {\sqrt{-L(0)}\,}^{-1}-B_{-2}(0)\nonumber\\
   && \quad\;\;\mbox{mod}\;\; OPS^{-3}(\partial \Omega).\nonumber\end{eqnarray*}

We will compute the first,  second and  third symbols of $B$. It is obvious that the principal symbols of $L(0)$ and $\frac{\partial L(x_{n+1})}{\partial x_{n+1}}(0)$ are $-\sum_{j,k=1}^n  h^{jk} (x, 0) \xi_j\xi_k =-\langle \xi,\xi\rangle$ and
$- \sum \frac{\partial h^{jk}}{\partial x_{n+1}}(x,0) \xi_j\xi_k$, respectively.
      Furthermore, we choose a normal coordinate system on $(\partial \Omega, h)$, centered at $x_0\in \partial \Omega$.  From (C.24) of \cite{Ta2} (see also, (4.68)--(4.70) of Appendix C of \cite{Ta2}), one has
\begin{eqnarray} \label{6610} &&\sum_{j,k=1}^n \frac{\partial h^{jk}}{\partial x_{n+1}} (x_0,0)  \xi_j\xi_k =2\langle A_\nu^* \xi, \xi\rangle,\end{eqnarray}
so that the principal symbol of $\frac{\partial L(x_{n+1})}{\partial x_{n+1}}(0) \, L(0)^{-1}$  is $ \frac{2\langle A^*_\nu \xi, \xi\rangle}{\langle \xi, \xi\rangle}$. It follows from (C.26) of \cite{Ta2} that  \begin{eqnarray} \label{6620}
&& |h|^{-1} \frac{\partial |h|}{\partial x_{n+1}} (x_0,0) =\sum_{j=1}^n  \frac{\partial h_{jj}}{\partial x_{n+1}}(x_0,0) =-2\, \mbox{Tr}\, A_\nu.
\end{eqnarray}
Consequently  \begin{eqnarray*} \label{y5} p_0^B= \frac{1}{2} \left(\mbox{Tr}\, A_\nu -
\frac{\langle A^*_\nu \xi,\xi\rangle}{\langle \xi,\xi\rangle}
\right)\end{eqnarray*}
and \begin{eqnarray}\label{y60} p_{-1}^B(x_0, \xi)&=&-\frac{1}{2} \left[\frac{1}{4|h|}\, \frac{\partial^2 |h|}{\partial x_{n+1}^2} (x_0,0) -\frac{3}{16}\left(\frac{1}{|h|}\,  \frac{\partial |h|}{\partial x_{n+1}}(x_0,0)\right)^2 \right.\\
 && \left.-\frac{5}{16} \left(\sum_{j,k=1}^n\frac{\partial h^{jk}}{\partial x_{n+1}}(x_0,0) \, \xi_j\xi_k\right)^2
\left(\sum_{j,k=1}^n h^{jk}(x_0,0) \, \xi_j\xi_k\right)^{-2}\right.\nonumber\\
 &&\left. +\frac{1}{4} \left(\sum_{j,k=1}^n\frac{\partial^2 h^{jk}}{\partial x_{n+1}^2}(x_0,0) \, \xi_j\xi_k\right)
\left(\sum_{j,k=1}^n h^{jk}(x_0,0) \, \xi_j\xi_k\right)^{-1}
\right]\nonumber\\
&& \quad \, \times \left(\sum_{j,k=1}^n h^{jk}(x_0,0) \, \xi_j\xi_k\right)^{-1/2} \nonumber\end{eqnarray}
\begin{eqnarray}
 &&
 =-\frac{1}{8} \left[\frac{\partial^2 |h|}{\partial x_{n+1}^2} (x_0,0) -3\left(\mbox{Tr}\, A_\nu\right)^2- \frac{5\langle A^*_\nu \xi, \xi\rangle^2}{\langle\xi,\xi\rangle^2} \right.\nonumber\\ && \quad \; \left. +\frac{1}{\langle \xi,\xi\rangle} \left(\sum_{j,k=1}^n\frac{\partial^2 h^{jk}}{\partial x_{n+1}^2}(x_0,0) \, \xi_j\xi_k\right)
 \right]\langle \xi, \xi\rangle^{-1/2}.\quad \;\nonumber\end{eqnarray}
  Note that  $\left(-\frac{1}{2}\,\frac{\partial h_{jk}}{\partial x_{n+1}}(x_0,0)\right)_{n\times n}$  is the matrix of Weingarten's map under a basis of $T_0(\partial \Omega)$ at $x_0$, and its eigenvalues are just the principal curvatures $\kappa_1,\cdots,\kappa_n$ of $\partial \Omega$ in the direction $\nu$.
  We can write
          \begin{eqnarray} \label{5,,1}   g_{jk} (x)  &=&  \delta_{jk} +\sum_{l=1}^{n+1} \frac{\partial g_{jk}}{\partial x_{l}}(0) \, x_{l}+
          \frac{1}{2} \sum_{l,m=1}^{n+1} \frac {\partial^2 g_{jk} }{\partial x_l\partial x_m} (0)\, x_l x_m \\
         && +\frac{1}{6} \sum_{l,m,v=1}^{n+1}
  \frac{\partial^3 g_{jk}}{\partial x_l \partial x_m\partial x_v}(0)\, x_l x_m x_v +O(|x|^4)  \; \; \mbox{near}\;\, 0. \nonumber\end{eqnarray}
            By taking $x=(0,\cdots, 0, x_{n+1}):=x_{n+1} e_{n+1}$, we find that  \begin{eqnarray*} &&|h(0+x_{n+1}e_{n+1})| = |g(0+x_{n+1}e_{n+1})|\\
           && \quad \; = \mbox{det}\, \left(\delta_{jk} + \frac{\partial g_{jk}}{\partial x_{n+1}}(0) x_{n+1} + \frac{1}{2}\,
\frac{\partial^2 g_{jk}}{\partial x_{n+1}^2} (0)x_{n+1}^2 +\frac{1}{6}\,
  \frac{\partial^3 g_{jk}}{\partial x_{n+1}^3}(0)\, x_{n+1}^3 + O(|x_{n+1}|^4)\right),\end{eqnarray*}
so that  \begin{eqnarray} \label{6601}\quad \quad\quad\; \frac{\partial^2 |h|}{\partial x_{n+1}^2} (0) &=& 2\sum_{1\le j<k\le n}
  \bigg( \frac{\partial h_{jj}}{\partial x_{n+1}} (0)\, \frac{\partial h_{kk}}{\partial x_{n+1}} (0) - \frac{\partial h_{jk}}{\partial x_{n+1}} (0)\, \frac{\partial h_{kj}}{\partial x_{n+1}} (0)\bigg)
+\sum_{j=1}^{n} \frac{\partial^2 h_{jj}}{\partial x_{n+1}^2} (0) \\
    &:=& 2Q_1 +\sum_{j=1}^{n} \frac{\partial^2 h_{jj}}{\partial x_{n+1}^2} (0),\nonumber
 \end{eqnarray}
  \begin{eqnarray} \label{66020} \frac{\partial^3 |h|}{\partial x_{n+1}^3} (0) &=& 6\sum_{1\le j<k<l\le n}
  \bigg( \frac{\partial h_{jj}}{\partial x_{n+1}} (0)\, \frac{\partial h_{kk}}{\partial x_{n+1}} (0) \frac{\partial h_{ll}}{\partial x_{n+1}} (0)
  + \frac{\partial h_{jk}}{\partial x_{n+1}} (0)\, \frac{\partial h_{kl}}{\partial x_{n+1}} (0) \frac{\partial h_{lj}}{\partial x_{n+1}} (0)
 \nonumber \\&& +\frac{\partial h_{jl}}{\partial x_{n+1}} (0)\, \frac{\partial h_{kj}}{\partial x_{n+1}} (0) \frac{\partial h_{lk}}{\partial x_{n+1}} (0)  - \frac{\partial h_{jl}}{\partial x_{n+1}} (0)\, \frac{\partial h_{kk}}{\partial x_{n+1}} (0) \frac{\partial h_{lj}}{\partial x_{n+1}} (0)
  \nonumber \\ \label{4.222} \quad   &&-\frac{\partial h_{jk}}{\partial x_{n+1}} (0)\, \frac{\partial h_{kj}}{\partial x_{n+1}} (0) \frac{\partial h_{ll}}{\partial x_{n+1}} (0)
   -\frac{\partial h_{jj}}{\partial x_{n+1}} (0)\, \frac{\partial h_{kl}}{\partial x_{n+1}} (0) \frac{\partial h_{lk}}{\partial x_{n+1}} (0)
  \bigg)  \\&& + 3\sum_{j\ne k}  \bigg( \frac{\partial h_{jj}}{\partial x_{n+1}} (0)\, \frac{\partial^2 h_{kk}}{\partial x_{n+1}^2} (0) - \frac{\partial h_{jk}}{\partial x_{n+1}} (0)\, \frac{\partial^2 h_{kj}}{\partial x_{n+1}^2} (0)\bigg)
+ \sum_{j=1}^{n} \frac{\partial^3 h_{jj}}{\partial x_{n+1}^3} (0)\nonumber\\
&:=& 6Q_2 +3Q_3 + \sum_{j=1}^{n} \frac{\partial^3 h_{jj}}{\partial x_{n+1}^3} (0). \nonumber
 \end{eqnarray}
 It is well-known (see, for example, \cite{Jo},  \cite{Spi2}, \cite{CLN}, \cite{GMa}) that $Q_1=4\sum_{1\le j<k\le n} \kappa_j\kappa_k =2\big(R_{\partial \Omega}- {\tilde{R}}+2{\tilde{R}}_{\nu\nu}\big)$ and $Q_2 =-8\sum_{1\le j<k<l\le n} \kappa_j\kappa_k\kappa_l$.
    A simple calculation shows that, at $x_0=0$, \begin{eqnarray*} \label{6604}\quad \, \Gamma_{(n+1)(n+1)}^s =0, \quad \, \Gamma_{(n+1)k}^s= \frac{1}{2}\sum_{l=1}^{n+1}g^{sl}\frac{\partial g_{lk}}{\partial x_{n+1}}, \quad
   \Gamma_{j(n+1)}^t=\frac{1}{2}\sum_{l=1}^{n+1} g^{tl}\, \frac{\partial g_{lj}}{\partial x_{n+1}},\end{eqnarray*}
where $\Gamma_{\beta \gamma}^\alpha$  are the Christoffel symbols. It follows from p.$\,$188 of \cite{Spi2} (or \cite{Gud}) that
\begin{eqnarray*} {\tilde{R}}_{j(n+1)k(n+1)} &=& -\frac{1}{2} \, \frac{\partial^2 g_{jk}}{\partial x_{n+1}^2}(0) + \sum_{s,t=1}^{n+1} g_{st} \left(\Gamma_{(n+1)k}^s  \Gamma_{j(n+1)}^t - \Gamma_{(n+1)(n+1)}^s \Gamma_{jk}^t\right) \\
 &=&-\frac{1}{2} \,\frac{\partial^2 h_{jk}}{\partial x_{n+1}^2}(0) +\frac{1}{4} \sum_{t,l=1}^n  \frac{\partial h_{tk}} { \partial x_{n+1}}(0) h^{tl}(0) \,\frac{\partial h_{lj}}{\partial x_{n+1}}(0)\\
  &=& -\frac{1}{2} \,\frac{\partial^2 h_{jk}}{\partial x_{n+1}^2}(0) +(A_\nu^2)_{jk},\end{eqnarray*}
 \begin{eqnarray*}  {\tilde R}_{j(n+1)k(n+1), (n+1)} &=& -\frac{1}{2} \, \frac{\partial^3 h_{jk}}{\partial x_{n+1}^3}(0) -\frac{1}{4}\, \sum_{m,l=1}^{n+1} \big(\frac{\partial h_{km}}{\partial x_{n+1}} (0)\, \frac{\partial h_{ml}}{\partial x_{n+1}}(0)\,\frac{\partial h_{lj}}{\partial x_{n+1}}(0)\big)
 \\ &&  +\frac{1}{4} \sum_{l=1}^{n+1} \big(\frac{\partial^2 h_{lk}}{\partial x_{n+1}^2} (0)\, \frac{\partial h_{lj}}{\partial x_{n+1}}
(0)+\frac{\partial h_{lk}}{\partial x_{n+1}}(0)\,\frac{\partial^2 h_{lj}}{\partial x_{n+1}^2}(0)\big)\\
 &=& -\frac{1}{2} \, \frac{\partial^3 h_{jk}}{\partial x_{n+1}^3} (0)+ 2\sum_{l=1}^{n+1} {\tilde R}_{l(n+1)k(n+1)} (A_\nu)_{lj}, \end{eqnarray*}
 so that  \begin{eqnarray} \label{6605} \left. \begin{array}{l}\frac{\partial^2 h_{jk}}{\partial x_{n+1}^2}(x_0,0)=-2{\tilde{R}}_{j(n+1)k(n+1)} +2\,(A_\nu^2)_{jk}\\
  \label{6605'} \frac{\partial^3 h_{jk}}{\partial x_{n+1}^3}(x_0,0) =-2{\tilde R}_{j(n+1) k(n+1), (n+1)} + 4 \sum_{l=1}^n  {\tilde R}_{l(n+1) k(n+1)} (A_\nu)_{lj}.\end{array} \right.\end{eqnarray}
   Combining these and (\ref{6601})--(\ref{4.222}) we have
 \begin{eqnarray} \label{6602} \\
   \left. \begin{array}{l}  \frac{\partial^2 |h|}{\partial x_{n+1}^2} (x_0,0) = 2Q_1 - 2\sum_{j=1}^n  {\tilde{R}}_{j(n+1)j(n+1)}
 +2\,\mbox{Tr}\, A_\nu^2,\\
   \frac{\partial^3 |h|}{\partial x_{n+1}^3} (x_0,0)
    =6Q_2+3Q_3+   \sum_{j=1}^{n} \big(-2{\tilde R}_{j(n+1)j(n+1), (n+1)}+ 4 \sum_{l=1}^n {\tilde R}_{l(n+1)j(n+1)} (A_\nu)_{lj} \big) \end{array}\right.\nonumber
  \end{eqnarray}
  with $$Q_3= \sum_{j\ne k} \left[ (-2A_\nu)_{jj}  \big(-2{\tilde R}_{k(n+1)k(n+1)} +2(A_\nu^2)_{kk}\big) -(-2A_\nu)_{jk} \big(-2{\tilde R}_{k(n+1)j(n+1)} +2(A_\nu^2)_{kj}\big)\right].$$
  From $hh^{-1}=I$, we get \begin{eqnarray*} \label{6606}&&\frac{\partial h}{\partial x_{n+1}} h^{-1} +h \frac{\partial h^{-1}}{\partial x_{n+1}}=0,\quad\quad\,\,  \frac{\partial^2 h}{\partial x_{n+1}^2} h^{-1} +2 \frac{\partial h}{\partial x_{n+1}} \frac{\partial h^{-1}}{\partial x_{n+1}}  +h\, \frac{\partial^2 h^{-1}}{\partial x_{n+1}^2} =0,\\
&& \frac{\partial^3 h}{\partial x_{n+1}^3} h^{-1} +3\frac{\partial^2 h}{\partial x_{n+1}^2}\, \frac{\partial h^{-1}}{\partial x_{n+1}}
 +3 \frac{\partial h}{\partial x_{n+1}}\, \frac{\partial^2 h^{-1}}{\partial x_{n+1}^2}+ h \frac{\partial^3 h^{-1}}{\partial x_{n+1}^3} =0,\end{eqnarray*}
which imply   \begin{eqnarray*} \label{08} \frac{\partial^2 h^{jk}}{\partial x_{n+1}^2} (x_0, 0)
 &=& -\frac{\partial^2 h_{jk}}{\partial x_{n+1}^2} (x_0, 0) -2\sum_{l=1}^n \frac{\partial h_{jl}}{\partial x_{n+1}}(x_0,0)\, \frac{\partial h^{lk}}{\partial x_{n+1}} (x_0, 0)\\
 &=& -\frac{\partial^2 h_{jk}}{\partial x_{n+1}^2} (x_0, 0) +
8 (A_\nu^2)_{jk},\end{eqnarray*}
\begin{eqnarray*} \label{66090} \frac{\partial^3 h^{jk}}{\partial x_{n+1}^3} (x_0, 0)
 &=& -\frac{\partial^3 h_{jk}}{\partial x_{n+1}^3} (x_0, 0) -3\sum_{l=1}^n \frac{\partial^2 h_{jl}}{\partial x_{n+1}^2}(x_0,0)\, \frac{\partial h^{lk}}{\partial x_{n+1}} (x_0, 0)\\ &&
 -3\sum_{l=1}^n \frac{\partial h_{jl}}{\partial x_{n+1}}(x_0,0)\, \frac{\partial^2 h^{lk}}{\partial x_{n+1}^2} (x_0, 0).\nonumber\end{eqnarray*}
Thus \begin{eqnarray} \label{6612} \quad \quad \qquad \left.\begin{array}{l} \frac{\partial^2 h^{jk}}{\partial x_{n+1}^2} (x_0,0)=2{\tilde{R}}_{j(n+1)k(n+1)} + 6 (A_\nu^2)_{jk},\\
 \frac{\partial^3 h^{jk}}{\partial x_{n+1}^3} (x_0, 0)
 = 2{\tilde{R}}_{j(n+1)k(n+1),(n+1)}  + 20 \sum_{l=1}^{n} {\tilde R}_{j(n+1)l(n+1)} (A_\nu)_{lk} +24 (A_\nu^3)_{jk}.\end{array}\right.\end{eqnarray}
  Combining   (\ref{4/04-1}), (\ref{y60}),  (\ref{6605}), (\ref{6602}) and (\ref{6612}), we obtain (\ref{y6}) and (\ref{4-080}).
Finally, it follows from p.$\,$363 of \cite{GS} that the Fourier transform of $2^{\lambda+\frac{n}{2}} \, \frac{\Gamma\big(\frac{\lambda+n}{2}\big)}{\Gamma\big(-\frac{\lambda}{2}\big)} |x|^{-\lambda -n}$ is $|\eta|^{\lambda}$. Since $(h^{jk})$ is  positive-definite, there exists a matrix $C$ such that $C^T C= (h^{jk})$. By setting $\eta=C\xi$ and by using the fact that $\widehat{f*g}= (2\pi)^{\frac{n}{2}} \hat{f}\cdot \hat{g}$, we obtain the  expressions (\ref{b.33})--(\ref{z111-1})  for ${\mathcal N}_g$.

\vskip 0.35 true cm

\noindent{\bf Remark 4.3.} \ \ The proofs of Theorems 4.1 and 4.2 also provide a general method to calculate all $p_m(x, \xi)$ ($m=1,0,-1,-2,\cdots$) and $p^B_{m}(x,\xi)$ ($m=0,-1,-2,\cdots$) for the pseudodifferential operators $-\sqrt{-\Delta_h}$ and $B$, respectively.

\vskip 1.49 true cm

\section{Series representation of the heat kernel associated to the Dirichlet-to-Neumann operator }

\vskip 0.45 true cm

Let $\Omega$ is a bounded domain with smooth boundary $\partial \Omega$ in $(n+1)$-dimensional Riemannian manifold $(\mathcal{M},g)$.  Define $-\Delta^0_h$ to be the Laplacian $-\Delta_h$ on $\partial \Omega$ with its coefficients {\it frozen} at $y\in \partial \Omega$, and
let $\sqrt{-\Delta_h^0}$  be the square root operator of $-\Delta_h^0$,
where $h$ is the induced metric on $\partial \Omega$ by $g$. Denote by
 $G_0(t, x,y)$ and $F_0(t,x,y)$   the fundamental solutions of
$\frac{\partial u}{\partial t}= \Delta_g^0u$ and $\frac{\partial v}{\partial t} = -\sqrt{-\Delta_h^0} v$ evaluated $t>0$,
$x\in \partial \Omega$, and the {\it same} point $y\in \partial \Omega$
at which the coefficients of $-\Delta_g^0$ and $\sqrt{-\Delta_h^0}$ are computed, respectively. Leu us denote \begin{eqnarray*} \left(G_0 \# \left((\Delta_h -\Delta_h^0)G_0\right)\right)(t,x,y) = \int_0^t \int_{\partial \Omega}  G_0(s, x, z) \left((\Delta_h -\Delta_h^0)G_0(
  t-s, z, y)\right) ds\,dz.\end{eqnarray*}
 Furthermore \begin{eqnarray*} G_0 \# \left((\Delta_h -\Delta_h^0)G_0\right)\#\cdots \#\left((\Delta_h -\Delta_h^0)G_0\right) \quad (m\mbox{-fold})\end{eqnarray*}
 can be constructed inductively.

\vskip 0.45 true cm

\noindent{\bf Theorem 5.1.} \ {\it  Suppose  that $F(t,x,y)$ satisfies  \begin{eqnarray*} \left\{ \begin{array}{ll} \frac{\partial F(t, x,y)}{\partial t} =-\sqrt{-\Delta_h}\,F(t,x,y) \quad \, & \mbox{for}\;\; (t,x,y)\in[0, \infty) \times \partial \Omega \times \partial \Omega, \\
F(0,x,y)=\delta(x-y) \quad \, & \mbox{for}\;\; x,y\in \partial \Omega. \end{array} \right.\end{eqnarray*}
Then $F(t,x,y)$ has the following series representation:
 \begin{eqnarray} \label{-b.21} && F(t, x,y)= \int_0^\infty \frac{te^{-t^2/4\mu}}{\sqrt{4\pi \mu^3}}   G_0(\mu,x,y) d\mu \\
   && \quad \quad  + \sum_{m\ge 1}
    \int_0^\infty \frac{te^{-t^2/4\mu}}{\sqrt{4\pi \mu^3}}   \left[  G_0\# \underset{m}{\underbrace{\left((\Delta_h -\Delta_h^0)G_0\right) \# \cdots \#
 \left((\Delta_h -\Delta_h^0) G_0\right)}}(\mu,x,y)\right]d\mu.\nonumber\end{eqnarray}}

\vskip 0.2 true cm

  \noindent  {\it Proof.} \ It is well-known (see, for example, (5.22) of p$\,$247 of \cite{Ta1}, \cite{Gri} or \cite{Gr}) that for $\lambda\ge 0$,
\begin{eqnarray*} e^{-t\sqrt{\lambda}} =\int_0^\infty \frac{t}{\sqrt{4\pi \mu^3}} e^{-t^2/4\mu} e^{-\mu\lambda} d\mu,\end{eqnarray*}
i.e., the Laplace transform of $\frac{t}{\sqrt{4\pi \mu^3}} e^{-t^2/4\mu}$ is $e^{-t\sqrt{\lambda}}$.
By applying the spectral theorem, we get that for all $t>0$,
\begin{eqnarray} \label{b.12} e^{-t\sqrt{-\Delta_h}}\phi(x) = \int_0^\infty  \frac{t}{\sqrt{4\pi \mu^3}} e^{-t^2/4\mu} e^{\mu \Delta_h}\phi(x)d\mu,  \quad \; \forall \phi\in H^{\frac{1}{2}}(\partial \Omega). \end{eqnarray}
 Therefore, we have  \begin{eqnarray*} \label{b.13} e^{-t\sqrt{-\Delta_h}} \delta(x-y)  = \int_0^\infty \frac{t}{\sqrt{4\pi \mu^3}} e^{-t^2/4\mu} \left( e^{\mu \Delta_h} \delta(x-y)\right) d\mu, \end{eqnarray*}
i.e., \begin{eqnarray} \label{b.14} F(t, x,y) = \int_0^\infty   \frac{t}{\sqrt{4\pi \mu^3}} e^{-t^2/4\mu}G(\mu, x,y) d\mu, \end{eqnarray}
where $G(t, x, y)$ is the fundamental solution of heat equation $\frac{\partial u}{\partial t}= \Delta_h u$ on $[0,\infty)\times \partial \Omega$.
 It is obvious (see, for example, \cite{MS} or p.$\,$4 of \cite{Frie}) that
  in the normal coordinates
 \begin{eqnarray} \label{-.5.1} G_0(t, x, y)= \left(\frac{1}{4\pi t}\right)^{n/2}  e^{- \sum_{j,k=1}^n (h_{jk}(y))(x_j-y_j-b_j(y)t)(x_k-y_k -b_k(y) t)/4t},\end{eqnarray}
where \begin{eqnarray*} b_k(y) = \sum_{j=1}^n \frac{1}{\sqrt{|h(y)|}}\, \frac{\partial (\sqrt|h|\,h^{jk})}{\partial x_j}(y), \;\;  k=1,2,\cdots, n.\end{eqnarray*}
 From the well-known estimate (see (3.5b) of \cite{MS}):
   \begin{eqnarray} \label{-5-2} && \;\;\big|G_0\# \underset{m}{\underbrace{\left((\Delta_h- \Delta_h^0)G_0\right) \#\cdots \#  \left((\Delta_g- \Delta_h^0)G_0\right)(t,x,y)}}\big|
 \le  \frac{c_2^m}{(m/2)!} \, t^{(m-n)/2} e^{ -c_3 (d(x,y))^2/4t},\end{eqnarray}
 one immediately get  that
  \begin{eqnarray} \label{-5-2}\,\,\quad \quad  G= G_0 + \sum_{m\ge 1}  G_0\# \underset{m}{\underbrace{\left((\Delta_h -\Delta_h^0)G_0\right) \# \cdots \#
 \left((\Delta_h -\Delta_h^0) G_0\right)}} .\end{eqnarray}
 Note that the sum of right-hand side converges rapidly to the unique fundamental solution $G$.
Simple calculation shows that
 \begin{eqnarray} \label{-b-d-1}  && \int_0^\infty \frac{te^{-t^2/4\mu}}{\sqrt{4\pi \mu^3}}   G_0(\mu,x,y)d\mu= \int_0^\infty \frac{te^{-t^2/4\mu}}{\sqrt{4\pi \mu^3}} \frac{1}{(4\pi \mu)^{n/2}}e^{-\sum_{j,k=1}^n h_{jk} (x_j-y_j)(x_k-y_k)/4\mu}d\mu\\
 && \quad \qquad \quad\;\; + \int_0^\infty \frac{te^{-t^2/4\mu}}{\sqrt{4\pi \mu^3}} \frac{1}{(4\pi \mu)^{n/2}}e^{\sum_{j,k=1}^n h_{jk} [(b_k(y))(x_j-y_j)+(b_j(y))(x_k-y_k)]/4}d\mu \nonumber\\
 && \quad \qquad \quad \;\;+ \int_0^\infty \frac{te^{-t^2/4\mu}}{\sqrt{4\pi \mu^3}} \frac{1}{(4\pi \mu)^{n/2}}e^{-\sum_{j,k=1}^n h_{jk}(b_j(y))(b_k(y)) \mu/4}d\mu \nonumber\\
 &&\quad \qquad \quad \le \frac{\Gamma(\frac{n+1}{2})}{\pi^{\frac{n+1}{2}}}\left[ \frac{t}{\big(t^2+\sum_{j,k=1}^n h_{jk}(x_j-y_j)(x_k-y_k)\big)^{(n+1)/2}}\right. \nonumber \\&& \left.\quad \quad \qquad\;\; +t^{-n} e^{\sum_{j,k=1}^n
 h_{jk} [(b_k(y))(x_j-x_k)+(b_j(y))(x_k-y_k)]} +t^{-n}\frac{{}_{}}{{}}\right] \nonumber\end{eqnarray}
 and, for $m\ge 1$, \begin{eqnarray} \label{-b-d-2}
  &&\bigg|\int_0^\infty \frac{te^{-t^2/4\mu}}{\sqrt{4\pi \mu^3}}\left[G_0\# \left((\Delta_h- \Delta_h^0)G_0\right) \#\cdots \#  \left((\Delta_h- \Delta_h^0)G_0\right)(\mu, x,y)\right]d\mu\bigg|\\
  && \quad \qquad\le
  \int_0^\infty \frac{te^{-t^2/4\mu}}{\sqrt{4\pi \mu^3}}   \left[ \frac{c_2^m}{(m/2)!} \mu^{(m-n)/2} e^{-c_3 (d(x,y))^2/4\mu}\right] d\mu \nonumber\\
  &&\quad \qquad =
  \frac{c_2^m \Gamma(\frac{n-m+1}{2})}{2^{m-n}[(m/2)!]\sqrt{\pi}} \, \frac{t}{\big(t^2  +c_3 (d(x,y))^2\big)^{(n-m+1)/2}}.\nonumber\end{eqnarray}
  Therefore, (\ref{-b.21}) holds in the strong sense.   $\quad \quad \square$

   \vskip 0.28 true cm

\noindent{\bf Remark 5.2.}  \ \  We can also give another series representation for the heat kernel $F(t, x, y)$
of  the square root  $-\sqrt{-\Delta_h}$ of the Laplacian on $\partial \Omega$. In fact, from  Proposition 13.3 of p.$\,$62 of \cite{Ta2}, one has
that \begin{eqnarray*}G(t,x,y) \sim
 \sum_{j\ge 0}  t^{(j-n)/2}  p_j\big(x, t^{-1/2} (x-y)\big) e^{- \sum_{l,k=1}^n h_{lk}(x) (x_l-y_l)(x_k-y_k)/4t}
, \end{eqnarray*}
where $p_j(x,z)$ is a polynomial in $z$ which is even or odd in $z$ according to the parity of $j$ (here $p_\beta (x,z)$, a polynomial of degree $|\beta|$ in $z$,  is determined by the following formula:
  $$ \left[\mbox{det} \big(4\pi \big(h^{jk}(x)\big)\big) \right]^{-1/2} D_{z}^\beta e^{-\sum_{j,k=1}^n h_{jk}(x)\,z_j z_k/4} = p_\beta (x,z) e^{-\sum_{j,k=1}^n h_{jk}(x)\,z_j z_k/4}).$$
Therefore
\begin{eqnarray*} \label{b.21} F(t,x,y) \sim
    \sum_{j=0}^\infty\int_0^\infty \frac{t e^{-t^2/4\mu}}{\sqrt{4\pi \mu^3}}\left[  \mu^{(j-n)/2}  p_j\big(x, \mu^{-1/2} (x-y)\big) e^{ -\sum_{l,k=1}^n h_{lk} (x)\,(x_l-y_l)( x_k-y_k)/4\mu} \right] d\mu.\end{eqnarray*}

\vskip 0.2 true cm

\noindent{\bf Theorem 5.3.} \ {\it Let $\Omega$ be a bounded domain with smooth boundary $\partial \Omega$ in $(n+1)$-dimensional Riemannian manifold $(\mathcal{M},g)$. Then the fundamental solution ${\mathcal{K}} (t, x,y)$ of the heat equation for the Dirichlet-to-Neumann operator on $\partial \Omega$,  defined by (\ref{55-2}), has the following series representation:
\begin{eqnarray}\label{-b.31} {\mathcal{K}} (t,x,y)&=& {\mathcal{K}}_V (t,x,y) +\sum_{0\le l<M} {\mathcal{K}}_{V_{-1-l}} (t, x,y)+
{\mathcal{K}}_{V'_M} (t,x,y)\\ &=&
\int_0^\infty \frac{te^{-t^2/4\mu}}{\sqrt{4\pi \mu^3}}  \left[ \frac{{}^{}}{{}_{}}G_0(\mu,x,y)\right. \nonumber\\
  &&\left.+ \sum_{m=1}^\infty  G_0\# \underset{m}{\underbrace{\big((\Delta_h-\Delta_h^0)G_0\big) \#\cdots  \# \big((\Delta_h-\Delta_h^0)G_0\big)}}(\mu,x,y)\right]d\mu
\nonumber\\
    && +\sum_{0\le l<M} {\mathcal{K}}_{V_{-1-l}} (t, x,y)+
{\mathcal{K}}_{V'_M} (t,x,y),\nonumber \end{eqnarray} where
\begin{eqnarray} \label {b.5} \\
\left.\begin{array}{ll}   {\mathcal{K}}_{V} (t,x,y)= \left(\frac{1}{2\pi}\right)^{n} \int_{{\Bbb R}^n}
    e^{i(x-y)\cdot \xi} e^{tp(x,\xi)} d\xi =F(t,x,y),\nonumber\\
      {\mathcal{K}}_{V_{-2}} (t,x,y)=\left(\frac{1}{2\pi}\right)^{n} \int_{{\Bbb R}^n}
    e^{i(x-y)\cdot \xi} t (p_0^B(x,\xi)) e^{tp_1(x,\xi)} d\xi,\nonumber\\
       {\mathcal{K}}_{V_{-3}} (t,x,y)=\left(\frac{1}{2\pi}\right)^{n} \int_{{\Bbb R}^n}
    e^{i(x-y)\cdot \xi} \left[ t (p_{-1}^B(x,\xi))\right.\nonumber\\
    \;\, \quad \qquad \quad\qquad \left.+ \frac{t^2}{2} \left((p_0^B(x, \xi))^2 +2p_0(x, \xi) p_0^B (x,\xi)\right)\right] e^{tp_1(x,\xi)} d\xi, \nonumber\\
     {\mathcal{K}}_{V_{-4}} (t,x,y)=\left(\frac{1}{2\pi}\right)^{n} \int_{{\Bbb R}^n}
    e^{i(x-y)\cdot \xi}\left[ t(p_{-2}^B (x, \xi)) \right.\nonumber\\ \left. \qquad \qquad \qquad \;\,  +t^2 (p_0p_{-1}^B +p_{-1} p_0^B + p_0^B p_{-1}^B) + \frac{t^3}{6}
    \big((p_0^B)^3  + 3 p_0^2 p_0^B +3 p_0 (p_0^B)^2 \big) \right]e^{tp_1(x,\xi)} d\xi,\nonumber
     \\   {\mathcal{K}}_{V'_{M}} (t,x,y)=\left(\frac{1}{2\pi}\right)^{n} \int_{{\Bbb R}^n}
    e^{i(x-y)\cdot \xi}  v'_M (t,x,\xi)  d\xi,\end{array}\right. \nonumber\end{eqnarray}
      and
      \begin{eqnarray} \label{-b-d-3} \big| {\mathcal{K}}_{V'_M} (t,x,y)\big| \le c_0e^{-c_1 t} \left\{ \begin{array}{ll} t\big[t+d(x,y)\big]^{M-1-n} \quad \; &\mbox{if}\;\;  1-M>-n,\\ t(|\log (d(x,y)+t)|+1) \quad \; &\mbox{if}\;\; 1-M=-n,\\
      t \quad \; &\mbox{if}\;\; 1-M<-n,\end{array}\right.\end{eqnarray}
 for some $c_0>0$ and any $c_1 <0$, where $M=2,3,4$. }

 \vskip 0.28 true cm

 \noindent  {\it Proof.} \  Let $p(x,\xi)$ and $p^B(x, \xi)$ are the symbols of $-\sqrt{-\Delta_h}$ and $B$, respectively. Then \begin{eqnarray*}p(x, \xi)\sim \sum_{k\ge 0} p_{1-k} (x, \xi), \quad \;\;  p^B(x, \xi) \sim \sum_{k\ge 0}
 p^B_{-k} (x, \xi).\end{eqnarray*}
 Let us write \begin{eqnarray*} r (x, \xi)&=& p(x, \xi) -p_1(x,\xi)-p_0(x, \xi)-p_{-1}(x, \xi),\\
 r^B(x,\xi)&=& p^B (x, \xi) - p^B_0(x, \xi)- p^B_{-1} (x, \xi)-p^B_{-2} (x, \xi).\end{eqnarray*}
 It is clear that $r(x,\xi)\in S^{-2}_{1,0}$ and $r^B(x, \xi)\in S^{-3}_{1,0}$.
 Let $\mathcal{C}$ is a suitable curve in the complex plane going in the positive direction around the spectrum of $-{\mathcal{N}}_g$; it can be taken as the boundary $W_{r_0, \epsilon}$ for a small $\epsilon$ (see, Section 2) (here we can use ${\mathcal{C}}_{\theta,R}$ consisting of the two rays $re^{i\theta}$ and $re^{-i\theta}$, $\theta= \theta_0+\epsilon$ ($0<\theta<\frac{\pi}{2}$) by a circular piece
 in the right-plane with radius $R\ge 2t |p_1(x,\xi)|$).
 We may choose $R$ large enough. Then for $\lambda\in {\mathcal{C}}_{\theta,R}$, we can write the symbol $q(x, \xi, \lambda)$ of the resolvent operator $(-{\mathcal{N}}_g-\lambda)^{-1}$ in local coordinates as
  \begin{eqnarray*} &&\left[-p_1(x, \xi) -p_0(x, \xi) -p_{-1}(x, \xi) -r(x,\xi) - p^B_0 (x, \xi) -p^B_{-1} (x, \xi) -p^B_{-2} (x, \xi)-r^B(x,\xi) -\lambda\right]^{-1}
  \\  &&  =\left(-p_1-\lambda\right)^{-1} \left[ 1+ \frac{-p_0-p_{-1} -r - p^B_0 -p^B_{-1} -p_{-2}^B -r^B}{-p_1-\lambda}\right]^{-1}\\  && =
  \left(-p_1-\lambda\right)^{-1} \left[1-\frac{- p_0 -p_{-1}- r-p^B_0  -p^B_{-1} -p_{-2}^B- r^B}{-p_1-\lambda} \right.\\
   && \left.\; + \left(\frac{- p_0 -p_{-1}- r-p^B_0  - p^B_{-1} -p_{-2}^B- r^B}{-p_1-\lambda}\right)^2\right.\\
   &&\left. \;  - \left(\frac{- p_0 -p_{-1}- r- p^B_{0} -p^B_{-1}-p_{-2}^B-r^B}{-p_1-\lambda}\right)^3\right.
  \\ &&\left.\;+\left(\frac{p_0-p_{-1} -r -p_0^B-p_{-1}^B -p_{-2}^B -r^B}{-p_1 -\lambda}\right)^4 - \cdots  \right]\\
      &&  =\left[ \frac{1}{-p_1-\lambda} - \frac{-p_0 -p_{-1}-r}{(-p_1-\lambda)^2}
   +\frac{(-p_0-p_{-1}-r)^2}{(-p_1-\lambda)^3} -\frac{(-p_0-p_{-1}-r)^3}{(-p_1-\lambda)^4}+\cdots \right]\\
     \\ && \; +\left[ - \frac{-p^B_0}{(-p_1-\lambda)^2}+
      \left(  -\frac{-p^B_{-1}}{(-p_1-\lambda)^2}
 + \frac{(p_0^B)^2 +2p_0p_0^B}{(-p_1-\lambda)^3}\right) \right.\\
 && \left. \; + \left( -\frac{-p_{-2}^B}{(-p_1-\lambda)^2} + \frac{2(p_0 p_{-1}^B +p_{-1}p_0^B +p_0^B p_{-1}^B)}{(-p_1-\lambda)^3}\right.\right.\\
  && \left.\left. \; +\frac{3p_0^2 p_0^B + 3p_0 (p_0^B)^2 +(p_0^B)^3}{(-p_1-\lambda)^4}\right) +\cdots \right].
   \end{eqnarray*}
  In the above asymptotic expansion of symbol, the ``$\cdots$'' denotes all the terms which belong to $S_{1,0}^{-5}$.
  Noting that \begin{eqnarray*}  \frac{1}{-p-\lambda}=\frac{1}{-p_1-\lambda} - \frac{-p_0-p_{-1} -r}{(-p_1-\lambda)^2}
   +\frac{(-p_0-p_{-1} -r)^2}{(-p_1-\lambda)^3} - \frac{(-p_0-p_{-1}-r)^3}{(-p_1-\lambda)^4}+\cdots,\end{eqnarray*}
  we get \begin{eqnarray} \label{b-0'}   && q(x, \xi, \lambda) = \left(-p-\lambda\right)^{-1} +\left[  \frac{p_0^B}{(-p_1-\lambda)^2}\right]
  +\left[ \frac{p_{-1}^B}{(-p_1-\lambda)^2} + \frac{(p_0^B)^2 +2p_0p_0^B}{(-p_1-\lambda)^3}\right]
  \\ && \quad \;\; + \left( \frac{p_{-2}^B}{(-p_1-\lambda)^2} + \frac{2(p_0 p_{-1}^B +p_{-1}p_0^B +p_0^B p_{-1}^B)}{(-p_1-\lambda)^3}\right.\nonumber\\
  && \quad \;\;\left. \; +\frac{3p_0^2 p_0^B + 3p_0 (p_0^B)^2 +(p_0^B)^3}{(-p_1-\lambda)^4}\right) +r'\nonumber\\
  && \quad =(-p(x, \xi)-\lambda)^{-1} + q_{-2} (x, \xi, \lambda)+q_{-3}(x, \xi,\lambda) + q_{-4} (x, \xi, \lambda)+r'(x, \xi, \lambda),\nonumber
   \end{eqnarray}
     with \begin{eqnarray*} && q_{-2} (x, \xi, \lambda) = \frac{p_0^B}{(-p_1-\lambda)^2}, \quad \, q_{-3}(x, \xi, \lambda)=\frac{p_{-1}^B}{(-p_1-\lambda)^2} + \frac{(p_0^B)^2 +2p_0p_0^B}{(-p_1-\lambda)^3}, \\
    && q_{-4} (x, \xi, \lambda)=  \frac{p_{-2}^B}{(-p_1-\lambda)^2} + \frac{2(p_0 p_{-1}^B +p_{-1}p_0^B +p_0^B p_{-1}^B)}{(-p_1-\lambda)^3} \\
    &&\quad \qquad \qquad \quad +\frac{3p_0^2 p_0^B + 3p_0 (p_0^B)^2 +(p_0^B)^3}{(-p_1-\lambda)^4},
     \quad \quad  \; r'\in S^{-4}_{1,0}.\end{eqnarray*}
                     The semigroup $e^{-t(-{\mathcal{N}}_g)}$ can be defined from $-{\mathcal{N}}_g$ by the Cauchy integral formula
   \begin{eqnarray} \label{b-1'}e^{t{\mathcal N}_g}= e^{-t(-{\mathcal{N}}_g)} =\frac{i}{2\pi} \int_{{\mathcal{C}}_{\theta,R}} e^{-t\lambda} (-{\mathcal{N}}_g-\lambda)^{-1} d\lambda,\end{eqnarray}
   where ${\mathcal{C}}_{\theta, R}$ is the curve in complex plane defined as before.
 By virtue of (\ref{b-0'}) and  (\ref{b-1'}), we get that in the local coordinate patch the symbol of $e^{t{\mathcal{N}}_g}$ is
 \begin{eqnarray*}  v(t, x,\xi) \sim  e^{tp(x,\xi)} + \sum_{1\le l<M}v_{-1-l} (t,x,\xi) +v'_M(t,x,\xi), \end{eqnarray*}
 where \begin{eqnarray*} &&v_{-1-l}(t,x,\xi) =\frac{i}{2\pi} \int_{{\mathcal{C}}_{\theta, R}} e^{-t \lambda} q_{-1-l} (x, \xi,\lambda) d\lambda, \quad \, (l=1,2,\cdots, M-1),\\
 && v'_{M}(t,x,\xi) =\frac{i}{2\pi} \int_{{\mathcal{C}}_{\theta, R}} e^{-t \lambda}\, r'(x,\xi,\lambda)d\lambda, \quad \, M=2,3,4,\cdots.\end{eqnarray*}
  In order to give an explicit expansion of the heat kernel, we rewrite $q_{-2}$, $q_{-3}$  and $q_{-4}$  as
  \begin{eqnarray*} q_{-2} (x, \xi, \lambda) = \frac{b_{1,1}(x, \xi)}{(-p_1-\lambda)^2}, \quad \, q_{-1-l}(x,\xi,\lambda)=\sum_{k=1}^{2l} \frac{b_{l,k}(x,\xi)}{(-p_1-\lambda)^{k+1}},\, (l=2,3)\end{eqnarray*}
  where \begin{eqnarray*} &&b_{1.1} (x,\xi)= p_0^B(x,\xi), \quad \, b_{1,2}(x, \xi)=0, \quad \; b_{2,1}(x,\xi) = p_{-1}^B(x,\xi),\\
&&b_{2, 2} (x, \xi)= \left(p_0^B(x,\xi)\right)^2+ 2p_0(x,\xi)\, p_0^B(x, \xi), \quad \; b_{2,3}(x, \xi)=b_{2,4}(x, \xi)=0,\\
 && b_{3, 1} (x, \xi)= \frac{p_{-2}^B}{(-p_1-\lambda)^2}, \quad \; b_{3,2}(x,\xi)= \frac{2(p_0 p_{-1}^B +p_{-1}p_0^B +p_0^B p_{-1}^B)}{(-p_1-\lambda)^3},\\
 && b_{3, 3} (x, \xi)=\frac{3p_0^2 p_0^B + 3p_0 (p_0^B)^2 +(p_0^B)^3}{(-p_1-\lambda)^4}, \;\; b_{3, 4} (x, \xi)=b_{3, 5} (x, \xi)=b_{3, 6} (x, \xi)=0.\end{eqnarray*}
Then, for  $t>0$, we find by setting $\varrho=t\lambda$ that
 \begin{eqnarray*}
    \frac{i}{2\pi} \int_{{\mathcal{C}}_{\theta, R}} e^{-t\lambda} \frac{b_{l,k} (x,\xi)}{ (-p_1(x,\xi) -\lambda)^{k+1}} d\lambda  & =& \frac{i}{2\pi} \int_{{\mathcal{C}}_{\theta, R}} e^{-\varrho} \frac{t^kb_{l,k} (x,\xi)}{ (-tp_1(x,\xi) -\varrho)^{k+1}} d\varrho \\
    &=& \frac{1}{k!} t^k b_{l, k}(x,\xi)\, e^{tp_1(x, \xi)},\end{eqnarray*}
  Thus,  \begin{eqnarray} \label{08-20.1} v(t,x,\xi) &=& e^{tp(x, \xi)}, \quad \;\;  v_{-1-l} (t, x,\xi) = \sum_{k=1}^{2l}  \frac{1}{k!} t^k b_{l, k}(x,\xi)\, e^{tp_1(x,\xi)}\\
    && \qquad \; \mbox{for}\;\; l=1,2,\cdots, M-1. \nonumber\end{eqnarray}
   We define  $V(t)$, $V_{-1-l}(t)$ and $V'_M(t)$ in local coordinates to be the pseudodifferential operators with symbols $v(t,x,\xi)$, $v_{-1-l} (t,x,\xi)$ and $v'_M(t, x,\xi)$, respectively.
    It follows that the heat kernel ${\mathcal{K}} (t, x, \xi)$  of $e^{t{\mathcal{N}}_g}$ is in local coordinates expanded according to the symbol expansion:
   \begin{eqnarray*} {\mathcal{K}}(t, x, y) =  F(t,x,y) +
   \sum_{1\le l<M} {\mathcal{K}}_{V_{-1-l}} (t, x, y)
   + {\mathcal{K}}_{V'_M} (t, x, y),\end{eqnarray*}
   where  ${\mathcal{K}}_{V} (t,x,y)$,  ${\mathcal{K}}_{V_{-1-l}} (t, x, y)$ and  ${\mathcal{K}}_{V'_{M}} (t,x,y)$ are as in (\ref{b.5}) because of (\ref{08-20.1}), and $M=2,3,4$.
        Finally, in view of $\gamma(-{\mathcal{N}}_g)=0$,  we find  by (2.24) of Theorem 2.5 of \cite{GG} that
      \begin{eqnarray*} \label{-b-d} \big| {\mathcal{K}}_{V'_M} (t,x,y)\big| \le c_0e^{-c_1 t} \left\{ \begin{array}{ll} t\big[t+d(x,y)\big]^{M-1-n} \quad \; &\mbox{if}\;\;  1-M>-n,\\ t(|\log (d(x,y)+t)|+1) \quad \; &\mbox{if}\;\; 1-M=-n,\\
      t \quad \; &\mbox{if}\;\; 1-M<-n,\end{array}\right.\end{eqnarray*}
 for some $c_0>0$ and any $c_1 <0$.    $\qquad\quad \;\; \square$

\vskip 0.26 true cm

\vskip 1.39 true cm

\section{Estimates of the pole}

\vskip 0.45 true cm

Let $u_k$ be the eigenfunction corresponding to the $k$-th Steklov eigenvalue $\lambda_k$ (i.e., ${\mathcal{N}}_g u_k =-\lambda_k u_k$), and denote by $dS(x) =\sqrt{|h|}dx$ the area element of $\partial \Omega$.
Since $Z=\int_{\partial \Omega} {\mathcal{K}}(t, x, x)dS(x)$ converges, $e^{t {\mathcal{N}}_g} : \phi \to \int_{\partial \Omega}
 {\mathcal{K}}\phi$ is a compact mapping of the (real) Hilbert space $H=L^2 (\partial \Omega, \sqrt{|h|} \,dx)$.
This implies \begin{eqnarray} \label {413} {\mathcal{K}}(t,x,y)=\sum_{k\ge 1} e^{-\lambda_k t} u_k(x)  u_k(y) \end{eqnarray}
with uniform convergence on compact figures of
$(0, \infty)\times \partial \Omega\times \partial \Omega$,
and the spur $Z$ is easily evaluated as (see, for example,  \cite{Min} or Chapter 4 of \cite{Gr})
\begin{eqnarray}  Z= \int_{\partial \Omega} \sum_{k\ge 1} e^{-\lambda_k t} u_k^2(x)dS(x)
= \sum_{k\ge 1} e^{-\lambda_k t}.\end{eqnarray}
Formula (\ref{-b.31}) can now be used to estimate the pole ${\mathcal{K}} (t,x,x)$  for $t\downarrow 0$,  up to terms of
magnitude $t^{3-n}$, then integrate the result over $\partial \Omega$ to get an estimate of $Z=\mbox{Tr}\; e^{t{\mathcal{N}}_g}=\int_{\partial \Omega} {\mathcal{K}}(t,x,x)dS(x)$.

 In order to state our main result regarding the relationship between the spectrum of the Dircichlet-to-Neumann (or Steklov〞Poincar谷) operator and various informations on $\partial \Omega$, we introduce the following notations:
 We denote by $\,{\tilde R}_{jkjk}(x)$ (respectively $R_{jkjk}(x)$), $\,{\tilde{R}}_{jj}(x)=\sum_{k=1}^{n+1} {\tilde{R}}_{jkjk}(x)$ (respectively
 ${R}_{jj}(x)=\sum_{k=1}^{n} R_{jkjk}(x)$), $\, {\tilde {R}}_{\Omega}$ (respectively $R_{\partial \Omega} (x)$) the
 the curvature tensor, the Ricci curvature tensor, the scalar curvature with respect to $\Omega$ (respectively $\partial \Omega$) at $x\in \partial \Omega$. Denote by  $\,\sum_{j=1}^n {\tilde R}_{j(n+1)j(n+1), (n+1)} (x)$ the covariant derivative of the curvature tensor with respect to $(\bar \Omega, g)$. In addition, we denote by  $\kappa_1 (x), \cdots, \kappa_n(x)$  the principal curvatures of $\partial \Omega$ at $x\in \partial \Omega$.

\vskip 0.29 true cm

\noindent{\bf Theorem 6.1.} \ {\it Suppose that $(\mathcal{M},g)$ is an $(n+1)$-dimensional, smooth Riemannian
manifold, and assume that $\Omega\subset \mathcal{M}$ is a bounded domain
with smooth boundary $\partial \Omega$.  Let $\mathcal {K}(t, x, y)$ be the heat kernel associated to the Dirichlet-to-Neumann operator.

   (a)  \  If $n\ge 1$, then
   \begin{eqnarray} \label{6.0.1}    && \int_{\partial \Omega} \mathcal{K} (t, x,x)dx=t^{-n}
\int_{\partial \Omega} a_0(n,x) \,dS(x) + t^{1-n}  \int_{\partial \Omega} a_1(n,x)\, dS(x) \\  &&\quad \quad \quad \quad + \left\{\begin{array}{ll} O(t^{2-n}) \quad \;\, \mbox{when}\;\; n>1,\\ O(t\log t) \quad \, \mbox{when} \;\; n=1,\end{array}\right. \quad \;\;\mbox{as}\;\; t\to 0^+;\nonumber\end{eqnarray}

   (b) \  If $n\ge 2$, then \begin{eqnarray} \label{6.0.1-2}    && \int_{\partial \Omega} \mathcal{K} (t, x,x)dx=t^{-n}
\int_{\partial \Omega} a_0(n,x) \,dS(x) + t^{1-n}  \int_{\partial \Omega} a_1(n,x)\, dS(x)\\
  && + t^{2-n} \int_{\partial \Omega} a_2(n,x) \,dS(x) + \left\{\begin{array}{ll} O(t^{3-n}) \quad \, \mbox{when}\;\; n>2,\\ O(t\log t) \; \;\;\mbox{when} \;\; n=2,\end{array}\right.    \quad \;\;\mbox{as}\;\; t\to 0^+;\nonumber\end{eqnarray}

    (c) \  If $n\ge 3$, then
\begin{eqnarray} \label{6.0.1-3}    && \int_{\partial \Omega} \mathcal{K} (t, x,x)dx=t^{-n}
\int_{\partial \Omega} a_0(n,x) \,dS(x) + t^{1-n}  \int_{\partial \Omega} a_1(n,x)\, dS(x) \\
  && \qquad \qquad  + t^{2-n} \int_{\partial \Omega} a_2(n,x) \,dS(x)  + t^{3-n} \int_{\partial \Omega} a_3(n,x) \,dS(x) \nonumber\\&&\qquad \qquad +
  \left\{\begin{array}{ll} O(t^{4-n}) \quad\, \, \mbox{when}\;\; n>3,\\ O(t\log t) \quad \mbox{when} \;\; n=3,\end{array}\right.
   \quad \;\;\mbox{as}\;\; t\to 0^+.\nonumber\end{eqnarray} Here
   \begin{eqnarray} \label{06-060-1} a_0(n,x)=\frac{\Gamma(\frac{n+1}{2})}{\pi^{\frac{n+1}{2}}},\end{eqnarray}
  \begin{eqnarray}\label{06-060}\quad a_1(n,x)= \left(\frac{1}{2\pi}\right)^{n}   \frac{(n-1)\Gamma(n)\,\mbox{vol}({\Bbb S}^{n-1})}{2n} \bigg(\sum_{j=1}^n \kappa_j(x)\bigg),
     \end{eqnarray}
    \begin{eqnarray}  \label{66.55}&&    a_2(n,x) = \frac{\Gamma(n-1) \mbox{vol}({\Bbb S}^{n-1}) }{8(2\pi)^n} \left[ \frac{3-n}{3n} R_{\partial \Omega} + \frac{n-1}{n} {\tilde R}_{\Omega}\right. \\ && \left. +\frac{n^3 -n^2 -4n +6}{n(n+2)} \big(\sum_{j=1}^n \kappa_j (x)\big)^2 + \frac{n^2 -n-2}{n(n+2)} \sum_{j=1}^n \kappa_j^2(x) \right]
      \nonumber \end{eqnarray} and
\begin{eqnarray} \label{6.6-1.} && a_3(n,x)=
  \bigg(\frac{1}{2\pi}\bigg)^{n}\frac{\Gamma(n-2)\,\mbox{vol} (S^{n-1}) }{8n}  \bigg[
   \frac{n^3 -2n^2 -7n+7}{2(n+2)} {\tilde R}_\Omega(x)\, \big(\sum_{j=1}^n \kappa_j(x)\big) \\ && \quad \, + \frac{-3n^4 -4 n^3 +59n^2 +75 n -180}{ 6(n+2)(n+4)} \big(\sum_{j=1}^n \kappa_j (x) \big)R_{\partial \Omega} \nonumber\\ && \quad \,+
   \frac{n^5 -20 n^3 +2n^2 +61 n -74}{6(n+2)(n+4)} \big(\sum_{j=1}^n \kappa_j(x)\big)^3 \nonumber\\ && \quad \,
    + \frac{n^4 +8n^3 +15n^2 +3n -32}{2(n+2) (n+4)} \big(\sum_{j=1}^n \kappa_j(x)\big)\big(\sum_{j=1}^n \kappa_j^2 (x)\big)\nonumber\\ && \quad \,
    +\frac{-6n^3-34n^2 +40}{3(n+2)(n+4)} \sum_{j=1}^n \kappa_j^3(x) + \frac{4n^2-6}{n+2} \sum_{j=1}^n\kappa_j(x) {\tilde R}_{jj}(x)\nonumber
    \\ && \quad \,-\frac{12n^3 +50 n^2 -6n -104}{3(n+2)(n+4)} \sum_{j=1}^n \kappa_j(x) R_{jj}(x) + (n-1) \sum_{j=1}^n {\tilde{R}}_{j(n+1)j(n+1),(n+1)}(x)\nonumber\\ &&\quad \,-\frac{n-2}{2 } {\tilde{R}}_{\Omega} (x)  +\frac{n-2}{2} R_{\partial \Omega} -\frac{n-2}{2} \sum_{j=1}^n \kappa_j^2(x)\bigg]\nonumber
.\end{eqnarray}}

\vskip 0.2 true cm

  \noindent  {\bf Proof.}    Put $x_0=0$ for simplicity, and (as in proof of Lemma 4.2) choose coordinates $x=(x_1,\cdots, x_n)$ on an open set in $\partial \Omega$ (centered at $x_0=0$) and then
 coordinates $(x, x_{n+1})$ on a neighborhood in $\bar\Omega$ such that
$x_{n+1}=0$ on $\partial \Omega$ and $|\nabla x_{n+1}|=1$ near $\partial \Omega$ while $x_{n+1}>0$ on $\Omega$ and such that $x$ is constant
on each geodesic segment in $\bar\Omega$ normal to $\partial \Omega$.  Clearly, the metric tensor on $\bar \Omega$ has
the form
  \begin{gather*} \label{11} \big(g_{jk} (x,x_{n+1}) \big)_{(n+1)\times (n+1)} =\begin{pmatrix} ( h_{jk} (x,x_{n+1}))_{n\times n}& 0\\
      0& 1 \end{pmatrix}. \end{gather*}
     Recall that $\left(-\frac{1}{2}\, \frac{\partial h_{jk}}{\partial x_{n+1}}(x_0)\right)_{n\times n}$  is the matrix of the Weingarten map under a basis of $T_0(\partial \Omega)$ at $x_0=0$, its eigenvalues are just the principal curvatures  $\kappa_1, \cdots, \kappa_n$   of $\partial \Omega$ at $x_0=0$ in the direction $\nu$. Therefore we can choose the principal curvature vectors $e_1, \cdots, e_n$ as an orthonormal basis of $T_0(\partial \Omega)$ at $x_0=0$  such that the symmetric matrix $\left(-\frac{1}{2}\, \frac{\partial h_{jk}}{\partial x_{n+1}}(x_0)\right)_{n\times n}$  becomes the diagonal matrix
      \begin{gather*}  \begin{pmatrix} \kappa_1  & 0  & \cdots  & 0 \\
     0  & \kappa_2  & \cdots & 0  \\
      \vdots & \vdots & \ddots & \vdots  \\
     0   & 0 & \cdots &\kappa_n
       \end{pmatrix}. \end{gather*}
          Furthermore, we take a geodesic normal coordinate system for $(\partial \Omega, h)$ centered at $x_0=0$, with respect to $e_1, \cdots, e_n$.
           As Riemann showed, one has (see p.$\,$555 of \cite{Ta2}, or \cite{Spi2})
            \begin{eqnarray} \label{7/14/1}&& h_{jk}(x_0)= \delta_{jk}, \; \; \frac{\partial h_{jk}}{\partial x_l}(x_0)
 =0, \;\; \frac{\partial^2 h_{jk}}{\partial x_l \partial x_m} (x_0) =-\frac{1}{3} R_{jlkm} -\frac{1}{3} R_{jmkl}\\
 && \qquad \qquad \mbox{for all} \;\; 1\le j,k,l,m \le n,\nonumber\end{eqnarray}
   so that \begin{eqnarray} \label {08-18.1} \frac{\partial^2 |h|}{\partial x_k \partial x_l} (x_0) =\sum_{m=1}^n \frac{\partial^2 h_{mm}}{\partial x_k\partial x_l} (x_0) =- \sum_{m=1}^n  \frac{R_{mkml}+ R_{mlmk}}{3}= -\sum_{m=1}^n \frac{2 R_{mkml}}{3}.\end{eqnarray}
   Also, from $hh^{-1}=I$ we get \begin{eqnarray}\label{7/14/2}
  \frac{\partial h^{jk}}{\partial x_l} (x_0)=0,\quad \frac{\partial^2 h^{jk}}{\partial x_l\partial x_m} (x_0) = -\frac{\partial^2 h_{jk}}{\partial x_l\partial x_m} (x_0)\quad\;\mbox{for all} \,\, 1\le j,k,l,m\le n.\end{eqnarray}

 In view of (\ref{-b-d-2}) and (\ref{-b-d-3}), we find by taking $x=y$ that  \begin{eqnarray} \label{+61}|{\mathcal{K}}_{V'_M} (t, x,x)| \le c_0\left\{ \begin{array}{ll} t^{M-n},  \quad \; & \mbox{if}\;\; 1-M>-n,\\
  t(1+\log t),\quad \; &\mbox{if}\;\;  1-M =-n,\\
   t, \quad\; &\mbox{if}\;\; 1-M<-n\end{array} \right.\end{eqnarray}
 and,   for $m\ge 1$, \begin{eqnarray} \label{-d-b-9}
  &&\bigg|\int_0^\infty \frac{te^{-t^2/4\mu}}{\sqrt{4\pi \mu^3}}\bigg[G_0\# \underset{m}{\underbrace{\left((\Delta_g- \Delta_g^0)G_0\right) \#\cdots \#  \left((\Delta_g- \Delta_g^0)G_0\right)}}(\mu, x,x)\bigg]d\mu\bigg|\\
  && \quad \qquad \quad\;\; \le   \frac{c_2^m \Gamma(\frac{n-m+1}{2})\,t^{m-n} }{2^{m-n}[(m/2)!]\sqrt{\pi}},\nonumber\end{eqnarray}
  where $M=2,3,4$ will be taken later (according to the term number of asymptotic expansion of the heat kernel). Consequently, by (\ref{-b.31}), the trace of the operator $e^{t {\mathcal{N}}_g}$ has the following asymptotic expansion
\begin{eqnarray}  \label {b.19}&&\mbox{Tr}\, e^{t {\mathcal{N}}_g} =
 \int_{\partial \Omega}\left[ \int_0^\infty \frac{te^{-t^2/4\mu}}{\sqrt{4\pi \mu^3}}  \left(\frac{{}_{}}{{}}G_0(\mu,x,x)+  \big(G_0\#\big((\Delta_h-\Delta_h^0)G_0\big)\big)(\mu,x,x)\right)d\mu \right.\\&& \qquad \qquad \quad \; \left.
+\sum_{1\le l<M} {\mathcal{K}}_{V_{-1-l}} (t, x,x)\right]dS(x) +{\mathcal{K}}_{V'_M} (t,x,x) \quad \mbox{as}\;\; t\to 0^+,\nonumber\end{eqnarray}
 i.e.,
 up to terms of magnitude $\le$ constant$\times t^{M-n}$, one is left with
 \begin{eqnarray}\label{-b-d-11} &&\int_{\partial \Omega}  {\mathcal{K}} (t,x,x) dS(x) =
  \int_{\partial \Omega}\left[ \int_0^\infty \frac{te^{-t^2/4\mu}}{\sqrt{4\pi \mu^3}}  \left(\frac{{}_{}}{{}_{}}G_0(\mu,x,x)\right.\right.\\
   &&\left.\left. \quad \quad +  \big(G_0\#\big((\Delta_h-\Delta_h^0)G_0\big)\big(\mu,x,x)\frac{{}_{}}{{}_{}}\right)d\mu
+\sum_{1\le l<M} {\mathcal{K}}_{V_{-1-l}} (t, x,x)\right]dS(x).\nonumber\end{eqnarray}
  According to the definition of $\Delta_h$, we get that at the fixed point $x_0=0$, \begin{eqnarray*}b_k(x_0)=\sum_{j=1}^n \frac{1}{\sqrt{|h|}} \, \frac{\partial (\sqrt{|h|}\, h^{jk})
}{\partial x_j}(x_0) =0,\end{eqnarray*} so that  \begin{eqnarray} \label{-b-d-12} \int_0^\infty \frac{te^{-t^2/4\mu}}{\sqrt{4\pi \mu^3}}\, G_0(\mu, x,x)d\mu=
\int_0^\infty \frac{te^{-t^2/4\mu}}{\sqrt{4\pi \mu^3}}\left(\frac{1}{4\pi \mu}\right)^{\frac{n}{2}}d\mu =\frac{\Gamma\big(\frac{n+1}{2}\big)
}{\pi^{\frac{n+1}{2}}}\,t^{-n}.\end{eqnarray}
  It follows from \cite{MS} that up to the desired precision,
the integrand $$G_0 (t-s, 0, x) \big((\Delta_h -\Delta_h^0) G_0 (s,x, 0)\big)\sqrt{\mbox{det}\; h}$$
can be replaced by the product of a factor
 $1+$ a line function $f$ of $x+o(|x|^2)$  and expression
 \begin{eqnarray} \label {j1}  && \frac{e^{-|x|^2/4(t-s)}}{[4\pi(t-s)]^{n/2}} \left[ \frac{1}{2} \frac{\partial^2 h^{ij}}{\partial x_k\partial x_l}(0) \, x_kx_l \frac{\partial^2}{\partial x_i\partial x_j}\right.\\
 && \left.\quad \quad \quad\;\;  +\left(\frac{\partial}{\partial x_k} \, \frac{1}{\sqrt{|h| }} \frac{\partial}{\partial x_i}\sqrt{|h|}\, h^{ij} \right)(0)
 x_k \frac{\partial}{\partial x_j}\right]\frac{e^{-|x|^2/4s}}{(4\pi s)^{n/2}}\nonumber\\
 && \quad \quad\;\; = (4\pi t)^{-n/2}\frac{e^{-|x|^2/4r}}{(4\pi r)^{n/2}} \left[ \frac{1}{2}\, \frac{\partial^2 h^{ij}}{\partial x_k\partial x_l}(0) x_kx_l\left(\frac{x_ix_j}{4s^2} -\frac{\delta_{ij}}{2s} \right)\right.\nonumber\\
  && \left.\quad \quad \quad \;\; -\left( \frac{\partial}{\partial x_k} \, \frac{1}{\sqrt{|h|} } \, \frac{\partial}{\partial x_i} \sqrt{|h|} \,h^{ij} \right)(0)
  \frac{x_kx_j}{2s} \right],\nonumber\end{eqnarray}
where $r= s(t-s)/t$. Now the factor alluded to above (\ref{j1}) can be replaced by $1$,
since $f\times$ (\ref{j1}) integrates to $0$
 while the last $2$ terms contribute $\le c t^{3-n}$.
     Consequently,  up to the desired precision, we find by p.$\,$52--53 of \cite{MS} and (\ref{7/14/1})--(\ref{7/14/2}) that at $x_0=0$,
\begin{eqnarray*}   && G_0\# \big((\Delta_h-\Delta^0_h)G_0\big) =\int_0^t ds \int_{{\Bbb R}^n} (\ref{j1}) \,dx \\
 && \quad \;\; = \frac{t}{(4\pi t)^{n/2}} \left[
-\frac{1}{6} \sum_{i,j=1}^n \frac{\partial^2 h^{ij}}{\partial x_i\partial x_j} - \frac{1}{2} \sum_{i=1}^n \frac{\partial^2 \sqrt{|h(x)|}}{\partial x_i^2} -\frac{1}{12} \sum_{i,j=1}^n  \frac{\partial^2 h^{ii}}{\partial x_j^2}
\right]\nonumber\\
&& \quad \;\; = \frac{t}{(4\pi t)^{n/2}}\left[
\frac{1}{6} \sum_{i,j=1}^n \frac{\partial^2 h_{ij}}{\partial x_i\partial x_j} - \frac{1}{4} \sum_{i=1}^n \frac{\partial^2 |h(x)|}{\partial x_i^2} +\frac{1}{12} \sum_{i,j=1}^n  \frac{\partial^2 h_{ii}}{\partial x_j^2}
\right]\nonumber\\
&&\quad \;\; = \frac{t}{(4\pi t)^{n/2}}  \sum_{i,j=1}^n \bigg[ -\frac{1}{18} R_{ijji} +\frac{1}{6} R_{jiji}  - \frac{1}{18}\, R_{ijij} \bigg]
\nonumber\\         &&\quad \;\; =  \frac{t}{(4\pi t)^{n/2}} \left(\frac{1}{6} \sum_{i,j=1}^n R_{jiji}\right) =\frac{t}{(4\pi t)^{n/2}}\,\frac{R_{\partial \Omega}}{6},\nonumber\end{eqnarray*}
    so that   \begin{eqnarray}   \label{j2} \quad \quad\quad  \int_0^\infty \frac{te^{-t^2/4\mu}}{\sqrt{4\pi \mu^3}} \left[\big(G_0\# \big((\Delta_g-\Delta_g^0)G_0\big)(\mu,x,x)\big)\right]d\mu &=& \frac{R_{\partial \Omega}(x)}{6} \int_0^\infty \frac{te^{-t^2/4\mu}}{\sqrt{4\pi \mu^3}} \, \frac{\mu}{(4\pi \mu)^{n/2}} d\mu
  \\  &=& \frac{\Gamma(\frac{n-1}{2})\,R_{\partial \Omega}(x)}{24\pi^{\frac{n+1}{2}}}\, t^{2-n}, \nonumber\end{eqnarray}
  where $R_{\partial \Omega} (x)$ is the scalar curvature of $(\partial \Omega, h)$. It is easy to calculate that $R_{\partial \Omega} = {\tilde{R}} -2{\tilde R}_{\nu\nu} +\frac{1}{2}Q_1$, where $Q_1=4\sum_{1\le j<k\le n} \kappa_j\kappa_k$ (see (\ref{6601})).
In the normal coordinate system on $\partial \Omega$ with respect to the principal curvature vectors $e_1, \cdots, e_n$ centered at $x_0=0$, the Weingarten map  $A_\nu$ is self-adjoint,  and (see p.$\,$159 of \cite{Jo})
\begin{eqnarray} \mbox{Tr}\, A_\nu = \sum_{j=1}^n \kappa_j, \quad \; \, \frac{\langle A^*_\nu \xi, \xi\rangle}{\langle\xi, \xi\rangle} =\frac{\sum_{j=1}^n \kappa_j\xi_j^2}{\sum_{j=1}^n \xi_j^2}.\end{eqnarray}
 By (\ref{b.40}), (\ref{y6}) and  (\ref{4-080}), we can easily see that
\begin{eqnarray*}
 p_0^B(x,\xi)= \frac{1}{2} \left( \sum_{j=1}^n\kappa_j(x) - \frac{\sum_{j=1}^n \kappa_j(x)\xi_j^2}{\sum_{j=1}^n \xi_j^2}\right),
 \end{eqnarray*}
\begin{eqnarray*}\label{+y6} p_{-1}^B(x, \xi)&=&\frac{1}{8\sqrt{\sum_{j=1}^n \xi_j^2}} \left[-2 Q_1 +2 \sum_{j=1}^n {\tilde R}_{j(n+1)j (n+1)}
 - 2 \sum_{j=1}^n \kappa_j^2(x)  +3 \bigg(\sum_{j=1}^n \kappa_j(x)\bigg)^2 \right.\\
   && \left. +\frac{5\big(\sum_{j=1}^n \kappa_j(x) \xi_j^2\big)^2}{\big(\sum_{j=1}^n \xi_j^2\big)^2} -\frac{1}{\sum_{j=1}^n \xi_j^2}\, \sum_{j,k=1}^n \bigg(2{\tilde{R}}_{j(n+1)k(n+1)} + 6
   \kappa_j^2(x) \delta_{jk}\bigg) \xi_j\xi_k\right], \end{eqnarray*}
   \begin{eqnarray*} \label{4-080'}   p_{-2}^B (x, \xi)= \frac{\sum_{j=1}^n \kappa_j(x)\xi_j^2}{8(\sum_{j=1}^n \xi_j^2)^2}   \bigg( 2Q_1 -2 \sum_{j=1}^n {\tilde R}_{j(n+1)j(n+1)} + 2 \sum_{j=1}^n \kappa_j^2(x)  -3 \big(\sum_{j=1}^n \kappa_j(x)\big)^2
    \end{eqnarray*}
\begin{displaymath}  \quad \quad\;\; -\frac{5\big(\sum_{j=1}^n \kappa_j(x) \xi_j^2\big)^2}{\big(\sum_{j=1}^n \xi_j^2\big)^2}
  +\frac{1}{\sum_{j=1}^n \xi_j^2}\, \sum_{j,k=1}^n \big(2{\tilde{R}}_{j(n+1)k(n+1)} + 6
   \kappa_j^2(x) \delta_{jk}\big) \xi_j\xi_k
     \bigg)\qquad \qquad \qquad \quad
    \nonumber \\ \end{displaymath}
     \begin{eqnarray*} &&  - \frac{1}{4 \sum_{j=1}^n \xi_j^2}\, \left\{ \frac{3}{2} (\sum_{j=1}^n \kappa_j(x) ) \big( 2Q_1 -2 \sum_{j=1}^n {\tilde R}_{j(n+1)j(n+1)} + 2 \sum_{j=1}^n \kappa_j^2(x)
      \big)-4(\sum_{j=1}^n \kappa_j(x))^3 \right. \\
      && \left.
     +\frac{1}{4} \bigg(-48\sum_{1\le j<k<l\le n}\kappa_j\kappa_k\kappa_l \, - 6\sum_{j\ne k} \kappa_j \big(-2{\tilde R}_{k(n+1)k(n+1)}\right.\end{eqnarray*} \begin{eqnarray*}
     &&\left.   +2\kappa_k^2\big)- \sum_{j=1}^{n} \big(2{\tilde R}_{j(n+1)j(n+1), (n+1)} - 4 \sum_{j=1}^n {\tilde R}_{j(n+1)j(n+1)} \kappa_j \big)  \bigg)\right.\\
      &&  \left. -\frac{1}{2} \bigg[\frac{\sum_{j=1}^n \kappa_j(x)\xi_j^2}{\sum_{j=1}^n \xi_j^2}
      +\sum_{j=1}^n\kappa_j(x)
      \bigg] \bigg[ \frac{1}{2 \sum_{j=1}^n \xi_j^2} \bigg(\sum_{j,k=1}^n \big(2 {\tilde R}_{j(n+1)k(n+1)}+6 \kappa_j^2\delta_{jk} \big)  \xi_j\xi_k \bigg)\right. \;\;\end{eqnarray*} \begin{eqnarray*}
      &&  \left.-\frac{2\big(\sum_{j=1}^n \kappa_j(x)\xi_j^2\big)^2}{\big(\sum_{j=1}^n \xi_j^2\big)^2 } +2 \bigg(\sum_{j=1}^n \kappa_j(x)
      \bigg)^2 -\frac{1}{2}
      \bigg(2Q_1- 2\sum_{j=1}^n {\tilde R}_{j(n+1)j(n+1)} +2 \sum_{j=1}^n \kappa_j^2(x)
      \bigg)\bigg]\right.\\
      && \left. +\frac{1}{4 \sum_{j=1}^n \xi_j^2}
            \bigg[\sum_{j,k=1}^n \bigg(2{\tilde R}_{j(n+1)k(n+1),(n+1)} +20 {\tilde R}_{j(n+1)k(n+1)}\kappa_{{}_k} +24 \kappa_j^3\delta_{jk}\bigg) \xi_j\xi_k\bigg]\right.\end{eqnarray*} \begin{eqnarray*}
             && \left.
            -\frac{3\sum_{j=1}^n \kappa_j(x)\xi_j^2
            }{2\big(\sum_{j=1}^n \xi_j^2\big)^2} \bigg[ \sum_{j,k=1}^n\bigg(2 {\tilde R}_{j(n+1)k(n+1)}
        +6 \kappa_j^2\delta_{jk}\bigg)\xi_j\xi_k\bigg]  \right.\\   &&\left. +  \frac{4
         \big(\sum_{j=1}^n \kappa_j(x)\xi_j^2\big)^3
         }{\big(\sum_{j=1}^n \xi_j^2\big)^3} + \frac{\sum_{j=1}^n \kappa_j(x)
         }{2} \bigg[\frac{1}{2\sum_{j=1}^n \xi_j^2}\bigg( \sum_{j,k=1}^n
      \big( 2{\tilde R}_{j(n+1)k(n+1)} +6 \kappa_j^2\delta_{jk} \big)\xi_j\xi_k\bigg) \right.\end{eqnarray*} \begin{eqnarray*}
      && \left.-\frac{2 \big(\sum_{j=1}^n \kappa_j(x)\xi_j^2\big)^2
      }{\big(\sum_{j=1}^n \xi_j^2\big)^2} +2 \bigg(\sum_{j=1}^n \kappa_j(x)
      \bigg)^2 - \frac{1}{2} \bigg (2Q_1-2 \sum_{j=1}^n {\tilde R}_{j(n+1)j(n+1)} + 2\sum_{j=1}^n \kappa_j^2(x)
      \bigg)\bigg] \right.\\
       &&  \left.-\frac{1}{4} \bigg[\frac{\sum_{j=1}^n \kappa_j(x)\xi_j^2}{\sum_{j=1}^n \xi_j^2}
      +\sum_{j=1}^n\kappa_j(x)
       \bigg] \bigg[ 2Q_1 -2\sum_{j=1}^n
       {\tilde R}_{j(n+1)j(n+1)} + 2 \sum_{j=1}^n \kappa_j^2(x)
               -4(\sum_{j=1}^n \kappa_j(x))^2 \bigg] \right\} \end{eqnarray*}
         Also, it follows from  (\ref{-5b.1'}), (\ref{-5b.2'}), (\ref{7/14/1})--(\ref{7/14/2}) that $p_0(x,\xi)$ and  $p_{-1} (x,\xi)$ have the following simple expressions in the  normal coordinates chosen before:
         \begin{eqnarray} \label{6-603;-} p_0(x,\xi)= 0,\end{eqnarray}
         \begin{eqnarray} \label{6-603;}  &&p_{-1}(x,\xi) = \frac{1}{4\big(\sum_{j=1}^n \xi_j^2 \big)^{\frac{3}{2}}} \bigg[ -\frac{1}{2}
\sum_{k,l=1}^n \frac{\partial^2 |h(x)|}{\partial x_k\partial x_l}\, \xi_k \xi_l - \sum_{j,k,l=1}^n \frac{\partial^2 h^{jk}}{\partial x_k \partial x_l}\, \xi_j\xi_l \quad\quad\\
&&\qquad\quad + \frac{1}{\sum_{j=1}^n \xi_j^2} \sum_{j,k,l,m=1}^n \frac{\partial^2 h^{jk}}{\partial x_m\partial x_l}\, \xi_j\xi_k\xi_l\xi_m -\frac{1}{2}\sum_{j,k,m=1}^n \frac{\partial^2 h^{jk}}{\partial x_m^2}\, \xi_j\xi_k\bigg] \quad \qquad\nonumber \end{eqnarray}
\begin{eqnarray} && \quad \;\; = \frac{1}{4\big(\sum_{j=1}^n \xi_j^2 \big)^{\frac{3}{2}}} \bigg[ \sum_{k,l,m=1}^n  \frac{R_{mkml}}{3} \xi_k\xi_l +\sum_{j,k,l=1}^n
  \frac{R_{jklk}}{3}\xi_j\xi_l\nonumber\\ && \;\; \quad\quad + \frac{1}{\sum_{j=1}^n \xi_j^2}\sum_{j,k,l,m=1}^n \frac{R_{jmkl} +  R_{jlkm}}{3} \xi_j\xi_k\xi_l\xi_m -
 \sum_{j,k,m=1}^n  \frac{R_{jmkm}}{3} \xi_j\xi_k \bigg]\;\nonumber\\
 &&\;\,\quad= \frac{1}{4\big(\sum_{j=1}^n \xi_j^2 \big)^{\frac{3}{2}}} \bigg[ \sum_{j,k,l=1}^n
  \frac{R_{jklk}}{3}\xi_j\xi_l\nonumber\\ &&\;\,\quad\quad + \frac{1}{\sum_{j=1}^n \xi_j^2}\sum_{j,k,l,m=1}^n \frac{R_{jmkl} +  R_{jlkm}}{3} \xi_j\xi_k\xi_l\xi_m  \bigg].
 \nonumber
\end{eqnarray}
              Therefore \begin{eqnarray} \label{6-602} \quad \;\quad\quad\; {\mathcal{K}}_{V_{-2}} (t,x,x)= \left(\frac{1}{2\pi}\right)^{n}  \int_{{\Bbb R}^n} \frac{t}{2} \left[\sum_{j=1}^n\kappa_j(x) -\frac{\sum_{j=1}^n \kappa_j(x)\xi_j^2}{\sum_{j=1}^n \xi_j^2} \right]e^{-t\, \sqrt{\sum_{j=1}^n \xi_j^2}}\,d\xi.\end{eqnarray}
         Next, from (\ref{-5b.1'}) and (\ref{b.40}), we see that $$ p_0^B (x, -\xi) =p_0^B (x, \xi), \quad \; p^B_{-1}(x,-\xi)= p_{-1}^B(x,\xi), \;\; \mbox{for all}\;\;  \xi\in {\Bbb R}^n,$$
 so that
 \begin{eqnarray}\label{6-603} &&{\mathcal{K}}_{V_{-3}} (t,x,x)= \big(\frac{1}{2\pi}\big)^{n} \int_{{\Bbb R}^n}
 \big[ tp^B_{-1}(x,\xi) +\frac{t^2}{2}\big(\big(p_0^B (x,\xi)\big)^2
  +2p_0(x,\xi) p_0^B(x,\xi)\big)\big] e^{tp_1(x,\xi)}d\xi \\ && \quad \;\;=\left(\frac{1}{2\pi}\right)^{\frac{n}{2}} \int_{{\Bbb R}^n}
 \left[ tp^B_{-1}(x,\xi) +\frac{t^2}{2}\big(p_0^B (x,\xi)\big)^2 \right] e^{-t\sqrt{\sum_{j=1}^n \xi_j^2}}d\xi\nonumber \\
  && \quad \;\; =
 \left(\frac{1}{2\pi}\right)^{\frac{n}{2}} \int_{{\Bbb R}^n}
    \left\{\frac{t}{8\big(\sum_{j=1}^n \xi_j^2\big)^{\frac{1}{2}}}\left[-2Q_1 +2 \sum_{j=1}^n {\tilde R}_{j(n+1)j (n+1)}(x)
    \right.\right.\nonumber
\\ &&\quad\;\;\left.\left.
    - 2 \sum_{j=1}^n \kappa_j^2(x)  +3\bigg( \sum_{j=1}^n \kappa_j(x)\bigg)^2 +\frac{5\big(\sum_{j=1}^n \kappa_j(x)\, \xi_j^2\big)^2}{\big(\sum_{j=1}^n \xi_j^2\big)^2}
     \right.\right.\nonumber
\\ &&\quad\;\;\left.\left.
     -\,\frac{1}{\sum_{j=1}^n \xi_j^2} \,\sum_{j=1}^n \bigg(2{\tilde{R}}_{j(n+1)j(n+1)} +6
   \kappa_j^2(x) \bigg) \xi_j^2\right] \right.\nonumber\\
  &&\quad \;\; \left.\, +\,\frac{t^2}{8} \bigg(\sum_{j=1}^n \kappa_j(x) - \frac{\sum_{j=1}^n \kappa_j(x)\,\xi_j^2}{\sum_{j=1}^n \xi_j^2}\bigg)^2 \right\}e^{-t\, \sqrt{\sum_{j=1}^n \xi_j^2}}\,d\xi,\nonumber
    \end{eqnarray}  and \begin{eqnarray}\label{6.-18} && {\mathcal{K}}_{V_{-4}} (t,x,x)= \left(\frac{1}{2\pi}\right)^{n} \int_{{\Bbb R}^n}
    \bigg[ t\,p_{-2}^B (x, \xi)  +t^2 \bigg(p_0(x, \xi)p_{-1}^B(x, \xi)\\ && \quad \;\;\quad +p_{-1}(x, \xi) p_0^B(x, \xi)  + p_0^B(x, \xi) p_{-1}^B(x, \xi)\bigg)  + \frac{t^3}{6}
    \bigg((p_0^B(x, \xi))^3 \nonumber\\ && \quad \;\;\quad + 3 (p_0(x, \xi))^2 p_0^B(x, \xi) +3 p_0(x, \xi) (p_0^B(x, \xi))^2 \bigg)      \bigg]e^{-t\, \sqrt{\sum_{j=1}^n \xi_j^2}}\,d\xi.\nonumber
   \\ && \quad \quad = \left(\frac{1}{2\pi}\right)^{\frac{n}{2}} \int_{{\Bbb R}^n}
    \bigg[ t\,p_{-2}^B (x, \xi)  +t^2 \bigg(p_{-1}(x, \xi) p_0^B(x, \xi) + p_0^B(x, \xi) p_{-1}^B(x, \xi)\bigg)  \nonumber
   \\ && \quad \;\;\quad + \frac{t^3}{6}
    \big(p_0^B(x, \xi)\big)^3  \bigg]e^{-t\, \sqrt{\sum_{j=1}^n \xi_j^2}}\,d\xi.\nonumber
    \end{eqnarray}
          By applying the spherical coordinates transform
  \begin{eqnarray*} \left\{ \begin{array}{ll} \xi_1= r\cos \phi_1\\
  \xi_2 = r\sin \phi_1 \cos \phi_2\\
  \xi_3 =r\sin \phi_1 \sin \phi_2 \cos \phi_3\\
  \cdots \cdots \cdots \cdots\\
   \xi_{n-1} = r\sin \phi_1 \sin \phi_2 \cdots \sin \phi_{n-2} \cos \phi_{n-1} \\
   \xi_{n} = r\sin \phi_1 \sin \phi_2 \cdots \sin \phi_{n-2} \sin \phi_{n-1},\end{array} \right.\end{eqnarray*}
  where $0\le r<+\infty, \, 0\le \phi_1\le \pi, \,\cdots,\, 0\le \phi_{n-2} \le \pi, \;  0\le \phi_{n-1} \le 2\pi$,
   we then have $d\xi =r^{n-1} dS=r^{n-1} \sin^{n-2} \phi_1 \sin^{n-3} \phi_2 \cdots \sin \phi_{n-2}.$
   Since   \begin{eqnarray*} \int_0^\pi \sin^{k} \phi \;d\phi = \frac{\Gamma(\frac{n+1}{2})}{\Gamma(\frac{n+2}{2})}\,\sqrt{\pi}= \left\{ \begin{array}{ll}\frac{(2m-1)!!}{(2m)!!}\, \pi, \quad \, k=2m,\\
      2  \frac{(2m)!!}{(2m+1)!!}, \quad \;\, k=2m+1,\end{array}\right. \quad \mbox{vol}\, ({\Bbb S}^{n-1})=\frac{2\pi^{\frac{n}{2}}}{\Gamma(\frac{n}{2})},\end{eqnarray*}
  and since the space ${\Bbb R}^n$ is symmetric about the coordinate axes $\xi_1, \cdots, \xi_n$, we immediately get  \begin{eqnarray}  \label{7/15/1}
   &&\int_{{\Bbb R}^n}  \big(\sum_{j=1}^n \xi_j^2\big)^{\frac{m}{2}} \xi_k^\alpha \, e^{-\sqrt{\sum_{j=1}^n \xi_j^2}}
     d\xi =\left\{ \begin{array}{ll} \Gamma (n+m)\,vol({\Bbb S}^{n-1})\,  &\mbox{for}\;\; \alpha=0\\
    0 \, &\mbox{for} \;\; \alpha=1,\end{array} \right. \;\; n\ge 1,\nonumber\\
&&\int_{{\Bbb R}^n}   \big(\sum_{j=1}^n \xi_j^2\big)^{\frac{m-2}{2}} \xi_k\xi_l \, e^{-\sqrt{\sum_{j=1}^n \xi_j^2}}   d\xi =\left\{ \begin{array}{ll} \frac{\Gamma (n+m)\,vol({\Bbb S}^{n-1})}{n}   \, &\mbox{for}\;\; k=l\\
0 \;\; &\mbox{for}\;\; k\ne l, \end{array} \right.\;\; n\ge 2, \nonumber\\
&&\int_{{\Bbb R}^n}  \big(\sum_{j=1}^n \xi_j^2\big)^{\frac{m-4}{2}} \xi_k^2\xi_l^2\, e^{-\sqrt{\sum_{j=1}^n \xi_j^2}}   d\xi =\left\{ \begin{array}{ll} \frac{3\,\Gamma (n+m)\,vol({\Bbb S}^{n-1})}{n(n+2)}   \, &
\mbox{for} \;\, k=l
\\  \frac{\Gamma (n+m)\,vol({\Bbb S}^{n-1}) }{n(n+2)}  \, &
\mbox{for}\;\; k\ne l,\, \end{array} \right. \;\; n\ge 2,\nonumber\\
&&\int_{{\Bbb R}^n}  \big(\sum_{j=1}^n \xi_j^2\big)^{\frac{m-4}{2}} \xi_k\xi_l\xi_m\xi_p\, e^{-\sqrt{\sum_{j=1}^n \xi_j^2}}   d\xi =0\quad\; klmp\;\; \mbox{comprises}\;\, \le 1 \;\; \mbox{pairs},\;\; n\ge 3,  \\
&&\int_{{\Bbb R}^n}  \big(\sum_{j=1}^n \xi_j^2\big)^{\frac{m-6}{2}} \xi_k^6 \, e^{-\sqrt{\sum_{j=1}^n \xi_j^2}}   d\xi = \frac{15\,\Gamma (n+m)\,vol({\Bbb S}^{n-1})}{n(n+2)(n+4)},\;\; n\ge 3\nonumber\\
 &&\int_{{\Bbb R}^n}  \big(\sum_{j=1}^n \xi_j^2\big)^{\frac{m-6}{2}} \xi_k^4 \xi_l^2 \, e^{-\sqrt{\sum_{j=1}^n \xi_j^2}}   d\xi = \frac{3\,\Gamma (n+m)\,vol({\Bbb S}^{n-1})}{n(n+2)(n+4)},\quad \; k\ne l,\;\; n\ge 3,\nonumber\\
 &&\int_{{\Bbb R}^n}  \big(\sum_{j=1}^n \xi_j^2\big)^{\frac{m-6}{2}} \xi_k^2 \xi_l^2 \xi_m^2\, e^{-\sqrt{\sum_{j=1}^n \xi_j^2}}   d\xi = \frac{\Gamma (n+m)\,vol({\Bbb S}^{n-1})}{n(n+2)(n+4)},\;\quad k\ne l\ne m,\;\; n\ge 3,\nonumber\\
  &&\int_{{\Bbb R}^n}  \big(\sum_{j=1}^n \xi_j^2\big)^{\frac{m-6}{2}} \xi_k\xi_l\xi_m\xi_p\xi_q\xi_v\, e^{-\sqrt{\sum_{j=1}^n \xi_j^2}}   d\xi =0 \quad \; klmpqv\,\; \mbox{comprises}\, \le 2 \,\; \mbox{pairs}, \, n\ge 3.   \nonumber\end{eqnarray}
     It follows from (\ref{6-602}) and (\ref{6-603})   that
  \begin{eqnarray*}
  {\mathcal{K}}_{V_{-2}} (t,x,x)&=&\left(\frac{1}{2\pi}\right)^{n} \frac{t^{1-n}}{2}  \int_0^\infty
     \frac{(tr)^{n-1} d(tr)}{e^{tr}} \\
     &&  \quad   \times \int_{{\Bbb S}^n}\bigg[ \sum_{j=1}^n \kappa_j (x)-\sum_{j=1}^{n-1} \kappa_j(x) \sin^2 \phi_1 \cdots \sin^2 \phi_{j-1}
     \cos^2 \phi_j \\
     && \quad \; - \kappa_n(x) \sin^2 \phi_1 \cdots \sin^2 \phi_{n-1}
     \bigg]  dS  \\
     &&=t^{1-n}\left(\frac{1}{2\pi}\right)^{\frac{n}{2}} \frac{(n-1)\Gamma(n)\,\mbox{vol} ({\Bbb S}^{n-1})
     }{2n} \bigg(\sum_{j=1}^n \kappa_j (x)\bigg)\\
   \\          && := t^{1-n} a_1(n,x) \end{eqnarray*}
              and
               \begin{eqnarray*}
  &&{\mathcal{K}}_{V_{-3}} (t,x,x) =
  \left(\frac{1}{2\pi}\right)^{n} \int_{{\Bbb S}^n}dS \left\{\frac{t^{2-n}}{8} \int_0^\infty
     \frac{(tr)^{n-2} d(tr)}{e^{tr}}  \left[
          -2Q_1(x)\frac{{}_{}}{{}_{}} \right.\right.\\ &&
         \left.\left. +2 \sum_{j=1}^n {\tilde R}_{j(n+1)j (n+1)}(x) - 2 \sum_{j=1}^n \kappa_j^2(x)
   +3\bigg( \sum_{j=1}^n \kappa_j(x)\bigg)^2 \right.\right.\\
   &&  \left.\left.
   + 5\bigg(\sum_{j=1}^{n-1} \kappa_j(x) \sin^2 \phi_1 \cdots \sin^2 \phi_{j-1}
     \cos^2 \phi_j  + \kappa_n(x) \sin^2 \phi_1 \cdots \sin^2 \phi_{n-1}
      \bigg)^2  \right.\right. \;\; \end{eqnarray*}
      \begin{eqnarray*}&&  \left. \left.
      - \sum_{j=1}^{n-1} \bigg( 2{\tilde{R}}_{j(n+1)j(n+1)}(x)
 +6 \kappa_j^2(x)\bigg) \sin^2 \phi_1 \cdots \sin^2 \phi_{j-1}
     \cos^2 \phi_j \right.\right.\qquad \quad \quad\;\qquad \qquad  \\
      && \left.\left. - 6\kappa_n^2(x) \sin^2 \phi_1 \cdots \sin^2 \phi_{n-1}\right]  + \frac{t^{2-n}}{8} \int_0^\infty \frac{(tr)^{n-1} d(tr)}{e^{tr}} \bigg( \sum_{j=1}^n \kappa_j(x)\right.  \end{eqnarray*}
        \begin{eqnarray*}
            && \left.
     -\sum_{j=1}^{n-1} \kappa_j(x)    \sin^2 \phi_1 \cdots \sin^2 \phi_{j-1} \cos^2 \phi_j
   -\kappa_n(x) \sin^2 \phi_1 \cdots \sin^2 \phi_{n-1}
   \bigg)^2\right\}  \qquad \quad\qquad\quad \\
   &=&t^{2-n}\left(\frac{1}{2\pi}\right)^{n}  \frac{ \Gamma(n-1)\,\mbox{vol}({\Bbb S}^{n-1})}{8}  \left[ -8\sum_{1\le j<k\le n} \kappa_j\kappa_k +\frac{2(n-1)}{n}  \sum_{j=1}^n {\tilde R}_{j(n+1)j(n+1)} (x)
  \right.\\  &&\left.
  +\,\frac{-2n^2 -8n -4}{n(n+2)}\sum_{j=1}^n \kappa_j^2(x) + \frac{n^3 +2n^2 +3n +8}{n(n+2)} \bigg(\sum_{j=1}^n \kappa_j(x)\bigg)^2 \right]t^{2-n}\\
    &=&t^{2-n}\left(\frac{1}{2\pi}\right)^{n}  \frac{ \Gamma(n-1) \,\mbox{vol} ({\Bbb S}^{n-1})}{8}  \left[ \frac{n^3 -2n^2 -5n +8}{n(n+2)}
    \bigg(\sum_{j=1}^n \kappa_j(x)\bigg)^2 \right.\\  &&\left.
    +\frac{2(n-1)}{n} \sum_{j=1}^n {\tilde R}_{j(n+1)j(n+1)} (x)
   +\,\frac{2n^2 -4}{n(n+2)}\sum_{j=1}^n \kappa_j^2(x) \right]\\
     &:=&t^{2-n} \tilde{a}_2(n,x). \end{eqnarray*}

 Furthermore, \begin{eqnarray} \label{82'} && \bigg(\frac{1}{2\pi}\bigg)^{n}\int_{{\Bbb R}^n} tp_{-2}^B(x,\xi)\, e^{-t\sqrt{\sum_{j=1}^n \xi_j^2}}d\xi  = t^{3-n} \bigg(\frac{1}{2\pi}\bigg)^{n} \Gamma(n-2)\, \mbox{vol}({\Bbb S}^{n-1})\qquad \quad\qquad\;\end{eqnarray}\begin{eqnarray*}
 &&  \times \bigg\{
 \frac{-n^2-4n+5}{8n(n+2)} \big(\sum_{j=1}^n \kappa_j(x)\big)\big(\sum_{j=1}^n {\tilde R}_{j(n+1)j(n+1)}\big) +\frac{-n^2 +6n-7}{8n(n+2)(n+4)} \big(\sum_{j=1}^n \kappa_j(x)\big)^3 \\
  &&+\bigg( \frac{1}{4}+\frac{-12n^2-18n}{8n (n+2) (n+4)}  \bigg) \sum_{j=1}^n \kappa_j^3(x) - \frac{n^3 +8n^2 -7n -2}{8n(n+2)(n+4)} \bigg(\sum_{j=1}^n \kappa_j(x)\bigg) \bigg(\sum_{j=1}^n \kappa_j^2(x)\bigg) \\
  && + \bigg(\frac{1}{2}+ \frac{-5n-1}{4n(n+2)} \bigg)\sum_{j=1}^n \kappa_j (x)\, {\tilde R}_{j(n+1)j(n+1)} + \frac{n-1}{8n} \sum_{j=1}^n {\tilde R}_{j(n+1)j(n+1), (n+1)} \bigg\},\end{eqnarray*}
\begin{eqnarray} \label{83'}  \\
&& \bigg(\frac{1}{2\pi}\bigg)^{n} \int_{{\Bbb R}^n} \frac{t^3}{6} \big[ (p_0^B (x,\xi))^3 +3 p_0^2(x,\xi)\, p_0^B(x,\xi)+ 3p_0(x,\xi) (p_0^B(x,\xi))^2\big]e^{-t\sqrt{\sum_{j=1}^n\xi_j^2}} \,d\xi \qquad \quad\qquad\nonumber\\ && \qquad \qquad =
 \bigg(\frac{1}{2\pi}\bigg)^{n} \int_{{\Bbb R}^n} \frac{t^3}{6} \big(p_0^B (x,\xi)\big)^3 e^{-t\sqrt{\sum_{j=1}^n\xi_j^2}} \,d\xi\qquad \qquad \nonumber\end{eqnarray}\begin{eqnarray*} &&=
  \frac{t^{3-n}}{48} \bigg(\frac{1}{2\pi}\bigg)^{n}\Gamma(n)\,\mbox{vol}({\Bbb S}^{n-1}) \bigg[\frac{n^3 +3n^2 -7n -13}{n(n+2)(n+4)} \big(\sum_{j=1}^n \kappa_j(x)\big)^3  \\ && \quad + \frac{6n+18}{n(n+2)(n+4)} \big(\sum_{j=1}^n \kappa_j(x)\big) \big(\sum_{j=1}^n \kappa_j^2(x) \big)
  -\frac{8}{n(n+2)(n+4)}\sum_{j=1}^n \kappa_j^3(x)\bigg]\end{eqnarray*} and
 \begin{eqnarray} \label{84'} &&\bigg(\frac{1}{2\pi}\bigg)^{n}\int_{{\Bbb R}^n} t^2p_0^B (x,\xi) \big[ p_{-1}(x,\xi) + p_{-1}^B(x,\xi)\big] e^{-t\sqrt{\sum_{j=1}^n\xi_j^2}} d\xi\qquad \quad \quad \qquad\quad\quad \;\,\end{eqnarray} \begin{eqnarray*} &=& \bigg(\frac{1}{2\pi}\bigg)^{n}\frac{\Gamma(n-1)\,\mbox{vol}(S^{n-1}) t^{3-n}}{8}  \bigg[\frac{n^2+9n+24}{3n(n+2)(n+4)} \big(\sum_{j=1}^n \kappa_j(x)\big)\big(\sum_{j,k=1}^n R_{jkjk}(x)\big)
 \end{eqnarray*} \begin{eqnarray*}&& +\frac{n^2+n-1}{n(n+2)} \big(\sum_{j=1}^n\kappa_j(x)\big) \big(\sum_{j=1}^n
 {\tilde R}_{j(n+1)j(n+1)}(x)\big) + \frac{-n^3 -6n^2 -3 n +15}{2n(n+2)(n+4)} \big(\sum_{j=1}^n \kappa_j(x)\big)^3 \\
   && + \frac{n^3 +3n^2-7n -7}{n(n+2)(n+4)} \big(\sum_{j=1}^n \kappa_j (x)\big)\big(\sum_{j=1}^n  \kappa_j^2(x)\big)   +\frac{6n+4}{n(n+2)(n+4)}
   \big(\sum_{j=1}^n \kappa_j^3(x)\big) \\&& +\frac{2}{n(n+2)} \sum_{j=1}^n \kappa_j(x) \, {\tilde R}_{j(n+1)j(n+1)} (x)  -\frac{2n+16}{3n(n+2)(n+4)} \sum_{j,k=1}^n \kappa_j(x)\, {R}_{jkjk}(x) \\
   && -\frac{1}{n} \sum_{j=1}^n {\tilde R}_{j(n+1)j(n+1)}(x) +\frac{1}{2n} \big(\sum_{j=1}^n \kappa_j(x)\big)^2-\frac{1}{n} \sum_{j=1}^n \kappa_j^2(x)\bigg].\end{eqnarray*}
  Note that for any $x\in \partial \Omega$  \begin{eqnarray} \label {0-1-0} R_{jkjk} (x)= {\tilde {R}}_{jkjk} (x)  +\kappa_j(x)\kappa_k(x),\quad \; 1\le j,k\le n.\end{eqnarray}
 Since  \begin{eqnarray*} &&{\tilde{R}}_\Omega (x) = \sum_{j=1}^{n+1} {\tilde{R}}_{jj} (x) =\sum_{j=1}^{n+1} \big(\sum_{k=1}^n {\tilde{R}}_{jkjk}
 (x) +{\tilde {R}}_{j(n+1)j(n+1)}\big)\\ &&
\quad \quad\quad  =2\sum_{j=1}^n {\tilde{R}}_{j(n+1)j(n+1)}(x) +\sum_{j,k=1}^n {\tilde{R}}_{jkjk},\end{eqnarray*}
 we obtain \begin{eqnarray} \label{0-2-0} \sum_{j=1}^n {\tilde{R}}_{j(n+1)j(n+1)}(x)
= \frac{1}{2} \left( {\tilde{R}}_\Omega (x)-R_{\partial \Omega} +\sum_{1\le j\ne k\le n}^n
\kappa_j(x)\kappa_k (x)\right).\end{eqnarray}
Inserting (\ref{82'})---(\ref{84'}) into (\ref{6.-18}) and then using (\ref{0-2-0})
     we obtain that
 \begin{eqnarray}\label{7/15/3} &&{\mathcal{K}}_{V_{-4}} (t,x,x) =t^{3-n}  \bigg(\frac{1}{2\pi}\bigg)^{n} \Gamma(n-2)\,\mbox{vol} (S^{n-1}) \bigg[\frac{n^3-2n^2 -7n +7}{8n(n+2)}
   \qquad \quad  \end{eqnarray}
      \begin{eqnarray*} &&  \times  \big(\sum_{j=1}^n \kappa_j(x)\big)\big(\sum_{j=1}^n{\tilde R}_{j(n+1)j(n+1)} (x)\big)
   + \frac{n^5 -3n^4 -26n^3 +47n^2 +124 n -158}{48n(n+2)(n+4)} \big(\sum_{j=1}^n \kappa_j(x)\big)^3\\
   &&  + \frac{3n^3+7n^2 -9n-16}{12n(n+2)(n+4)}  \sum_{j=1}^n \kappa_j^3(x)
   +\frac{n^4+n^3 -16n^2 -3n+22}{8n(n+2)(n+4)}  \big(\sum_{j=1}^n \kappa_j(x)\big) \big(\sum_{j=1}^n \kappa_j^2(x)\big)\\&&
   + \frac{2 n^2 -3}{4n(n+2)} \sum_{j=1}^n \kappa_j(x)\, {\tilde R}_{j(n+1)j(n+1)} (x)+ \frac{n-1}{8n} \sum_{j=1}^n {\tilde R}_{j(n+1)j(n+1), (n+1)}(x)\\&&+ \frac{n^3 +7n^2 +6n -48}{24n(n+2)(n+4)} \big(\sum_{j=1}^n  \kappa_j (x)\big)\big(\sum_{j,k=1}^n   R_{jkjk} (x)\big)
   + \frac{-n^2 -6n +16}{12n(n+2)(n+4)} \sum_{j,k=1}^n \kappa_j(x)\,  R_{jkjk} (x)
\\ &&    -  \frac{n-2}{8n} \sum_{j=1}^n {\tilde R}_{j(n+1)j(n+1)}(x) +\frac{n-2}{16n} \big(\sum_{j=1}^n \kappa_j(x)\big)^2-\frac{n-2}{8n}\sum_{j=1}^n \kappa_j^2(x)\bigg]\end{eqnarray*}
      \begin{eqnarray*} && = t^{3-n} \bigg(\frac{1}{2\pi}\bigg)^{n}\frac{\Gamma(n-2)\,\mbox{vol} (S^{n-1}) }{8n}  \bigg[
   \frac{n^3 -2n^2 -7n+7}{2(n+2)} {\tilde R}(x)\, \big(\sum_{j=1}^n \kappa_j(x)\big) \\ && \quad \, + \frac{-3n^4 -4 n^3 +59n^2 +75 n -180}{ 6(n+2)(n+4)} \big(\sum_{j=1}^n \kappa_j (x) \big)R_{\partial \Omega} \\ && \quad \,+
   \frac{n^5 -20 n^3 +2n^2 +61 n -74}{6(n+2)(n+4)} \big(\sum_{j=1}^n \kappa_j(x)\big)^3 \\ && \quad \,
    + \frac{n^4 +8n^3 +15n^2 +3n -32}{2(n+2) (n+4)} \big(\sum_{j=1}^n \kappa_j(x)\big)\big(\sum_{j=1}^n \kappa_j^2 (x)\big)\\ && \quad \,
    +\frac{-6n^3-34n^2 +40}{3(n+2)(n+4)} \sum_{j=1}^n \kappa_j^3(x) + \frac{4n^2-6}{n+2} \sum_{j=1}^n\kappa_j(x) \big(\sum_{k=1}^{n+1} {\tilde R}_{jkjk}(x)\big)
    \\ && \quad \,-\frac{12n^3 +50 n^2 -6n -104}{3(n+2)(n+4)} \sum_{j=1}^n \kappa_j(x) \big(\sum_{k=1}^nR_{jkjk}(x)\big) + (n-1) \sum_{j=1}^n {\tilde{R}}_{j(n+1)j(n+1),(n+1)}(x)\\ && \quad \,-\frac{n-2}{2} {\tilde{R}}_{\Omega} (x)  +\frac{n-2}{2} R_{\partial \Omega} -\frac{n-2}{2} \sum_{j=1}^n \kappa_j^2(x)\\
  && := t^{3-n} a_3(n,x).
    \end{eqnarray*}
               Put  $a_2(n,x)=\frac{\Gamma(\frac{n-1}{2})\, R_{\partial \Omega}(x)}{24\pi^{\frac{n+1}{2}}}+{\tilde{a}}_2(n,x)$. In view of (\ref{0-1-0}) we get that
               \begin{eqnarray*} a_2(n,x) &=& \frac{R_{\partial \Omega}} {24 \pi^{\frac{n+1}{2}}} \, \frac{2^{1-n}\Gamma(n-1) vol({\Bbb S}^{n-1})}{\pi^{\frac{n-1}{2}} } \\ && +\left(\frac{1}{2\pi}\right)^{n}  \frac{ \Gamma(n-1) \,\mbox{vol} ({\Bbb S}^{n-1})}{8}  \left[ \frac{n^3 -2n^2 -5n +8}{n(n+2)}
    \bigg(\sum_{j=1}^n \kappa_j(x)\bigg)^2 \right.\\  &&\left.
    +\frac{n-1}{n} \left({\tilde{R}}_{\Omega} (x) -R_{\partial \Omega} (x) +\sum_{1\le j\ne k\le n}^n \kappa_j(x)\kappa_k (x) \right) \sum_{j=1}^n
   +\,\frac{2n^2 -4}{n(n+2)}\sum_{j=1}^n \kappa_j^2(x) \right]\end{eqnarray*} \begin{eqnarray*} &=& \frac{\Gamma(n-1) vol({\Bbb S}^{n-1}) }{8(2\pi)^n} \left[ \frac{3-n}{3n} R_{\partial \Omega} + \frac{n-1}{n} {\tilde R}_{\Omega}\right. \\ && \left. +\frac{n^3 -n^2 -4n +6}{n(n+2)} \big(\sum_{j=1}^n \kappa_j (x)\big)^2 + \frac{n^2 -n-2}{n(n+2)} \sum_{j=1}^n \kappa_j^2(x) \right].\qquad \qquad \qquad \quad \end{eqnarray*}
                     From the above arguments, (\ref{-b-d-11}), (\ref{-b-d-12}), (\ref{j2}) and (\ref{+61}), we find that
        \begin{eqnarray*} \label{j6} &&\int_{\partial \Omega}   {\mathcal{K}} (t,x,x)dS(x)  = \frac{\Gamma(\frac{n+1}{2})}{\pi^{\frac{n+1}{2}}} t^{-n}  \int_{\partial \Omega} 1 \,dS(x)
   + t^{1-n} \int_{\partial \Omega} a_1(n,x) dS(x)\\
     && \quad \quad \quad + \cdots + t^{M-1-n}\int_{\partial \Omega} a_{M-1}(n,x)dS(x)
           + {\mathcal{K}}_{V'_M} (t,x,x) \quad \;\mbox{as}\;\; t\to 0^+.\nonumber\end{eqnarray*}
   Finally,  if we choose $M=2$, $M=3$ and $M=4$ then (\ref{6.0.1}), (\ref{6.0.1-2}), (\ref{6.0.1-3}) follow, respectively.
         $\square$.

\vskip 0.28 true cm

\noindent{\bf Remark 6.2.} \ \  (i) \  Our asymptotic expansion is sharp because the following asymptotic expansion holds  for $t\to 0^+$ (see, (4.2.62) of \cite{Gr}):
\begin{eqnarray} \label{7/16-5}\quad\quad {\mathcal{K}} (t,x,x)\sim \sum_{j\in {\Bbb N},\; j-n\notin {\Bbb N}_+}
  c_{j-n} (x) t^{j-n} + \sum_{j-n\in {\Bbb N}_+}
  c_{j-n} (x) t^{j-n}\log\, t +\sum_{l\in {\Bbb N}_+}
  r_l(x) t^{l},\end{eqnarray}
  in the usual sense: the difference between ${\mathcal{K}}(t,x,x)$ and the sum of terms up to $j-n=N \in {\Bbb N}_+$, $l=N$, is $O(t^{N+1})$.
Consequently (see (4.2.64) of \cite{Gr}),
\begin{eqnarray} \label{7/16-6} &&\mbox{Tr}\, e^{t{\mathcal{N}}_g} \sim t^{j-n}\sum_{j\in {\Bbb N},\; j-n\notin {\Bbb N}_+}
  \int_{\partial \Omega} c_{j-n} (x)dS(x)\,   \\ && \qquad \qquad + t^{j-n}\log\, t\sum_{j-n\in {\Bbb N}_+}
  \int_{\partial \Omega}  c_{j-n} (x)dS(x)   +t^{l}\sum_{l\in {\Bbb N}_+}
  \int_{\partial \Omega} r_l(x)dS(x).\nonumber \end{eqnarray}

 (ii) \   By our method, we can also get the asymptotic expansion for any integer $M>4$:
\begin{eqnarray} &&   \int_{\partial \Omega}  {\mathcal{K}}(t,x,x) dS(x)= t^{-n}\frac{\Gamma(\frac{n+1}{2})}{\pi^{\frac{n+1}{2}}} \int_{\partial \Omega} 1\, dS(x) + t^{1-n} \int_{\partial \Omega} a_1(n,x)\,dS(x)\\
  && \qquad\quad \quad +t^{2-n} \int_{\partial \Omega} a_2(n,x)\, dS(x)   +\cdots+ t^{M-1-n}  \int_{\partial \Omega}
 a_{M-1}(n,x)\, dS(x) \nonumber\\
  && \qquad\quad \quad+ \left\{\begin{array}{ll} O(t^{M-n}) \quad \, \mbox{when}\;\; n>M-1,\nonumber\\ O(t\log t) \;\quad\; \mbox{when} \;\; n=M-1,\end{array}\right.\quad \;\;\, \mbox{as}\;\; t\to 0^+, \nonumber \end{eqnarray}
if we further calculate the lower-order symbol equations for the operators $-\sqrt{-\Delta_h}$ and $B$, respectively (note that ${\mathcal{N}}_g=-\sqrt{-\Delta_h} +B$).

(iii) \ The above asymptotic
expansion shows that one can ``hear'' $\int_{\partial \Omega} a_{M-1}(n,x)dS(x)$ ($M=1,2,3,\cdots$)  by ``hearing'' all of the Steklov eigenvalues.

(iv)  \ We can obtain the same result when
 we directly calculate the first three symbols of $-\sqrt{-\Delta_h}$ and $B$ instead of using ${\mathcal{K}}_V(t,x,y)$.

(v) \  By applying the Tauberian theorem (see, for example, Theorem 15.3 of p.$\,$30 of \cite{Kor}) for the first term on the right side of (\ref{6.0.1}), we immediately get Sandgren's asymptotic formula (\ref{1-4}) for the Steklov eigenvalues.

\vskip 0.38 true cm

In particular, we have the following:
\vskip 0.15 true cm

\noindent{\bf Corollary 6.3.} \  {\it Let $(\mathcal{M},g)$ be an $(n+1)$-dimensional, smooth Riemannian
manifold, and let $\Omega\subset \mathcal{M}$ be a bounded domain
with smooth boundary $\partial \Omega$. Let $h$ be the induced metric on $\partial \Omega$ by $g$. If $F(t,x,y)$ is a fundamental solution of $\frac{\partial u}{\partial t} =-\sqrt{-\Delta_h}\,u$ on $[0,+\infty)\times (\partial \Omega)$, then
\begin{eqnarray} && \int_{\partial \Omega} F(t,x,x)dS(x)=\frac{\Gamma(\frac{n+1}{2})}{\pi^{\frac{n+1}{2}}} t^{-n} \int_{\partial \Omega} 1\, dS(x) + \frac{\Gamma(\frac{n-1}{2})}{12\pi^{\frac{n+1}{2}}}t^{2-n}  \int_{\partial \Omega} R_{\partial \Omega}\, dS(x)\nonumber
 \\
 && \quad \quad \quad\; + \frac{\Gamma(\frac{n-3}{2})}{720\pi^{\frac{n+1}{2}}} t^{4-n} \int_{\partial \Omega} (10A- B+2C)\, dS(x)  + \frac{\Gamma(\frac{n-5}{2})}{4^3\pi^{\frac{n+1}{2}}}t^{6-n}D+o(t^{8-n}) \quad \;\mbox{as}\;\; t\to 0^+\nonumber\end{eqnarray}
 with \begin{eqnarray*} A= \big( \sum_{1\le i<j\le n} R_{ijij} \big)^2 =(R_{\partial \Omega})^2,  \, \quad B= \sum_{1\le j,k\le n} \big(\sum_{1\le i\le n} R_{ijik}\big)^2, \, \quad
   C= \sum_{1\le i,j,k,l\le n}  \big(R_{ijkl}\big)^2,\end{eqnarray*}
   \begin{eqnarray*} &&D = \int_{\partial \Omega} D(n,x)dS(x)= \frac{(4\pi)^{-n/2}}{7!} \int_{\partial \Omega} \big(-\frac{142}{9} R_{ijij,k} R_{mlml,k} -\frac{26}{9} R_{ijik,l}R_{mjmk,l}\qquad \qquad\qquad
   \qquad \quad\qquad \quad \qquad \qquad \\
     && \;\;-\frac{7}{9} R_{ijkm,l} R_{ijkm,l} -\frac{35}{9} R_{ijij}R_{mlml}R_{pqpq} + \frac{14}{3} R_{ijij} R_{mlmp} R_{qlqp} -\frac{14}{3} R_{ijij} R_{mlpq} R_{mlpq} \\
     && \;\;+4R_{ijik} R_{jlml} R_{kpmp} -\frac{20}{9} R_{ijik} R_{lpmp} R_{jlkm}+\frac{8}{9} R_{ijik} R_{jlmp} R_{klmp} -\frac{8}{3} R_{ijkl} R_{ijmp} R_{klmp}\big)dS(x),\end{eqnarray*}
   where $R_{ijkl}$ and $R_{ijkl,m}$ denote the curvature tensor and its covariant derivative on $(\partial \Omega, h)$,respectively.}

\vskip 0.24 true cm

\noindent {\bf Proof.}  It follows from (1.5a) and (7.2) of \cite{MS} and p.613 of \cite{Gil} that,  as $t\to 0^+$,
\begin{eqnarray*} && \int_{\partial \Omega} G(t, x,x)dS(x)= \int_{\partial \Omega} \left[\left( G_0(t,x,x) + \big(G_0\#\big((\Delta_h -\Delta_h^0)G_0\big)\big) (t,x,x) \frac{{}_{}}{} \right.\right.\\ && \,\left. \left. +
 G_0\# \big((\Delta_h -\Delta_h^0)G_0\big) \# \big((\Delta_h -\Delta_h^0) G_0\big)\right)(t,x,x)\right.\\ && \,\left.+
 G_0\# \big((\Delta_h -\Delta_h^0)G_0\big)\# \big((\Delta_h -\Delta_h^0)G_0\big) \# \big((\Delta_h -\Delta_h^0) G_0\big)(t,x,x)
 \right] dS(x) +o(t^{4-\frac{n}{2}})\nonumber\\
 && =\frac{1}{(4\pi t)^{n/2}} \left[ \int_{\partial \Omega} 1\, dS(x)  + \frac{t}{3} \int_{\partial \Omega} R_{\partial \Omega}\, dS(x)  +\frac{t^2}{180} \int_{\partial \Omega} (10 A - B + 2C)dS(x)\right. \nonumber
 \\  &&\left. \, +  D t^3 +o(t^4)\right].\nonumber\end{eqnarray*}
 Combing this and (\ref{-b.21}),  we immediately get an  asymptotic expansion:
 \begin{eqnarray*} && \int_{\partial \Omega} F(t,x,x)dS(x)= \int_0^\infty \frac{te^{-\frac{t^2}{4\mu}}}{\sqrt{4\pi \mu^3}} \left\{ \frac{1}{(4\pi \mu)^{n/2}} \left[ \int_{\partial \Omega} 1\, dS(x) + \frac{\mu}{3} \int_{\partial \Omega} R_{\partial \Omega}\, dS(x) \right.\right. \\
   && \quad \left.\left.  +\frac{\mu^2}{180} \int_{\partial \Omega} (10 A - B + 2C)dS(x) + D\mu^3 +o(\mu^4)\right]\right\} d\mu\nonumber\\
   &&  =\frac{\Gamma(\frac{n+1}{2})}{\pi^{\frac{n+1}{2}}} t^{-n} \int_{\partial \Omega} 1\, dS(x) + \frac{\Gamma(\frac{n-1}{2})}{12\pi^{\frac{n+1}{2}}}t^{2-n}  \int_{\partial \Omega} R_{\partial \Omega}\, dS(x)\nonumber
 \\
 && \quad  + \frac{\Gamma(\frac{n-3}{2})}{2880\,\pi^{\frac{n+1}{2}}} t^{4-n} \int_{\partial \Omega} (10A- B+2C)\, dS(x)
 +\frac{\Gamma(\frac{n-5}{2})}{4^3\pi^{\frac{n+1}{2}}}t^{6-n}D
 +o(t^{8-n}) \quad \;\mbox{as}\;\; t\to 0^+.\nonumber\end{eqnarray*}

\vskip 0.28 true cm

\noindent{\bf Remark 6.4.} \ \ For $n=2$, it follows from \cite{MS} that $10A-B+2C= 12 (R_{\partial \Omega})^2$. Therefore by applying the classical Gauss-Bonnet formula for the Euler characteristic $E$ of the (two-dimensional) boundary $\partial \Omega$ (i.e.,
$\int_{\partial \Omega} R_{\partial \Omega}=2\pi E$), we get \begin{eqnarray*} &&\int_{\partial \Omega} F(t,x,x)dS(x)
= \frac{\Gamma(\frac{3}{2})}{\pi^{3/2}t^{2}} \, \mbox{vol}\, (\partial \Omega) + \frac{\Gamma(\frac{1}{2})}{6\pi^{1/2}}  E + \frac{\Gamma(-\frac{1}{2})\,t^{2}}{240\pi^{3/2}}  \int_{\partial \Omega} (R_{\partial \Omega})^2  \\
&&\quad \quad \qquad +\frac{\Gamma(-\frac{3}{2})D}{4^3\pi^{3/2}}t^4+o(t^{6}) \quad \mbox{as}\;\; t\to 0^+,\nonumber\end{eqnarray*}
in particular, the Euler characteristic of boundary $\partial \Omega$ is audible by hearing all the eigenvalues of $-\sqrt{-\Delta_h}$ on $\partial \Omega$.
\vskip 1.68  true cm

\centerline {\bf  Acknowledgments}

\vskip 0.38 true cm
 I wish to express my sincere gratitude to Professor
 L. Nirenberg and Professor Fang-Hua Lin for their
  support and help.  This
research was supported by SRF for ROCS, SEM (No. 2004307D01)
   and NNSF of China (11171023/A010801).

  \vskip 1.65 true cm

\vskip 0.32 true cm

\end{document}